\documentclass[10pt,a4paper]{article}

\usepackage{subcaption}
\DeclareCaptionFormat{custom}
{\textbf{#1#2}{\small #3}}
\usepackage{amsmath,amsfonts}
\usepackage{graphicx}
\usepackage{epstopdf}
\usepackage[section]{algorithm}
\usepackage{algpseudocode}
\ifpdf
  \DeclareGraphicsExtensions{.eps,.pdf,.png,.jpg}
\else
  \DeclareGraphicsExtensions{.eps}
\fi
\usepackage{multirow}	
\usepackage{arydshln}	

\usepackage{amsthm}		
\newtheorem{theorem}{Theorem}[section]
\newtheorem{lemma}{Lemma}[section]

\newtheorem{proposition}{Proposition}[section] 
\newtheorem{corollary}{Corollary}[section] 
\newtheorem{property}{Property}[section] 
\newtheorem{definition}{Definition}[section]

\usepackage{color}
\definecolor{gray}{rgb}{0.65,0.65,0.65}
\newcommand{\gr}[1]{\textcolor{gray}{#1}}
\definecolor{rouge}{rgb}{1.0,0.0,0.0}
\definecolor{bleu}{rgb}{0.0,0.0,1.0}
\definecolor{vert}{rgb}{0.2,0.7,0.2}

\usepackage{hyperref}
\usepackage[capitalise,nameinlink]{cleveref}[0.19]
\crefname{fact}{Fact}{Facts} 
\crefname{property}{Property}{Properties} 
\crefname{appendix}{Appendix}{Appendices}
\crefformat{equation}{\textup{#2(#1)#3}}
\Crefformat{equation}{\textup{#2(#1)#3}}
\crefrangeformat{equation}{\textup{#3(#1)#4--#5(#2)#6}}
\Crefrangeformat{equation}{\textup{#3(#1)#4--#5(#2)#6}}
\crefmultiformat{equation}{\textup{#2(#1)#3}}{ and \textup{#2(#1)#3}}{, \textup{#2(#1)#3}}{, and \textup{#2(#1)#3}}
\Crefmultiformat{equation}{\textup{#2(#1)#3}}{ and \textup{#2(#1)#3}}{, \textup{#2(#1)#3}}{, and \textup{#2(#1)#3}}
\crefrangemultiformat{equation}{\textup{#3(#1)#4--#5(#2)#6}}{ and \textup{#3(#1)#4--#5(#2)#6}}{, \textup{#3(#1)#4--#5(#2)#6}}{, and \textup{#3(#1)#4--#5(#2)#6}}
\Crefrangemultiformat{equation}{\textup{#3(#1)#4--#5(#2)#6}}{ and \textup{#3(#1)#4--#5(#2)#6}}{, \textup{#3(#1)#4--#5(#2)#6}}{, and \textup{#3(#1)#4--#5(#2)#6}}

\newcommand{\comment}[1]{}

\providecommand{\card}[1]{\lvert#1\rvert}
\providecommand{\abs}[1]{\lvert#1\rvert}
\providecommand{\set}[1]{\{#1\}}
\providecommand{\minus}[1]{\backslash\set{#1}}

\begin{document}

\title{Deriving differential approximation results for $k$-CSPs from combinatorial designs\footnote{Parts of this work have been published in conference proceedings \cite{CT18,CT18-E}.}}

\author{Jean-Fran\c{c}ois Culus\\
	\texttt{jean-francois.culus@st-cyr.terre-net.defense.gouv.fr}\\
	CREC, Saint-Cyr, France\\
	MEMIAD, Université des Antilles, France\\
	\\
	Sophie Toulouse\footnote{Corresponding author}\\ 
	\texttt{sophie.toulouse@lipn.univ-paris13.fr}\\
	Université Sorbonne Paris Nord, CNRS, LIPN, France
}

\maketitle

\hrule
\begin{abstract} 
{\noindent}Inapproximability results for $\mathsf{Max\,k\,CSP\!-\!q}$ have been traditionally established using balanced $t$-wise independent distributions, which are closely related to orthogonal arrays, a famous family of combinatorial designs. 
In this work, we investigate the role of these combinatorial structures in the context of the differential approximability of $\mathsf{k\,CSP\!-\!q}$, providing new structural insights and approximation bounds.

We first establish a direct connection between the average differential ratio on $\mathsf{k\,CSP\!-\!q}$ instances and orthogonal arrays. This allows us to derive the new differential approximability bounds of 
	$1/q^k$ for $(k +1)$-partite instances, 
	$\Omega(1/n^{\lfloor k/2\rfloor})$ for Boolean instances, 
	$\Omega(1/n)$ when $k =2$, 
	and $\Omega(1/n^{k -\lceil\log_{\Theta(q)}k\rceil})$ when $k, q\geq 3$.
We then introduce families of array pairs, called {\em alphabet reduction pairs of arrays}, that are still related to balanced $k$-wise independence. Using these pairs of arrays, we establish a reduction from $\mathsf{k\,CSP\!-\!q}$ to $\mathsf{k\,CSP\!-\!k}$ (where $q >k$), with an expansion factor of $1/(q -k/2)^k$ on the differential approximation guarantee. Combining this with a 1998 result by Yuri Nesterov, we conclude that $\mathsf{2\,CSP\!-\!q}$ is approximable within a differential factor of $0.429/(q -1)^2$.
Finally, using similar Boolean array pairs, {\em called cover pairs of arrays}, we prove that every Hamming ball of radius $k$ provides a $\Omega(1/n^k)$-approximation of the instance diameter. 

Thus, our work highlights the relevance of combinatorial designs for establishing structural differential approximation guarantees for CSPs.


\end{abstract}
\hrule

\section{Introduction}
Many combinatorial optimization problems can be formulated as {\em Constraint Satisfaction Problems}, CSPs, over a finite domain. In a CSP over an alphabet $\Sigma$, we have a set $\set{x_1 ,\ldots, x_n}$ of $\Sigma$-valued variables and a set $\set{C_1, \ldots, C_m}$ of constraints. 
Each constraint $C_i$ applies some predicate $P_i:\Sigma^{k_i}\rightarrow\set{0, 1}$ to a subset of variables. The goal is to find an assignment that optimizes the number of satisfied constraints.

CSPs include many well-known problems. For example, the {\em Maximum Satisfiability Problem} ($\mathsf{Max\,Sat}$) is the Boolean CSP where the objective is to satisfy as many disjunctive clauses $(\ell_{i_1} \vee\ldots\vee \ell_{i_{k_i}})$ as possible, where a literal $\ell_j$ represents either the Boolean variable $x_j$ or its negation $\bar x_j$. 
Another important CSP is $\mathsf{Lin\!-\!q}$, for which constraints are equations of the form
\[(\alpha_{i, 1} x_{i_1} +\ldots+ \alpha_{i, k_i} x_{i_{k_i}}\equiv\alpha_{i, 0}\bmod{q})\] 
where the coefficients $\alpha_{i,j}$ are elements of $\mathbb{Z}_q:=\mathbb{Z}/q\mathbb{Z}$.

\begin{figure}\footnotesize
$I =\min_{x =(x_1, x_2, x_3, x_4)\in\Sigma_3^4} v(I, x)$ where:
\begin{align}\nonumber
v(I, x)	&=
	(x_1 +x_2\equiv 1\bmod{3}) 
	+2.3\times(x_1 =x_3)
	+1.4\times(x_2 =1 \wedge x_4 =2)  
	+7\times x_4
\end{align}
For this instance, we have:\\
\begin{tabular}{l}
$\bullet$	$k_1 =2$, $w_1 =1$,		$J_1 =(1, 2)$, 
		$C_1 =P_1(x_{J_1}) =P_1(x_1, x_2) =1$ if $x_1 +x_2\equiv 1\bmod{3}$ and 0 otherwise,\\
$\bullet$	$k_2 =2$, $w_2 =2.3$,	$J_2 =(1, 3)$, 
		$C_2 =P_2(x_{J_2}) =P_2(x_1, x_3) =1$ if $x_1 =x_3$ and 0 otherwise,	\\
$\bullet$	$k_3 =2$, $w_3 =1.4$,	$J_3 =(2, 4)$,
		$C_3 =P_3(x_{J_3}) =P_3(x_2, x_4) =1$ if $(x_2, x_4) =(1, 2)$ and 0 otherwise,		\\
$\bullet$	$k_4 =1$, $w_4 =7$,		$J_4 =(4)$,
		$C_4 =P_4(x_{J_4}) =P_4(x_4) =x_4$.
\end{tabular}
\caption{Illustration of an instance $I$ of $\mathsf{CSP\!-\!q}$ and the notations used to describe such an instance, where $q =3$, $n =4$ and $m =4$.}\label{ex-CSP}
\end{figure}

In the most general case, for each $i\in[m] :=\set{1, 2 ,\ldots, m}$, the constraint $C_i$ is associated with a positive weight $w_i$, and the functions $P_i$ can take real values (see {\em e.g.} \cite{R08,A10} for the latter generalization). The goal is then to optimize an objective function of the form 
$$\begin{array}{rl}
v(I, x)	&=\sum_{i =1}^m w_i P_i(x_{J_i}) =\sum_{i =1}^m w_i P_i(x_{i_1} ,\ldots, x_{i_{k_i}})
\end{array}$$
over $\Sigma^n$, where for all $i\in[m]$, $J_i =(i_1 ,\ldots, i_{k_i})$ is a subsequence of $(1 ,\ldots, n)$, $k_i \leq n$, $P_i:\Sigma^{k_i}\rightarrow \mathbb{R}$, and $w_i >0$. 

We define the alphabet $\Sigma_q$ of size $q$ as the set $\Sigma_q :=\set{0, 1, 2 ,\ldots, q -1}$. When working with algebraic operations such as addition and multiplication, we consider working modulo $q$ and identify $\Sigma_q$ with the ring $\mathbb{Z}_q$, which forms a field when $q$ is prime. 
We denote by $\mathsf{CSP\!-\!q}$ the case where the variables take values in $\Sigma_q$.  \Cref{ex-CSP} shows an example when $q =3$ and $n =m =4$.

For a family $\mathcal{F}$ of functions, we denote by $\mathsf{CSP(\mathcal{F})}$ the CSP where the functions $P_i$ are elements of $\mathcal{F}$ (for all $i\in[m]$). For example, $\mathsf{Lin\!-\!2}$ corresponds to 
$$\mathsf{CSP(\set{XNOR^k,XOR^k\,|\,k\in\mathbb{N}\backslash\set{0}})}$$
where, for a positive integer $k$, $XNOR^k$ and $XOR^k$ refer to the $k$-ary Boolean predicates that are true for entries with an even and an odd number of non-zero coordinates, respectively.

In this paper, we focus on $k$-CSPs, a subclass of CSPs where each constraint depends on at most $k$ variables. 
For a specific CSP $\mathsf{\Pi}$ ({\em e.g.}, $\mathsf{\Pi} =\mathsf{Lin\!-\!q}$), we denote by $\mathsf{k\,\Pi}$ ({\em e.g.}, $\mathsf{k\,Lin\!-\!q}$) and $\mathsf{Ek\,\Pi}$ ({\em e.g.}, $\mathsf{Ek\,Lin\!-\!q}$) its restriction to instances where each constraint depends on {\em at most} and {\em exactly} $k$ variables, respectively. Its restrictions to instances where the goal is to maximize or to minimize are denoted by $\mathsf{Max\,\Pi}$ and $\mathsf{Min\,\Pi}$, respectively.

When constraints can take both positive and negative values, maximizing and minimizing a CSP become equivalent problems. In fact, flipping the sign of each constraint $C_i$ is equivalent to inverting the optimization objective, turning a maximization problem into a minimization one, and vice versa. Thus, $\mathsf{Max\,CSP\!-\!q}$ and $\mathsf{Min\,CSP\!-\!q}$ are equivalent optimization problems. 

\subsection{Approximation measures}

Even when $q =k =2$, $\mathsf{Max\,2\,Sat}$ and $\mathsf{Min\,2\,Sat}$ are $\mathbf{NP\!-\!hard}$ \cite{GJS76, KKM94}. 
Therefore, an important issue in optimization CSPs is to characterize their computational complexity by studying their {\em approximability}. 
Given an optimization CSP $\mathsf{\Pi}$, we denote its instance set by $\mathcal{I}_\mathsf{\Pi}$. For an instance $I\in\mathcal{I}_\mathsf{\Pi}$, we denote the best and worst solution values on $I$ by $\mathrm{opt}(I)$ and $\mathrm{wor}(I)$, respectively. An approximation measure quantifies how close the value $v(I, x)$ of an approximate solution $x$ is to the optimal value.
 %
A widely used approximation measure is the {\em standard approximation ratio}, which directly compares $v(I, x)$ to $\mathrm{opt}(I)$. More formally, the standard ratio of a solution $x$ on an instance $I$ is defined as\footnote{The standard ratio is also commonly defined as the inverse ratio $\max\left\{v(I, x)/\mathrm{opt}(I), \mathrm{opt}(I)/v(I,x)\right\}$.}:
$$\min\left\{v(I, x)/\mathrm{opt}(I), \mathrm{opt}(I)/v(I, x)\right\}$$

A solution $x$ is said to be {\em $\rho$-standard approximate} on $I$ for some $\rho\in(0, 1]$ if this ratio is at least $\rho$. A {\em $\rho$-standard approximation algorithm} is a polynomial-time algorithm $\mathcal{A}$ that, given any instance $I$ of $\mathsf{\Pi}$, returns a solution that is at least $\rho(I)$-standard approximate. $\mathsf{\Pi}$ is {\em approximable within a standard factor of $\rho$} if such an algorithm exists.

A natural way to analyze the average solution value of an instance $I$ of $\mathsf{CSP\!-\!q}$ is through the expected value $\mathbb{E}_X[v(I, X)]$ of a random solution, where $X =(X_1 ,\ldots, X_n)$ is a vector of pairwise independent random variables, each uniformly distributed over $\Sigma_q$. 
H{\aa}stad and Venkatesh in \cite{HV04} introduced an approximation measure, which we refer to as the {\em gain approximation measure}. This measure is based on the {\em optimal advantage over a random assignment}, defined as the difference $\abs{\mathrm{opt}(I) -\mathbb{E}_X[v(I, X)]}$, which quantifies how much an optimal solution outperforms a purely random assignment.
The gain ratio of a solution $x$ on $I$ is the ratio:
$$\frac{v(I, x) -\mathbb{E}_X[v(I, X)]}{\mathrm{opt}(I) -\mathbb{E}_X[v(I, X)]}$$

This measure was motivated by the fact that for many CSPs, for all constant $\varepsilon >0$, finding  solutions with value at least $\mathbb{E}_X[v(I, X)] +\varepsilon\times\sum_{i =1}^m w_i$ in almost satisfiable instances is $\mathbf{NP\!-\!hard}$. For example, $\mathsf{E3\,Lin\!-\!2}$ is such a CSP \cite{H97}.
 %
The {\em differential approximation measure} is based on the distance to a worst solution value instead of the mean solution value. Namely, the differential ratio of $x$ on $I$ is the ratio:
$$\frac{v(I, x) -\mathrm{wor}(I)}{\mathrm{opt}(I) -\mathrm{wor}(I)}$$

The distance $|\mathrm{opt}(I) -\mathrm{wor}(I)|$ between the extreme values $\mathrm{opt}(I)$ and $\mathrm{wor}(I)$ is known as the {\em diameter} of $I$. 
The differential ratio gained prominence in approximation theory due to its stability under affine transformations of the objective function, meaning that rescaling or shifting the objective function does not affect the differential ratio of the solutions \cite{ABMV77,AAP80,BR95,DP96}. 

The notions of {\em $\rho$-differential} and {\em $\rho$-gain} 
{\em approximate solutions}, {\em approximation algorithms}, and {\em approximable problems} are defined analogously to their counterparts for the standard approximation measure.
These approximation measures follow a hierarchical relationship. 
For a maximization instance $I$ where $v(I,.)$ is non-negative, every solution $x$ satisfies:
$$\begin{array}{rll}															 \displaystyle
\frac{v(I, x)}{\mathrm{opt}(I)}													&\displaystyle
	\geq \frac{v(I, x) -\mathrm{wor}(I)}{\mathrm{opt}(I) -\mathrm{wor}(I)}		&\displaystyle
	\geq \frac{v(I, x) -\mathbb{E}_X[v(I, X)]}{\mathrm{opt}(I) -\mathbb{E}_X[v(I, X)]}
\end{array}$$

In particular, for all positive integers $q$ and $k$, if $\mathsf{k\,CSP\!-\!q}$ is $\rho$-gain approximable, then it is $\rho$-differential approximable and, if it is, $\mathsf{Max\,k\,CSP\!-\!q}$ is $\rho$-standard approximable {\em on instances with non-negative solution values}.

\subsection{Differential approximability of CSPs}\label{sec-dapx}

In this paper we investigate three questions concerning the differential approximability of $\mathsf{k\,CSP\!-\!q}$, about which, unlike the standard approximation, only a few facts are known. 
 %
The {\em Conjunctive Constraint Satisfaction Problem}, $\mathsf{CCSP}$ for short, is the Boolean CSP where  constraints are conjunctive clauses. For all constants $\varepsilon>0$, the restriction of $\mathsf{Max\,CCSP}$ to unweighted instances is $\mathbf{NP\!-\!hard}$ to approximate within standard approximation ratio $1/m^{1 -\varepsilon}$, where we recall that $m$ represents the number of constraints in the CSP instance. This is due to the standard inapproximability bound of \cite{H99,Z07} for the {\em Maximum Independent Set problem}, which extends by reduction to $\mathsf{Max\,CCSP}$ \cite{BR95}. 
$\mathsf{Max\,Sat}$ is inapproximable within any constant differential factor assuming $\mathbf{P}\neq\mathbf{NP}$, as Escoffier and Paschos argue in \cite{EP05}. However, the same authors also noted that the conditional expectation technique \cite{J74} provides $1/m$-differential approximate solutions on unweighted instances of $\mathsf{Sat}$ \cite{EP05} (see \cref{sec-average} for more details). 
For $\mathsf{Lin\!-\!2}$, H{\aa}stad and Venkatesh show that combining this technique with exhaustive search allows to approximate the optimal gain over a random assignment within a factor of $\Omega(1/m)$. $\mathsf{Lin\!-\!2}$ is therefore in particular $\Omega(1/m)$-differentially approximable. 

The differential approximability bound of $\Omega(1/m)$ for $\mathsf{Lin\!-\!2}$ extends to $\mathsf{k\,CSP\!-\!q}$ for all constant integers $k$ and $q$, using a binary encoding of the variables and the discrete Fourier transform \cite{CT18-E}.
When $q =k =2$, $\mathsf{2\,CSP\!-\!2}$ admits a $(2 -\pi/2)$-differential approximation algorithm (where $2 -\pi/2 >0.429$), which combines the semidefinite programming-based algorithm of Goemans and Williamson \cite{GW95} with derandomization techniques such as in \cite{MR99}.
This result is due to Nesterov, who established this approximation guarantee for {\em Unconstrained Binary Quadratic Programming} in \cite{N98}. The approximability bound of $2 -\pi/2$ extends by reduction to $\mathsf{3\,CSP\!-\!2}$, up to a multiplicative factor of $1/2$ on the approximation guarantee \cite{CT18}. 

Whether $\mathsf{k\,CSP\!-\!q}$ admits a constant differential approximation factor remains an open question for  $q\geq 3$ or $k\geq 4$. However, standard inapproximability bounds are known. 
The {\em primary hypergraph} of a CSP instance $I$ contains a vertex $j$ for each variable $x_j$ of $I$, and a hyperedge $e_i =(i_1 ,\ldots, i_{k_i})$ for each constraint $C_i =P_i(x_{i_1} ,\ldots, x_{i_{k_i}})$ of $I$. Assimilating $I$ to its primary hypergraph, a {\em strong coloring} of $I$ is a partition $V_1\sqcup\ldots\sqcup V_\nu$ of $[n]$ such that the support $J_i =(i_1 ,\ldots, i_{k_i})$ of any constraint intersects each color set $V_c$ in at most one index. We say that $I$ is {\em $\nu$-partite} if such a partition of size $\nu$ exists. The smallest integer $\nu$ for which $I$ is $\nu$-partite is called the {\em strong chromatic number} of $I$ (see {\em e.g.} \cite{AK13}).
Let $q\geq 2$ and $k\geq 3$ be two integers. In \cite{C13}, Chan establishes that the restriction of $\mathsf{Max\,k\,CSP\!-\!q}$ to $k$-partite instances with non-negative solution values is $\mathbf{NP\!-\!hard}$ to standardly approximate within constant ratio better than
	$(q -1) k/q^{k -1}$ if $q$ is a prime power, 
	$O(k/q^{k -1})$ if $k\geq q$, 
	and $O((q -1)k/q^{k -1})$ otherwise. 
With respect to the differential approximation measure, these inapproximability bounds hold in $k$-partite instances of $\mathsf{k\,CSP\!-\!q}$. 
 %
Moreover, assuming $\mathbf{P}\neq\mathbf{NP}$, for all constant $\varepsilon >0$, the 6-gadget reducing $\mathsf{E3\,Lin\!-\!2}$ to $\mathsf{E2\,Lin\!-\!2}$ in \cite{H97} implies a differential inapproximability bound of $7/8 +\varepsilon$ for bipartite instances of $\mathsf{E2\,Lin\!-\!2}$ (see \cref{sec-ap-bib} for more details).

\subsection{Approximability of CSPs and balanced $t$-wise independence}

Many inapproximability bounds for $k$-CSPs, including the ones of \cite{C13}, rely on {\em balanced $t$-wise independent distributions} or {\em balanced $t$-wise independent subsets} (notably see \cite{AM08, AK13}). 
 %
Let $q\geq 1$, $t\geq 1$, and $\nu\geq t$ be three integers. A probability distribution $\mu$ on $\Sigma_q^\nu$ is said to be {\em balanced $t$-wise independent} if, for any $t$ coordinates $(Y_{c_1}, \dots, Y_{c_t})$ of a vector $Y$ from the probability space $(\Sigma_q^\nu, \mu)$, every $t$-tuple $(v_1, \dots, v_t)$ of values appears with probability $1/q^t$.
 %
By extension (see {\em e.g.} \cite{C13}), a subset $\mathcal{U}$ of $\Sigma_q^\nu$ is said to be {\em balanced $t$-wise independent} if, for each sequence $J = (c_1, \dots, c_t)$ of $t$ indices from $[\nu]$ and each $v \in \Sigma_q^t$, $\mathcal{U}$ contains exactly $|\mathcal{U}|/q^t$ vectors $u$ such that $(u_{c_1}, \dots, u_{c_t}) = (v_1, \dots, v_t)$.

For example, the predicate $ZeroSum^{\nu, q}$ evaluates to $1$ for $(y_1, \dots, y_\nu)\in\Sigma_q^\nu$ such that $y_1 + \dots + y_\nu\equiv 0\bmod{q}$. Fixing any $\nu -1$ variables to any $\nu -1$ values $v_1, \dots, v_{\nu -1}$ uniquely determines the remaining variable to be $(-v_1 -\dots - v_{\nu -1})\bmod{q}$ to satisfy the equation. The $q^{\nu -1}$ accepting entries of $ZeroSum^{\nu, q}$ thus form a balanced $(\nu -1)$-wise independent subset of $\Sigma_q^\nu$.
Furthermore, consider the probability distribution on $\Sigma_q^\nu$ that assigns the probability $1/q^{\nu-1}$ to the accepting entries of $ZeroSum^{\nu, q}$, and $0$ to all other vectors of $\Sigma_q^\nu$. This distribution is clearly $(\nu -1)$-wise independent. 

For a function $P$ on $\Sigma_q^k$, we denote by $r_P$ the average value of $P$ over $\Sigma_q^k$, i.e.: 
$$\begin{array}{rl}
	r_P	&:= \sum_{y\in\Sigma_q^k} P(y)/q^k
\end{array}$$
Furthermore, for $v\in\Sigma_q^k$, we define the {\em shifted function $P_v$} that assigns to each $y\in\Sigma_q^k$ the value of $P$ taken at $y +v$. Formally:
$$\begin{array}{rll}
	P_v(y_1 ,\ldots, y_k)	&= P\left((y_1 +v_1) \bmod{q} ,\ldots, (y_k +v_k) \bmod{q}\right),	&y_1 ,\ldots, y_k \in\Sigma_q
\end{array}$$

For instance, in \cref{ex-CSP}, $P_1$ is the function $ZeroSum^{2, 3}_v$ with $v =(2, 0)$. Moreover, let $AllZeros^{k, q}$ denote the predicate on $\Sigma_q^k$ that accepts the single entry $(0 ,\ldots, 0)$. Then, in the same example, we have $P_3= AllZeros^{2, 3}_v$ with $v =(2, 1)$. 
 %
The inapproximability bounds from \cite{C13} mentioned in \cref{sec-dapx} actually follow from the Theorem below\footnote{We present a simplified version of the result of \cite{C13}, which applies to CSPs over finite abelian groups.}:
\begin{theorem}[\cite{C13}]\label{thm-inapx-kCSP}
Let $k\geq 3$ and $q\geq 2$ be two fixed integers, and $P$ be a predicate on $\mathbb{Z}_q^k$ such that the set $P^{-1}(1)$ of its accepting entries forms a balanced pairwise independent subgroup of $\mathbb{Z}_q^k$. Then $\mathsf{Max\,CSP(\set{P_v\,|\,v\in\mathbb{Z}_q^k})}$ is $\mathbf{NP\!-\!hard}$ to standardly approximate within constant ratio greater than $r_P$, even in $k$-partite instances. 
\end{theorem}

As an example, consider the predicate $ZeroSum^{k, q}$ where $q\geq 2$ and $k\geq 3$. The set of its accepting entries forms a subgroup of $\mathbb{Z}_q^k$. Since this subgroup is balanced $(k -1)$-wise independent with $k -1\geq 2$, it is also balanced pairwise independent. By applying \Cref{thm-inapx-kCSP}, we obtain a differential inapproximability bound of $1/q$ for $k$-partite instances of $\mathsf{k\,CSP\!-\!q}$, assuming $\mathbf{P}\neq\mathbf{NP}$. 
In particular, this result implies an approximability upper bound of 1/2 for $\mathsf{3\,CSP\!-\!2}$. 
For other values of $(k, q)$, Chan establishes the bounds of $O(k/q^{k -1})$ or $O(k/q^{k -2})$, as mentioned earlier. These bounds are obtained using more sophisticated predicates, mostly derived from infinite families of linear codes. They also involve a standard approximation-preserving reduction to CSPs over smaller alphabets of size a prime power (see \cref{sec-reduc-intro} for more details).  

With a slight abuse of terminology, we say that {\em a function} is {\em balanced $t$-wise independent} 
if its mean value remains unchanged when any $t$ of its variables are fixed at arbitrary values. Formally, such a function $P$ of $k\geq t$ variables, each with domain $\Sigma_q$, must satisfy:
\begin{align}\label{eq-I_q^t}
\textstyle\sum_{y\in\Sigma_q^k: y_J =v} P(y)/q^{k -t} 	&=r_P,
	&J =(j_1 ,\ldots, j_t)\in[k]^t,\,j_1 <\ldots<j_t,\ v\in\Sigma_q^t	
\end{align}

From now on, we denote the set of such functions by $\mathcal{I}_q^t$. 
Balanced $t$-wise independent functions on $\Sigma_q^k$ are a natural extension of balanced $t$-wise independent distributions on $\Sigma_q^k$, which satisfy \cref{eq-I_q^t}. Conversely, let $P$ be a function on $\Sigma_q^k$ with minimum value $P_*$. Then $P$ satisfies \cref{eq-I_q^t} {\em if and only if} the function
\begin{align}\nonumber
	y	&\displaystyle\mapsto\tilde{P}(y) 
			:=\frac{P(y) -P_*}{\sum_{u\in\Sigma_q^k} (P(u) -P_*)},	&\forall y\in\Sigma_q^k	
\end{align}
defines a balanced $t$-wise independent distribution on $\Sigma_q^k$ (by construction, $\tilde{P}$ takes values in $[0, 1]$, has a mean value of $1/q^k$, and satisfies \cref{eq-I_q^t} {\em if and only if} $P$ does).
For example, if $P=ZeroSum^{k, q}$, then its normalized version $\tilde{P}$ corresponds to $ZeroSum^{k, q}/q^{k -1}$.

Balanced $t$-wise independent functions are also a natural extension of balanced $t$-wise subsets. Consider here that the predicates consisting in accepting exactly the vectors of such subsets of $\Sigma_q^k$ necessarily satisfy \cref{eq-I_q^t}. For instance, for all integers $k, q\geq 2$, $ZeroSum^{k, q}$ is such a predicate with $t =k -1$; hence, for all integers $k, q\geq 2$, $ZeroSum^{k, q}\in\mathcal{I}_q^{k -1}$. 

\subsection{Outline}

We identify new connections between balanced $t$-wise independence and optimization CSPs over $q$-ary alphabets, which allow to establish new positive and conditional differential approximation results for $k$-CSPs. 
In the inapproximability results of \cite{C13}, balanced $t$-wise independence restricts the functions used to express the constraints of the CSP. In this article, balanced $t$-wise independence essentially concerns distributions on the solution set of the CSP instance. More specifically, we manipulate arrays or pairs of arrays to model multisets of solutions, where each row represents  a solution of the CSP instance. Balanced $t$-wise independence precisely constrains the frequency of the rows in the arrays. 
In particular, our results involve a famous family of combinatorial designs called {\em Orthogonal Arrays} (OAs for short), which can be thought of as rational balanced $t$-wise independent measures (see \cref{sec-def_OA+DS} for a complete definition of OAs). 
 %
The paper is organized as follows:

$\bullet$ In \cref{sec-average}, we study the differential ratio reached at the average solution value of a $\mathsf{k\,CSP\!-\!q}$ instance. We show a connection between this ratio and OAs or related designs, involving the strong chromatic number of the instance (\cref{thm-E-OA}). 
Orthogonal arrays and linear codes from the literature then allow to deduce that the average differential ratio is $\Omega(1)$ on instances with a bounded strong chromatic number, 
$\Omega(1/n^{k/2})$ when $q =2$, 
and $\Omega(1/n^{k -\lceil\log_{\Theta(q)} k\rceil})$ in all other cases.

$\bullet$ In \cref{sec-reduc}, we introduce array pairs with entries from $\Sigma_q$, called {\em alphabet reduction pairs of arrays}, which can be viewed as a constrained decomposition of balanced $k$-wise independent functions on $\Sigma_q^q$. We show that these pairs of arrays allow finding solutions of CSPs over an alphabet of size $q$ by solving CSPs over an alphabet of a smaller size $p$, provided that each constraint depends on at most $p$ variables (\cref{thm-reduc-CD}). 
By constructing such pairs (\cref{thm-gamma_qpk}), we show that whenever $\mathsf{k\,CSP\!-\!k}$ is approximable within some differential factor $\rho$, for any $q >k$, $\mathsf{k\,CSP\!-\!q}$ is approximable within differential factor $\rho/(q -k/2)^k$.
Thus, in particular, it follows from \cite{N98} that for all constant integers $q \geq 2$, $\mathsf{2\,CSP\!-\!q}$ is differentially approximable within a constant factor.

$\bullet$ In \cref{sec-vois}, we consider pairs of arrays similar to those of \cref{sec-reduc}, but with Boolean coefficients, called {\em cover pairs of arrays}. We relate these array pairs to the highest differential ratio of a solution over an arbitrary Hamming ball of fixed radius. 
Interpreting the pairs of arrays constructed in the preceding section as cover pairs of arrays, we show that, for $\mathsf{k\,CSP\!-\!q}$, every Hamming ball of radius $k$ contains a pair of solutions whose difference in value is a fraction $\Omega(1/n^k)$ of the instance diameter.	
In fact, the obtained pairs of arrays allow to express the value of any solution as a linear combination of solution values over an arbitrary Hamming ball of radius $k$ (\cref{thm-vois-id}).

As is customary, we discuss the results obtained and the prospects they offer in a concluding section. Technical arguments and side issues are compiled in separate appendices, \cref{app-pv} and \cref{app-supp}, respectively.

\subsection{Conventions and notations used in the rest of the paper}

For an instance $I$ of a CSP, unless otherwise specified, $n$, $m$ and $\nu$ will always refer to its number of variables, its number of constraints, and its strong chromatic number.

For a positive integer $q$, arithmetic operations on elements of $\Sigma_q$ are always performed modulo $q$, and arithmetic operations over $\Sigma_q^\nu$ are interpreted componentwise modulo $q$.

\medskip{\bf{\em Uniform shifts.}} 
For a symbol $a\in\Sigma_q$, we denote by $\mathbf{a} =(a ,\ldots, a)$ the vector whose coordinates are all equal to $a$ (the dimension depends on the context). In particular, for a function $P$ of variables with domain $\Sigma_q$, $P_{\mathbf{a}}$ refers to its translation by the vector $\mathbf{a}$.
For example, if we consider the third constraint of \cref{ex-CSP}, then for $(x_2, x_4)\in\Sigma_3^2$, we have: 
$$\begin{array}{rlll}	
	{P_3}_{\mathbf{1}}(x_2, x_4) 
				&=AllZeros^{2, 3}_{(2, 1) +\mathbf{1}}(x_2, x_4) 
				&=AllZeros^{2, 3}_{0, 2}(x_2, x_4)
				&=(x_2 =0 \wedge x_4 =1)\\
	{P_3}_{\mathbf{2}}(x_2, x_4) 
				&=AllZeros^{2, 3}_{(2, 1) +\mathbf{2}}(x_2, x_4) 
				&=AllZeros^{2, 3}_{1, 0}(x_2, x_4)
				&=(x_2 =2 \wedge x_4 =0)\\
\end{array}$$

{\bf{\em Function families $\mathcal{E}_q$ and $\mathcal{O}_q$.}} 
 %
Given a positive integer $k$, $XOR^k$ is notable in that, for any two $k$-dimensional Boolean vectors $y$ and $\bar y$, we have either $XOR^k(y) =XOR^k(\bar y)$, or $XOR^k(y) +XOR^k(\bar y) =1$, depending on whether $k$ is even or odd. This follows from the fact that the number of non-zero coordinates in $\bar y$ has the same parity as the number of non-zero coordinates in $y$ {\em if and only if} $k$ is even.

$\mathcal{E}_q$ and $\mathcal{O}_q$ generalize such Boolean predicates to $q$-ary alphabets.
Functions in $\mathcal{E}_q$ are stable under a uniform shift of all their variables, while functions in $\mathcal{O}_q$ have the property that their mean value over any $q$ successive shifts $y, y +\mathbf{1}, \ldots, y +\mathbf{q -1}$ coincides with their overall mean value.
Formally, given a positive integer $k$, a function $P:\Sigma_q^k\rightarrow\mathbb{R}$ belongs to $\mathcal{E}_q$ and $\mathcal{O}_q$ if it satisfies the following relations \cref{eq-E_q} and \cref{eq-O_q}, respectively:
\begin{align}
			\label{eq-E_q}	P_{\mathbf{a}}(y) &=P(y),	&&y\in\Sigma^k_q,\ a\in\Sigma_q\\[-2pt]
\textstyle	\label{eq-O_q}	\sum_{a =0}^{q -1}P_{\mathbf{a}}(y)/q 	&=r_P,		&&y\in\Sigma^k_q
\end{align}

For example, let $AllEqual^{k, q}$ refer to the predicate on $\Sigma_q^k$ that accepts the entries $(y_1 ,\ldots, y_k)$ with $y_1 =\ldots= y_k$. For instance, constraint $C_2$ in \cref{ex-CSP} consists of $AllEqual^{2, 2}(x_1, x_3)$. It is obvious that $AllEqual^{k, q}$ belongs to $\mathcal{E}_q$. 
 %
Now consider the equation $(y_1 +\ldots+ y_k \equiv 0 \bmod{q})$ over $\Sigma_q^k$. 
For any $(y_1 ,\ldots, y_k)\in\Sigma_q^k$ and any $a\in\Sigma_q$, we have:
$$\begin{array}{rll}
	(y_1 +a) +\ldots+ (y_k +a) 	&=(y_1 +\ldots+ y_k) +k\times a
\end{array}$$ 
We deduce that $ZeroSum^{k, q}\in\mathcal{E}_q$ if $k$ is a multiple of $q$ (in which case $ka\equiv 0\bmod{q}$), and that $ZeroSum^{k, q}\in\mathcal{O}_q$ if $k$ and $q$ are mutually prime (in which case $ka \equiv -(y_1 +\ldots+ y_k)\bmod{q}$ for a single $a\in\Sigma_q$).

Let $P\in\mathcal{E}_q$. From \cref{eq-E_q}, we deduce that fixing any variable of $P$ to an arbitrary value does not affect its mean value, which remains $r_P$. Therefore, every function in $\mathcal{E}_q$ is an element of $\mathcal{I}^1_q$.
For some insight into the function families $\mathcal{E}_q$ and $\mathcal{O}_q$ and the corresponding CSPs, we invite the reader to refer to \cref{sec-func}.

\medskip{\bf{\em Arrays.}}
 %
Let $\nu$ and $q \geq 2$ be two positive integers. 
An array with $\nu$ columns on the symbol set $\Sigma_q$ is a multisubset of $\Sigma_q^\nu$. 
An array is {\em simple} if no word $u\in\Sigma_q^\nu$ occurs more than once as a row in it.
Given an array $M$, we denote by $M_r$ the row indexed by $r$, and by $M^c$ the column indexed by $c$.

For an $R\times\nu$ array $M$ on $\Sigma_q$, we denote by $\mu^M$ the empirical frequency of the words of $\Sigma_q^\nu$ in $M$, i.e.:
\begin{align}\nonumber
	\mu^M(u)	&\textstyle:= \left\vert\left\{r\in[R]\,|\,M_r =u\right\}\right\vert/R,	&u\in\Sigma_q^\nu
\end{align}
$\mu^M$ defines a probability distribution on $\Sigma_q^\nu$. 
In addition, we consider the function $\mu^M_E$ which smooths $\mu^M$ by averaging over uniform shifts. Namely, $\mu^M_E$ assigns to each $u\in\Sigma_q^\nu$ a fraction $1/q$ of the total frequency in $M$ of words of the form $u +\mathbf{a}$:
\begin{align}\nonumber
	\mu^M_E(u)	&\textstyle:=\sum_{a =0}^{q -1} \mu^M(u +\mathbf{a})/q, 	&u\in\Sigma_q^\nu
\end{align}

By construction, $\mu^M_E$ also defines a probability distribution on $\Sigma_q^\nu$, which belongs to $\mathcal{E}_q$. 

\medskip{\bf{\em Neighborhoods.}}  
Given two positive integers $q$ and $\kappa$, the {\em Hamming distance} between two vectors $x, y\in\Sigma_q^\nu$, denoted by $d_H(x, y)$, is the number of coordinates on which $x$ and $y$ differ, i.e.:
\begin{align}\nonumber
	d_H(x, y)			&=\card{\set{j\in[\nu]\,:\,x_j \neq y_j}},		&x, y\in\Sigma_q^\nu 
\end{align}

For $d\in\set{0, 1 ,\ldots, n}$ and $x\in\Sigma_q^n$, the {\em Hamming ball of radius $d$ centered at $x$}, denoted by $B^d(x)$, is the set of vectors $y\in\Sigma_q^n$ that are at Hamming distance at most $d$ from $x$, i.e.:
\begin{align}\nonumber
	B^d(x)		&=\set{y\in\Sigma_q^n\,:\,d_H(x, y)\leq d},		&x\in\Sigma_q^n 
\end{align}

In particular, $x\in B^d(x)$. Furthermore, for an integer $a\in\Sigma_q$, we denote by $B^d_\mathbf{a}$ the function that associates with each $x\in\Sigma_q^n$ the shift by $\mathbf{a}$ of $B^d(x)$. Equivalently, $B^d_\mathbf{a}$ associates with each $x\in\Sigma_q^n$ the Hamming ball of radius $d$ centered at $x +\mathbf{a}$. Namely:
\begin{align}\nonumber
	B^d_\mathbf{a}(x)	&=\set{y +\mathbf{a}\,|\,y\in B^d(x)} 
						 =B^d(x +\mathbf{a}),		&x\in\Sigma_q^n,\ a\in\Sigma_q 
\end{align}

Finally, we denote by $\tilde{B}^d$ the function that associates with each $x\in\Sigma_q^n$ the union of the uniform shifts of $B^d(x)$, i.e.:
\begin{align}\nonumber
	\tilde{B}^d(x)	&=B^d(x) \cup B^d_\mathbf{1}(x) \cup\ldots\cup B^d_\mathbf{q -1}(x)
					 =\cup_{a =0}^{q -1} B^d(x +\mathbf{a}),		&x\in\Sigma_q^n 
\end{align}

\subsection{Obtained approximation bounds}%

We summarize in \cref{tab-dapx_lit,tab-dapx_new} the resulting knowledge of the differential approximability of $\mathsf{k\,CSP\!-\!q}$, its restrictions $\mathsf{k\,CSP(\mathcal{I}_q^t)}$ and $\mathsf{k\,CSP(\mathcal{E}_q)}$, and the restriction $\mathsf{CSP(\mathcal{O}_q)}$ of $\mathsf{CSP\!-\!q}$.

\begin{table}\footnotesize{
\begin{tabular}{p{1.7cm}|p{10.5cm}}
Restriction	&Approximation bound\\\hline
\multirow{2}{*}{$\mathsf{CSP(\mathcal{O}_q)}$}	&$1/q$ (trivial)\\
[2pt]	&$\neg\,1/q +\varepsilon$, even in tripartite instances of 
		$\mathsf{E3\,CSP(\mathcal{O}_q\cap\mathcal{I}_q^2)}$ \cite{C13}\\
[5pt]
\multirow{2}{*}{\begin{minipage}{1.7cm}
$\mathsf{2\,CSP\!-\!2}$,\\\hspace*{0.25cm}$\mathsf{3\,CSP(\mathcal{E}_2)}$
\end{minipage}}		&$2 -\pi/2\ (>0.429)$ \cite{N98,CT18}\\
[2pt]	&$\neg\,7/8 +\varepsilon$ in bipartite instances of $\mathsf{Lin\!-\!2}$,
		due to the gadget of \cite{H97} from $\mathsf{E3\,Lin\!-\!2}^*$\\
[5pt]
\multirow{2}{*}{$\mathsf{3\,CSP\!-\!2}$}	&$1 -\pi/4\ (>0.214)$ using \cite{N98}, by reduction to $\mathsf{2\,CSP\!-\!2}$ \cite{CT18}\\
[2pt]	&$\neg\,1/2 +\varepsilon$, even in 3-partite instances of $\mathsf{E3\,Lin\!-\!2}$ \cite{C13}\\
[5pt]
$\mathsf{k\,CSP\!-\!q}$		&$\Omega(1/m)$ using \cite{HV04}, by reduction to $\mathsf{Lin\!-\!2}$ \cite{CT18-E}\\
[5pt]
\begin{minipage}{1.7cm}
$\mathsf{k\,CSP\!-\!q}$\\\hspace*{0.25cm}with $k\geq 3$
\end{minipage}	&\begin{minipage}{10.5cm}
	$\neg\,O(k/q^{k -1}) +\varepsilon$ if $k\geq q$, $\neg\,(q -1)k/q^{k -1} +\varepsilon$ if $q$ is a prime power, and\\\hspace*{0.25cm}$\neg\,O((q -1)k/q^{k -1}) +\varepsilon$ otherwise, even in $k$-partite instances of $\mathsf{Ek\,CSP(\mathcal{I}_q^2)}$ \cite{C13}
				\end{minipage}\\
\end{tabular}\\
$^*$ see \cref{sec-ap-bib} for more details.	
\caption{\label{tab-dapx_lit}
Differential approximability bounds that are already known for $\mathsf{k\,CSP\!-\!q}$ and $\mathsf{CSP(\mathcal{O}_q)}$, where $k\geq 2$ and $q\geq 2$. Inapproximability bounds hold for all constant $\varepsilon >0$, provided that $\mathbf{P}\neq\mathbf{NP}$.} 
}\end{table}

\begin{table}\footnotesize{
{\em Lower bounds for the average differential ratio (\cref{sec-average}):}\\
\begin{tabular}{l|l|l}
Lower bound	&Restriction	&Conditions	on $\nu$, $q$, $k$, and $t$	\\\hline
$1/q$			&$\mathsf{CSP(\mathcal{O}_q)}$		&	\\
\hline&&\\[-8pt]
$1/q^{k -1}$	&$\mathsf{k\,CSP(\mathcal{E}_q)}$	&$q$ or $k$ is odd and $\nu\leq k +1$	\\
\hline&&\\[-8pt]
$1/q^{\min\{\nu -t, k\}}$	&$\mathsf{k\,CSP(\mathcal{I}_q^t)}$	&$\nu\leq k +t +1$	\\
\hline&&\\[-8pt]
\multirow{4}{*}{$1/q^k$}
	&\multirow{2}{*}{$\mathsf{k\,CSP\!-\!q}$}
			 &$\nu\leq k +1$, or $q$ prime power, $q >k$, and $\nu\leq q +1$,		\\[1pt]
			&&or $q$ power of 2, $q >3 =k$, and $\nu\leq q +2$							\\
\cline{2-3}&&\\[-8pt]
	&\multirow{2}{*}{$\mathsf{k\,CSP(\mathcal{I}_q^t)}$}
			 &$q$ prime power, $q >k$ and $\nu\leq q +1 +t$,		\\[1pt]
			&&or $q$ power of 2, $q >3 =k$, and $\nu\leq q +2 +t$	\\
\hline&&\\[-8pt]
\multirow{2}{*}{$1/(2q -2)^k$}	&\multirow{2}{*}{$\mathsf{k\,CSP\!-\!q}$}
		 &$\lceil\log_2 q\rceil >k$ and $\nu\leq 2^{\lceil\log_2 q\rceil} +1$,		\\[1pt]
		&&or $\lceil\log_2 q\rceil >3 =k$ and $\nu\leq 2^{\lceil\log_2 q\rceil} +2$	\\
\hline&&\\[-8pt]
$\Omega(1/\nu^{\lfloor k/2\rfloor})$				&$\mathsf{k\,CSP\!-\!2}$	&\\
\hline&&\\[-8pt]
$\Omega(1/\nu^{k -\lceil\log_{\Theta(q)} k\rceil})$	&$\mathsf{k\,CSP\!-\!q}$	&$q\geq 3$
\\\multicolumn{3}{c}{}\\[-0pt]
\end{tabular}

{\em Differential approximability bounds derived from \cite{N98} by reduction to $\mathsf{2\,CSP\!-\!2}$ (\cref{sec-reduc}):}\\
\begin{tabular}{l|l}
Restriction							&Approximation bound\\\hline
$\mathsf{2\,CSP\!-\!q}$ for $q\geq 3$							&$(2 -\pi/2)/(q -1)^2$\\[1pt]
$\mathsf{2\,CSP(\mathcal{E}_q)}$ for $q\in\set{3, 4, 5, 7, 8}$	&$(2 -\pi/2)/q$
\\\multicolumn{2}{c}{}\\[-0pt]
\end{tabular}

{\em Approximability bounds related to Hamming balls of fixed radius (\cref{sec-vois}):}\\
\begin{tabular}{l|l}
Restriction		&Approximation guarantee\\\hline
$\mathsf{CSP(\mathcal{O}_q)}$
	&For all solutions $x$, the best differential ratio over $\set{x +\mathbf{a}\,|\,a\in\Sigma_q}$ is at least $1/q$.\\
[8pt]
\begin{minipage}{0.11\textwidth}
$\mathsf{2\,CSP\!-\!2}$,\\\hspace*{0.25cm}$\mathsf{3\,CSP(\mathcal{E}_2)}$
\end{minipage}								&\begin{minipage}{0.8\textwidth}
Local optima w.r.t. $\tilde{B}^1$ are $1/O(\nu)$-differential approximate 
	and for all solutions $x$, the best differential ratio over $\tilde{B}^1(x)$ is $1/O(\nu\times n)$.
\end{minipage}\\
[11pt]
$\mathsf{Ek\,CSP(\mathcal{I}_q^{k -1})}$	&\begin{minipage}{0.8\textwidth}
The differential ratio of local optima w.r.t. $B^1$ is at least the average differential ratio
and for all solutions $x$, the highest differential ratio over $B^1(x)$ is $\Omega(1/n)$ times this ratio.
\end{minipage}\\
[11pt]
$\mathsf{k\,CSP\!-\!q}$						&\begin{minipage}{0.8\textwidth}
For all solutions $x$, the maximum distance between two solution values over $B^k(x)$ is a fraction $\Omega(1/n^k)$ of the instance diameter.
\end{minipage}
\end{tabular}
\caption{\label{tab-dapx_new}
	New differential approximability bounds for $\mathsf{k\,CSP\!-\!q}$, $\mathsf{CSP(\mathcal{O}_q)}$, $\mathsf{k\,CSP(\mathcal{E}_q)}$ and $\mathsf{k\,CSP(\mathcal{I}_q^t)}$, where $k\geq 2$, $q\geq 2$, and $t\in[k -1]$.} 
}\end{table}

\section{Differential approximation quality of a random assignment}\label{sec-average}

Solutions with value at least $\mathbb{E}_X[v(I, X)]$ are computationally easy to find, using the {\em conditional expectation technique} \cite{J74}. The method, when applied to an instance $I$ of $\mathsf{Max\,CSP\!-\!q}$, consists in associating a (new) random variable $X_j$ to each variable $x_j$ of $I$, and then iteratively fixing variables $x_j, j =1 ,\ldots, n$ to a symbol $a\in\Sigma_q$ that maximizes the conditional expectation:
	$$\mathbb{E}_X \left[v(I, X)\mid \,(X_1, X_2 ,\ldots, X_{j -1}, X_j) =(x_1, x_2 ,\ldots, x_{j -1}, a) \right].$$ 

By proceeding in this way, provided that the variables $X_j$, $j\in [n]$ are independently distributed, we  obtain a solution $x$ with value:
$$\begin{array}{rlcl}
	v(I, x)		&= 		\mathbb{E}_X \left[v(I, X)\,|\,X = x \right]\\
				&\geq	\mathbb{E}_X \left[(X_1 ,\ldots, X_{n -1}) =(x_1 ,\ldots, x_{n -1}) \right]	
				&\geq	\ldots	\geq	&\mathbb{E}_X[v(I, X)]
\end{array}$$

In particular, the method returns a solution with value at least the average solution value when the variables $X_j$, $j\in [n]$ are uniformly distributed.  
This naturally raises two questions: {\em is it possible to compute solutions in polynomial that outperform $\mathbb{E}_X[v(I, X)]$}, and {\em what is the gain of $\mathbb{E}_X[v(I, X)]$ over the worst solution value?} 

Approximating the optimum advantage over a random assignment provides a way to address the former question. 
We here address the second question. We specifically seek lower bounds for the average differential ratio on instances of $\mathsf{k\,CSP\!-\!q}$ and its restrictions $\mathsf{k\,CSP(\mathcal{E}_q)}$ and $\mathsf{k\,CSP(\mathcal{I}_q^t)}$. Such lower bounds provide in particular an estimate of the differential approximation guarantee offered by the conditional expectation technique.

\begin{figure}[t]
\begin{center}
\includegraphics[scale=0.8]{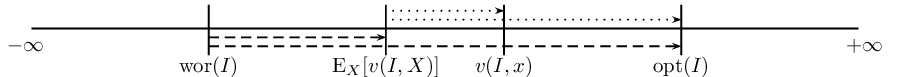}
\end{center}
\caption{\label{fig-mes}Quantities involved in the gain ratio achieved by a given solution $x$ (in dotted lines) and the average differential ratio (in dashed lines), on an instance $I$ where the goal is to maximize.}
\end{figure}

Note that the two questions are complementary, and the second has potential to offer insights into the first. In particular, one might think that the further $\mathbb{E}_X[v(I, X)]$ is from $\mathrm{wor}(I)$, the harder it is to get away from it. \Cref{fig-mes} shows the quantities involved in the gain approximation measure and the average differential ratio.

\subsection{Previous related works and preliminary remarks}

We discuss three restrictions for which some lower bound on the average differential ratio is either already established or obvious.

\smallskip
$\bullet$ {\bf{\em Unweighted instances.}}
In \cite{EP05}, Escoffier and Paschos analyze the differential ratio of solutions returned by the conditional expectation technique on unweighted instances of $\mathsf{Sat}$. They observed that on such an instance $I$ on which the goal is to maximize, provided that $\mathrm{opt}(I)\neq\mathrm{wor}(I)$, we have:
$$\begin{array}{rll}
	\lceil\mathbb{E}_X[v(I, X)]\rceil
		&\geq \mathrm{wor}(I) +1
		&\geq \mathrm{wor}(I) +\left(\mathrm{opt}(I) -\mathrm{wor}(I)\right)/m
\end{array}$$

Thus, the conditional expectation technique provides a $(1/m)$-differential approximate solution on $I$. We observe that for all integers $q\geq 2$, the same argument holds for instances of $\mathsf{CSP\!-\!q}$ with integer solution values and a polynomially bounded diameter. However, the argument does no extend to instances with arbitrary solution values.

\smallskip
$\bullet$ {\bf{\em Submodular pseudo-Boolean optimization.}}
A real-valued function that depends on Boolean variables is called {\em pseudo-Boolean}.
Given a positive integer $n$, a pseudo-Boolean function $P:\set{0, 1}^n\rightarrow\mathbb{R}$ is {\em submodular} {\em if and only if} it satisfies:
$$\begin{array}{rll}
P(y) +P(z)	&\geq P(y_1 \vee z_1 ,\ldots, y_n \vee z_n) +P(y_1 \wedge z_1 ,\ldots, y_n \wedge z_n),
				&y, z\in \set{0, 1}^n
\end{array}$$

Feige et al. proved in \cite{FMV11} that for any submodular pseudo-Boolean function $P$, if $x^*$ is a maximizer of $P$, then:
\begin{align}\label{eq-E-SM}
	\mathbb{E}_X[P(X)]	&\geq\textstyle P(x^*)/4 +P(\bar x^*)/4 +P(\mathbf{0})/4 +P(\mathbf{1})/4
\end{align}

Since a conical combination of submodular pseudo-Boolean functions remains submodular, inequality \cref{eq-E-SM} holds if $P$ is the objective function $v(I, .)$ of an instance of $\mathsf{Max\,CSP\!-\!2}$ where the functions $P_i$ are all submodular. 
Considering that none of the solution values $v(I, \bar x^*)$, $v(I, \mathbf{0})$, $v(I, \mathbf{1})$ can be less than $\mathrm{wor}(I)$, we deduce that the average differential ratio on such instances is at least $1/4$. 
For submodular instances of $\mathsf{Max\,CSP(\mathcal{E}_2)}$, since $v(I, \bar x^*) =\mathrm{opt}(I)$, this ratio is even lower bounded by $1/2$.

The {\em Maximum Directed Cut} problem, $\mathsf{Max\,Di\,Cut}$ is the restriction of $\mathsf{Max\,2\,CCSP}$ to clauses of the form $(x_{i_1} \wedge \bar x_{i_2})$, while the Boolean {\em Not-All-Equal Satisfiability} problem, $\mathsf{NAE\,Sat}$ is the restriction of $\mathsf{CSP\!-\!2}$ to constraints of the form $\neg(\ell_{i_1} =\ldots= \ell_{i_{k_i}})$. 
Submodular CSPs notably include $\mathsf{Max\,Di\,Cut}$, the restriction 
--- known as the {\em satisfiability problem with no mixed clause} --- 
of $\mathsf{Max\,Sat}$ to constraints of the form $(x_{i_1} \vee\ldots\vee x_{i_{k_i}})$ or $(\bar x_{i_1} \vee\ldots\vee \bar x_{i_{k_i}})$, and the restriction 
--- known as the {\em monotone not-all-equal satisfiability problem} ---  
of $\mathsf{Max\,NAE\,Sat}$ to constraints of the form $\neg (x_{i_1} =\ldots= x_{i_{k_i}})$. 
	
\smallskip
$\bullet$ {\bf{\em Restriction to functions of $\mathcal{O}_q$.}}
\label{sec-E-O_q} 
Consider a function $P\in\mathcal{O}_q$ over $\Sigma_q^n$. 
By definition of $\mathcal{O}_q$, the expected value of $P$ satisfies: 
\begin{align}\label{eq-E-O_q}
	\mathbb{E}_X[P(X)] =r_P	&=\left(P(x) +P(x +\mathbf{1}) +\ldots+ P(x +\mathbf{q -1})\right)/q,	&x\in\Sigma_q^n
\end{align}

Taking \cref{eq-E-O_q} at a maximizer $x^*$ and at a minimizer $x_*$ of $P$, we deduce that  $\mathbb{E}_X[P(X)]$ is a $1/q$-differential approximation of both $P(x^*)$ and $P(x_*)$. Since $\mathcal{O}_q$ is stable by linear combinations, this holds in particular when $P$ is the objective function of an instance of $\mathsf{CSP(\mathcal{O}_q)}$. Thus, for $\mathsf{CSP(\mathcal{O}_q)}$, the average differential ratio is at least $1/q$.

\medskip
Imposing uniform weights, or requiring that the functions expressing the constraints are submodular or belong to $\mathcal{O}_q$, is quite restrictive. Let us shift our focus to an instance $I$ of $\mathsf{Max\,Ek\,CSP\!-\!q}$, where $k, q\geq 2$. On $I$, the average solution value can be expressed as:
$$\begin{array}{rll}
\mathbb{E}_X[v(I, X)]
	&=\sum_{i =1}^m w_i \times\sum_{v\in\Sigma_q^k} P_i(v)/q^k
	&=\sum_{i =1}^m w_i r_{P_i}
\end{array}$$
(For example, when $I$ is an instance of $\mathsf{Lin\!-\!2}$, this value is equal to $\sum_{i =1}^m w_i\times 1/2$.)
For an optimal solution $x^*$ of $I$, we observe:
\begin{align}
\textstyle\sum_{i =1}^m w_i r_{P_i}		&\textstyle= \sum_{i =1}^m w_i \times
	\left(P_i(x^*_{J_i}) +\sum_{v\in\Sigma_q^k\minus{x^*_{J_i}}} P_i(v)\right)/q^k
\ \ \Leftrightarrow\nonumber\\[-2pt]
\mathbb{E}_X[v(I, X)]					&\textstyle= 
	v(I, x^*)/q^k +\sum_{i =1}^m w_i\times\sum_{v\in\Sigma_q^k\minus{x^*_{J_i}}} P_i(v)/q^k					\label{eq-E-q}
\end{align}

Thus, on $I$, provided that $w_i P_i$ is non-negative for all $i\in[m]$, the ratio of $\mathbb{E}_X[v(I, X)]$ to $\mathrm{opt}(I)$ is at least $1/q^k$. This ratio is even larger if $I$ is an instance of $\mathsf{Ek\,CSP(\mathcal{E}_q)}$ or $\mathsf{Ek\,CSP(\mathcal{I}_q^t)}$.
 %
First suppose that there is some $t\in[k -1]$ for which, for all $i\in[m]$, $P_i$ is balanced $t$-wise independent. For a constraint $C_i =P_i(x_{J_i})$ of $I$, we denote by $L_i =(i_1 ,\ldots, i_t)$ the sequence of the $t$ first elements of $J_i$, and by $R_i =(i_{t +1} ,\ldots, i_k)$ the remaining subsequence of $J_i$. By substituting in \cref{eq-I_q^t}, we obtain for the average solution value on $I$ the following expression:
\begin{align}
\textstyle\sum_{i =1}^m w_i r_{P_i}						&\textstyle=	
	\sum_{i =1}^m w_i \times\sum_{v\in\Sigma_q^{k -t}} P_i(x^*_{L_i}, v)/q^{k -t}
\qquad\text{by \cref{eq-I_q^t}}		\nonumber\\[-2pt]	&\textstyle=
	\sum_{i =1}^m w_i \left(P_i(x^*_{L_i}, x^*_{R_i}) 
			+\sum_{v\in\Sigma_q^{k -t}\minus{x^*_{R_i}}} P_i(x^*_{L_i}, v)\right)/q^{k -t}	
\ \ \Leftrightarrow\nonumber\\[-2pt]
\mathbb{E}_X[v(I, X)]									&\textstyle= v(I, x^*)/q^{k -t} 
	+\sum_{i =1}^m w_i\times\sum_{v\in\Sigma_q^{k -t}\minus{x^*_{R_i}}} P_i(x^*_{L_i}, v)/q^{k -t}																								\label{eq-E-I_q^t}
\end{align}

Since $\mathcal{E}_q\subseteq\mathcal{I}_q^1$, equality \cref{eq-E-I_q^t} holds in particular for $t =1$ if $I$ is an instance of $\mathsf{Ek\,CSP(\mathcal{E}_q)}$.
 %
It follows from \cref{eq-E-q,eq-E-I_q^t} that, for any instance of $\mathsf{Max\,Ek\,CSP\!-\!q}$, $\mathsf{Max\,Ek\,CSP(\mathcal{I}_q^t)}$, and $\mathsf{Max\,Ek\,CSP(\mathcal{E}_q)}$ where the constraints and their weights are all non-negative, the average standard ratio is at least $1/q^k$, $1/q^{k -t}$, and $1/q^{k -1}$, respectively. 
However, similar deductions cannot be made for the average differential ratio. 
Specifically, in the most general case, we cannot claim that the quantities
$$\textstyle\sum_{i =1}^m w_i \sum_{v\in\Sigma_q^k\minus{x^*_{J_i}}} P_i(v)/q^k\text{ and }
\sum_{i =1}^m w_i \sum_{v\in\Sigma_q^{k -t}\minus{x^*_{R_i}}} P_i(x^*_{L_i}, v)/q^{k -t}$$
appearing in the right-hand sides of the equalities \cref{eq-E-q,eq-E-I_q^t} are greater than or equal to a fraction $1 -1/q^k$ and $1 -1/q^{k -t}$, respectively, of $\mathrm{wor}(I)$.

The $1/q$ differential ratio for $\mathsf{CSP(\mathcal{O}_q)}$ results from the fact that the average value over {\em all} solutions can be computed using just $q$ specific solutions, and that for any solution $x$, there exists such a set of $q$ solutions containing $x$. Inspired by this singular case, we evaluate the average differential ratio for  $\mathsf{k\,CSP\!-\!q}$ and its restrictions $\mathsf{k\,CSP(\mathcal{E}_q)}$ and $\mathsf{k\,CSP(\mathcal{I}_q^t)}$. 
We adopt a kind of neighborhood approach: we associate with each solution $x$ of $I$ a multiset $\mathcal{X}(I, x)$ of solutions having the same mean solution value as the set of solutions, with relatively small size $R$, in which $x$ appears a certain number $R^* >0$ of times. 
Considering $\mathcal{X}(I, x^*)$ where $x^*$ is optimal, we find that the average differential ratio on $I$ is at least $R^*/R$. 

\subsection{Partition-based solution multisets}\label{sec-E-S}

Given an instance $I$ of $\mathsf{CSP\!-\!q}$, we introduce the following framework to construct our multisets of solutions $\mathcal{X}(I, x)$, $x\in\Sigma_q^n$:

\smallskip
$\bullet$ {\bf Solution multiset association.} Given a partition $\mathcal{V} =\set{V_1 ,\ldots, V_\nu}$ of $[n]$, a solution $x\in\Sigma_q^n$, and a vector $u\in\Sigma_q^\nu$, we define the solution $y(\mathcal{V}, x, u)$ by: 
\begin{align}\nonumber
	\left(y(\mathcal{V}, x, u)_{V_1}, \ldots, y(\mathcal{V}, x, u)_{V_\nu}\right)
	&=	\left(x_{V_1} +\mathbf{u_1}, \ldots, x_{V_\nu} +\mathbf{u_{\nu}}\right)
\end{align}
That is, the solution $y(\mathcal{V}, x, u)$ is obtained from $x$ by shifting each of its coordinates in $V_c$ by $u_c$, for each $c \in[\nu]$. In particular, if $u$ is the zero vector, then $y(\mathcal{V}, x, \mathbf{0}) =x$. 

We then consider arrays with $\nu$ columns and entries in $\Sigma_q$. With such an $R\times \nu$ array $M$, we associate the solution multiset:
\begin{align}\nonumber
	\begin{array}{rl}
		\mathcal{X}(I, x)	&=\left(y(\mathcal{V}, x, M_r)\,|\,r\in[R]\right)
	\end{array}
\end{align}

$\bullet$ {\bf Conditions.}
To ensure that the solution values $v(I, y(\mathcal{V}, x, M_r))$, $r\in[R]$ cover $\mathrm{opt}(I)$ provided that $x$ is optimal, we require that $M$ contains at least one row of zeros. 
Note that for any $(u, a)\in\Sigma_q^\nu\times\Sigma_q$, the expressions $y(\mathcal{V}, x, u +\mathbf{a})$ and $y(\mathcal{V}, x, u) +\mathbf{a}$ refer to the same solution. 
Thus, when considering the restriction $\mathsf{CSP(\mathcal{E}_q)}$ of $\mathsf{CSP\!-\!q}$, we only require that $M$ contains at least one row of the form $\mathbf{a}$. 

Since our ultimate goal is to connect the mean solution value to the value of an optimal solution, the mean of the solution values over $(y(\mathcal{V}, x, M_r)\,|\,r\in[R])$ should match the mean of the solution values over $\Sigma_q^n$. Formally, we need $M$ to satisfy:
\begin{align}\label{eq-E-M=E}
	\textstyle\sum_{r =1}^R v(I, y(\mathcal{V}, x, 	M_r))/R		&\textstyle=\mathbb{E}_X[v(I, X)],	&x\in\Sigma_q^n	
\end{align}

When such a case occurs, the average differential ratio on $I$ is at least $\mu^M(\mathbf{0})$, the frequency of $\mathbf{0}$ in $M$. 
Indeed, let $x^*$ be an optimal solution of $I$. We assume {\em without loss of generality} (w.l.o.g.) that the goal on $I$ is to maximize. Then we have:
$$\begin{array}{rl}
	\mathbb{E}_X[v(I, X)]
		&= 		\sum_{r =1}^R v(I, y(\mathcal{V}, x^*, M_r)) / R\qquad\text{by \cref{eq-E-M=E}}\\
		&\geq 	\mu^M(\mathbf{0})\times v(I, x^*) +\left(1 -\mu^M(\mathbf{0})\right)\times\mathrm{wor}(I)
\end{array}$$
If $v(I, .)\in\mathcal{E}_q$, we similarly obtain the lower bound of $\sum_{a =0}^{q -1}\mu^M(\mathbf{a})$, equivalently expressed as $q\times\mu^M_E(\mathbf{0})$, for the average differential ratio.

\medskip
Suppose that $\mathcal{V}$ is either given or can be computed in polynomial time. Then choosing a solution of maximum value over $\set{y(\mathcal{V}, \mathbf{0}, M_r)\,|\,r\in[R]}$ gives the same differential approximation guarantee as the average solution value. 
For example, on an instance $I$ of $\mathsf{CSP(\mathcal{O}_q)}$, relation \cref{eq-E-O_q} suggests to consider the partition $\mathcal{V} =\set{[n]}$ of $[n]$, and the array $M$ on $\Sigma_q$ defined by:
$$\begin{array}{rl}	
	M &\displaystyle	= \begin{pmatrix}
							0\\
							1\\
							\vdots\\
							q -1
						\end{pmatrix}
\end{array}$$
For this array, we have $\mu^M(\mathbf{0}) =1/q$. Solutions of the form $\mathbf{a}$ that achieve the best objective value are therefore $1/q$-differential approximate. 

Note that given any $q\geq 2$, the predicate on $\Sigma_q^3$ that accepts the solutions of the equation $(y_1 +y_2 -y_3 \equiv 0\bmod{q})$ belongs to $\mathcal{O}_q\cap\mathcal{I}_q^2$. Since the $q^2$ accepting entries of this predicate form a subgroup of $\mathbb{Z}_q^3$, it follows from \cite{C13} that $\mathsf{E3\,CSP(\mathcal{O}_q\cap\mathcal{I}_q^2)}$ is $\mathbf{NP\!-\!hard}$ to differentially approximate within any constant factor greater than $1/q$ in tripartite instances. 
Thus, assuming $\mathbf{P}\neq\mathbf{NP}$, no polynomial-time algorithm can outperform the trivial strategy of choosing a solution with maximum value over $\set{\mathbf{0} ,\ldots, \mathbf{q -1}}$.

\subsection{Connecting the average differential ratio to combinatorial designs}

\begin{definition}[see {\em e.g.} \cite{HSS99}]\label{def-OA}
Let $q$, $t$, $\nu\geq t$, and $R$ be four positive integers, where $R$ is a multiple of $q^t$. 
Then an $R\times\nu$ array $M$ with coefficients in $\Sigma_q$ is a {\em $q$-levels Orthogonal Array of strength $t$ with $\nu$ factors and $R$ runs}, an $OA(R, \nu, q, t)$ for short, if given any sequence $J =(c_1 ,\ldots, c_t)$ of pairwise distinct column indices, the rows of the subarray $M^J =(M^{c_1} ,\ldots, M^{c_t})$ coincide with each $u\in\Sigma_q^t$ equally often. Formally, $M$ must satisfy:
\begin{align}\label{eq-OA}
	|\set{r\in[R]\,|\,M_r^J = v}| 	&=R/q^t,	
		&J =(j_1 ,\ldots, j_t)\in[\nu]^t,\,j_1 <\ldots<j_t,\ v\in\Sigma_q^t
\end{align}%
\end{definition}%

Equivalently, an array $M$ is an $OA(R, \nu, q, t)$ {\em if and only if} the distribution $\mu^M$ it defines on $\Sigma_q^\nu$ is balanced $t$-wise independent, meaning that $\mu^M\in\mathcal{I}_q^t$.
Thus, rational-valued balanced $t$-wise independent distributions over $\Sigma_q^\nu$ and orthogonal arrays of strength $t$ with $\nu$ factors on the symbol set $\Sigma_q$ are essentially equivalent, interpreting the frequency of words from $\Sigma_q^\nu$ in the array as a distribution over $\Sigma_q^\nu$. 
In particular, balanced $t$-wise independent subsets $\mathcal{Y}$ of $\Sigma_q^\nu$ correspond exactly to simple orthogonal arrays of strength $t$ with $\nu$ columns and coefficients in $\Sigma_q$. 
For instance, the accepting entries of $ZeroSum^{t +1, q}$ form the rows of a simple $OA(q^t, t +1, q, t)$ on the symbol set $\Sigma_q$. The corresponding array for $q =3$ and $t =2$ is shown on the left side of \cref{tab-OA-DS-triv}.

\begin{definition}[see {\em e.g.} \cite{HSS99}]\label{def-DS}
Let $q\geq 2$, $t >0$, $\nu\geq t$, and $R$ be four integers, where $R$ is a multiple of $q^{t -1}$. 
An $R\times\nu$ array $M$ with coefficients in $\Sigma_q$ is a {\em Difference scheme of strength $t$} based on $(\mathbb{Z}_q, +)$, a $D_t(R, \nu, q)$ for short, if given any sequence $J =(c_1 ,\ldots, c_t)$ of pairwise distinct column indices, the rows of the subarray $M^J$ lie equally often on each subset $\set{u, u +\mathbf{1} ,\ldots, u +\mathbf{q -1}}$, $u\in\Sigma_q^t$ of words. Formally, $M$ must satisfy:
\begin{align}\label{eq-OA_E}
\begin{array}{rl}
\card{\set{r\in[R]\,|\,M_r^J\in\set{v, v+\mathbf{1} ,\ldots, v+\mathbf{q -1}}}} 
&=R/q^{t -1},\qquad\qquad\\\multicolumn{2}{r}{\qquad\qquad\qquad
	J =(j_1 ,\ldots, j_t)\in[\nu]^t,\,j_1 <\ldots<j_t,\ v\in\set{0}\times\Sigma_q^{t -1}}
\end{array}
\end{align}%
\end{definition}%

Equivalently, an array $M$ is a $D_t(R, \nu, q)$ {\em if and only if} the probability distribution $\mu^M_E$ it defines on $\Sigma_q^\nu$ is balanced $t$-wise independent, meaning that $\mu^M_E\in\mathcal{E}_q^t\cap\mathcal{I}_q^t$. 
\Cref{tab-OA-DS-triv} presents the trivial $D_t(q^{t -1}, t, q)$ for $q =3$ and $t =2$. Difference schemes can be viewed as a slight relaxation of orthogonal arrays. For some insight into such arrays and their connections to orthogonal arrays, we invite the reader to refer to \cite{HSS99}.

\label{sec-def_OA+DS} 
\begin{table}\footnotesize{\begin{center} 
\setlength\arraycolsep{3pt}
$\begin{array}{cc|cc}
\begin{array}{ccc}
	M^1		&M^2	&M^3\\\hline
	0		&0		&0\\
	0		&1		&2\\
	0		&2		&1\\
	1		&0		&2\\
	1		&1		&1\\
	1		&2		&0\\
	2		&0		&1\\
	2		&1		&0\\
	2		&2		&2
\end{array}&&&\begin{array}{cc}
	N^1		&N^2	\\\hline
	0		&0\\
	0		&1\\
	0		&2\\
	\\
	\\
	\\
	\\
	\\
	\\
\end{array}
\end{array}$
\caption{\label{tab-OA-DS-triv}
An $OA(3^2, 3, 3, 2)$ (on the left) and a $D_2(3^1, 2, 3)$ (on the right).
The rows of the OA are the accepting entries of $ZeroSum^{3, 3}$.}
\end{center}}\end{table}

\smallskip
Returning to an instance $I$ of $\mathsf{k\,CSP\!-\!q}$, our goal is to construct pairs $(\mathcal{V}, M)$ that satisfy \cref{eq-E-M=E}.
A sufficient condition for the average solution value over $(y(\mathcal{V}, x, M_r)\,|\,r\in[R])$ to match the average solution value over $\Sigma_q^n$ is that, for each constraint $P_i(x_{J_i})$ of $I$, its average value over $(y(\mathcal{V}, x, M_r)_{J_i}\,|\,r\in[R])$ equals $r_{P_i}$. Thus, our focus will be on finding pairs $(\mathcal{V}, M)$ that satisfy:
\begin{align}\label{eq-E-P_i}
	\textstyle\sum_{r =1}^R P_i(y(\mathcal{V}, x, M_r)_{J_i})/R	&= r_{P_i},	&i\in[m],\ x\in\Sigma_q^n
\end{align}

Let $P_i(x_{J_i}) =P_i(x_{i_1} ,\ldots, x_{i_{k_i}})$ be a constraint of $I$, and $x\in\Sigma_q^n$. A sufficient condition for a pair $(\mathcal{V}, M)$ to satisfy \cref{eq-E-P_i} at $(i, x)$ is that, over the solution multiset $(y(\mathcal{V}, x, M_r)\,|\,r\in[R])$, $P_i$ is taken the same number of times at each of its possible entries. 
This condition requires that any two different indices $j, h\in J_i$ belong to two different sets of $\mathcal{V}$. Otherwise, the difference $y(\mathcal{V}, x, M_r)_j -y(\mathcal{V}, x, M_r)_h$ would remain constant over all $r\in [R]$, namely equal to $x_j -x_h$. Therefore, we assume that $\mathcal{V} =\set{V_1 ,\ldots, V_{\nu}}$ is a strong coloring of $I$, and that $M$ is an array with $\nu$ columns.

We denote by $H =(c_1 ,\ldots, c_{k_i})$ the sequence of color indices satisfying $(i_1 ,\ldots, i_{k_i})\in V_{c_1} \times\ldots\times V_{c_{k_i}}$. Let $r\in[R]$. When evaluating $v(I, .)$ at $y(\mathcal{V}, x, M_r)$, $P_i$ is taken at $y(\mathcal{V}, x, M_r)_{J_i} =x_{J_i} +M^H_r$. We observe that the vectors $x_{J_i} +M^H_r$, $r\in[R]$ coincide equally often with each $v\in\Sigma_q^{k_i}$ {\em if and only if} the words $M^H_r$, $r\in[R]$ coincide equally often with each $v\in\Sigma_q^{k_i}$. Since $\card{H} =k_i\leq k$, we deduce that $(\mathcal{V}, M)$ satisfies \cref{eq-E-P_i} provided that $M$ is an orthogonal array of strength $k$.

Now assume that $P_i\in\mathcal{E}_q$. By definition, this implies that $P_i$ evaluates to the same value for any two entries $(y_1 ,\ldots, y_{k_i})$ and $(y_1 +a ,\ldots, y_{k_i} +a)$, where $a\in\Sigma_q$. In this case, a sufficient condition for $(\mathcal{V}, M)$ to satisfy \cref{eq-E-P_i} at $(i, x)$ is that the vectors $x_{J_i} +M^H_r$, $r\in[R]$ lie with equal frequency in the $q^{k_i -1}$ subsets $\set{v, v +\mathbf{1} ,\ldots, v +\mathbf{q -1}}$, $v\in\set{0}\times\Sigma_q^{k_i -1}$ of $\Sigma_q^{k_i}$. This occurs {\em if and only if} the vectors $M^H_r$, $r\in[R]$ also lie with equal frequency in each of these subsets. Thus $(\mathcal{V}, M)$ satisfies \cref{eq-E-P_i} in the case where $M$ is a difference scheme of strength $k$.

Finally, we assume that $P_i$ is balanced $t$-wise independent, where $t$ is some integer in $[k -1]$. This implies that we can fix the value of up to $t$ variables with index $j\in J_i$ (for instance, setting them to $x_j$), and still, when averaging the value taken by $P_i$ over all possible assignments for the remaining variables, obtain the average value of $P_i$. Therefore, instead of $\mathcal{V}$, we consider a new partition $\mathcal{U} =\set{V_1 ,\ldots, V_{\nu -t}, U_{\nu -t +1}}$ of $[n]$, where the last subset aggregates the remaining color sets, i.e.: $U_{\nu -t +1} =V_{\nu -t +1}\cup\ldots\cup V_\nu$. Furthermore, we define $M$ as an array with $\nu -t +1$ columns, the last of which contains only zeros.
Under these assumptions, given $r\in[R]$, $y(\mathcal{V}, x, M_r)_{J_i}$ is the vector $(v_{i_1} ,\ldots, v_{k_i})$ of $\Sigma_q^{k_i}$ defined for each $s\in[k_i]$ by 
$v_{i_s} =x_{i_s}$ if $c_s >\nu -t$ (in which case $i_s\in U_{\nu -t +1}$) and $v_{i_s} =x_{i_s} +M_r^{c_s}$ otherwise. We denote by $L$ the subsequence of indices of $H$ that lie in $[\nu -t]$, by $\kappa$ its length. Note that $\card{H\backslash L}\leq t$. Therefore, a sufficient condition for $(\mathcal{U}, M)$ to satisfy \cref{eq-E-P_i} at $(i, x)$ is that the vectors $M^L_r$, $r\in[R]$ coincide equally often with each $v\in\Sigma_q^\kappa$. Since $\kappa\leq\min\set{k_i, \nu -t}\leq\min\set{k, \nu -t}$, we conclude that $(\mathcal{U}, M)$ satisfies \cref{eq-E-P_i} if the $\nu -t$ first columns of $M$ form an orthogonal array of strength $\min\set{k, \nu -t}$. 

\smallskip
For example, on a $k$-partite instance $I$ of $\mathsf{k\,CSP(\mathcal{I}_q^{k -1})}$, we consider the partition $\mathcal{U} =\set{V_1, [n]\setminus V_1}$ of $[n]$ where $V_1$ is one of the color sets of a strong coloring of $I$, and the array $M$ on $\Sigma_q$ defined by:
\begin{align}\nonumber
	M 	&=\begin{pmatrix}
			0		&0\\
			1		&0\\
			\vdots	&\vdots\\
			q -1		&0
		\end{pmatrix}
\end{align} 
It follows from the previous analysis that on $I$, the average solution value is $1/q$-differential approximate. Since $ZeroSum^{k, q}$ belongs to $\mathcal{I}_q^{k -1}$ and fulfills the requirements of \cref{thm-inapx-kCSP}, for all integers $q\geq 2$ and $k\geq 3$, this constant factor is optimal assuming $\mathbf{P}\neq\mathbf{NP}$. 

Shifting each row of an $OA(R, \nu, q, t)$ $M$ on $\Sigma_q$ by $-u^*$, where $u^*$ is a maximizer of $\mu^M$, yields a new $OA(R, \nu, q, t)$ in which $\mathbf{0}$ is of maximum frequency. 
Therefore, given an $OA(R, \nu, q, t)$ $M$, we can always assume {\em w.l.o.g.} that $\mu^M$ is maximized at $\mathbf{0}$.
Likewise, given an $R\times\nu$ difference scheme $M$ on $\Sigma_q$, we can assume {\em w.l.o.g.} that $\mu^M_E$ is maximized at $\mathbf{0}$. 
In other words, as far as such arrays are concerned, we are interested in the {\em maximum frequencies} rather than the frequency of a particular word $v$ or a specific family $\set{v, v +\mathbf{1} ,\ldots, v +\mathbf{q -1}}$ of words. 
We introduce the following numbers:

\begin{definition}\label{def-E-OA}
For three positive integers $q$, $\nu$, and $t\in[\nu]$, we define $\rho(\nu, q, t)$ as the largest number $\rho$ for which there exists an orthogonal array $M$ on the symbol set $\Sigma_q$ with $\nu$ factors, strength $t$, and maximum frequency $\rho$.

Similarly, we define $\rho_E(\nu, q, t)$ as the largest number $\rho$ for which there exists a difference scheme $M$ with $\nu$ columns, strength $t$, and entries in $\Sigma_q$ such that:
\begin{align}\nonumber
	\max\nolimits_{v\in\set{0}\times\Sigma_q^{\nu -1}} \left\{
		q\times\mu^M_E(v) :=\mu^M(v) +\mu^M(v +\mathbf{1}) +\ldots+ \mu^M(v +\mathbf{q -1})	
	\right\} 	&=\rho
\end{align} 
\end{definition}%

\Cref{tab-OA-DS-ex-q=2,tab-OA-DS-ex-q=3} show a few arrays that achieve either $\rho(\nu, q, t)$ or $\rho_E(\nu, q, t)$. 
The preceding discussion establishes the following connection between these numbers and the average differential ratio on $\mathsf{k\,CSP\!-\!q}$ instances:

\begin{theorem}\label{thm-E-OA}
For all integers $q\geq 2$, $k\geq 2$ and $t\in[k -1]$, on any instance of $\mathsf{k\,CSP\!-\!q}$, $\mathsf{k\,CSP(\mathcal{E}_q)}$, and $\mathsf{k\,CSP(\mathcal{I}_q^t)}$ with a strong chromatic number $\nu\geq k$, the average differential ratio is at least $\rho(\nu, q, k)$, $\rho_E(\nu, q, k)$, and $\rho(\nu -t, q, \min\set{k, \nu -t})$, respectively. 
\end{theorem}%

\begin{table}\footnotesize{\begin{center}
\setlength\arraycolsep{3pt}
$\begin{array}{cp{5pt}cp{5pt}cp{5pt}c}
\rho_E(3, 3, 2) =1/3	&&\rho(3, 3, 2) =1/9	&&\rho_E(4, 3, 2) =1/5	&&\rho(4, 3, 2) =1/9\\[3pt]
\begin{array}{ccc}
	M^1 &M^2&M^3\\\hline
	0	&0	&0\\
	0	&1	&2\\
	0	&2	&1\\
	\\\\\\\\\\\\
	\\\\\\\\\\\\
\end{array}
&&\begin{array}{ccc}
	N^1 &N^2&N^3\\\hline
	0	&0	&0\\
	0	&1	&2\\
	0	&2	&1\\
	1	&1	&1\\
	1	&2	&0\\
	1	&0	&2\\
	2	&2	&2\\
	2	&0	&1\\
	2	&1	&0\\
	\\\\\\\\\\\\
\end{array}
&&\begin{array}{cccc}	
	O^1 &O^2&O^3&O^4\\\hline
	0	&0	&0	&0\\
	0	&0	&0	&0\\
	0	&0	&0	&0\\
	0	&0	&1	&2\\
	0	&0	&2	&1\\
	0	&1	&0	&2\\
	0	&1	&1	&2\\
	0	&1	&2	&0\\
	0	&1	&2	&1\\
	0	&1	&2	&2\\
	0	&2	&0	&1\\
	0	&2	&1	&0\\
	0	&2	&1	&1\\
	0	&2	&1	&2\\
	0	&2	&2	&1\\
\end{array}
&&\begin{array}{cccc}
	P^1 &P^2&P^3&P^4\\\hline	
	0	&0	&0	&0\\
	0	&1	&2	&2\\
	0	&2	&1	&1\\
	1	&1	&1	&0\\
	1	&2	&0	&2\\
	1	&0	&2	&1\\
	2	&2	&2	&0\\
	2	&0	&1	&2\\
	2	&1	&0	&1\\
	\\\\\\\\\\\\
\end{array}
\end{array}$
\caption{\label{tab-OA-DS-ex-q=3}
Arrays realizing $\rho_E(\nu, 3, 2)$ and $\rho(\nu, 3, 2)$ for $\nu\in\set{3, 4}$:
for the OAs $N$ and $P$, we know from \cref{thm-E-2CSPs} that $\rho(3, 3, 2) =\rho(4, 3, 2) =1/9$; 
for the DS $M$, we deduce from \cref{eq-OA-E-F} that $\rho_E(3, 3, 2)\leq 3\times\rho(3, 3, 2) =1/3$;
the DS $O$ was computed by computer (see \cref{sec-cd} for more details). 
With respect to relations \cref{eq-OA-E-F}, we observe that the inequalities
	$3\times\rho(4, 3, 2) >\rho_E(4, 3, 2)$ and $\rho_E(4, 3, 2) >\rho(3, 3, 2)$ are strict,
	while the inequality $3\times\rho(3, 3, 2)\geq\rho_E(3, 3, 2)$ is not.
}\end{center}}\end{table}

\subsection{Lower bounds for $\rho(\nu, q, k)$ and $\rho_E(\nu, q, k)$}\label{sec-rho+E}

We identify lower bounds for $\rho(\nu, q, t)$ and $\rho_E(\nu, q, t)$ induced by orthogonal arrays and difference schemes, mostly simple, from the literature. 
$F(\nu, q, t)$ refers to the minimum number of runs required for an orthogonal array of strength $t$, with $\nu$ factors and entries in a set of $q$ distinct symbols \cite{HSS99}. 
Similarly, we denote by $E(\nu, q, t)$ the minimum number of rows in a difference scheme with $\nu$ columns and strength $t$ based on $(\mathbb{Z}_q, +)$. 
Note that for all triples $(\nu, q, t)$ of positive integers, we have the trivial inequalities:
\begin{align}\nonumber
\begin{array}{rlcrl}
	\rho(\nu, q, t) &\geq 1/F(\nu, q, t),	&&\rho_E(\nu, q, t) &\geq 1/E(\nu, q, t)
\end{array}
\end{align}

First, we present useful known relations between the numbers $F(\nu, q, t)$ and $E(\nu, q, t)$, which naturally extend to the numbers $\rho(\nu, q, t)$ and $\rho_E(\nu, q, t)$. These relations are based on the following two properties:

\begin{property}[see {\em e.g.} \cite{HSS99}]\label{pty-OA-DS}
Let $M$ be an $R\times\nu$ array on $\Sigma_q$. 
We consider the arrays:
\begin{align}
A(M)	&:=(M_r^{[\nu -1]}\,|\,r\in[R]: M_r^\nu =0)									\label{eq-OA-OA-}\\
B(M)	&:=\left((M_r, 0)\,|\,r\in[R]\right)										\label{eq-OA-DS+}\\
C(M)	&:=\left(M_r +\mathbf{a}\,|\,r\in[R],\,a\in\set{0, 1,\ldots, q -1}\right)	\label{eq-DS-OA_E}
\end{align}

The following facts hold:
\begin{enumerate}
	\item \label{it-OA-OA-} if $M$ is an $OA(R, \nu, q, t)$, then $A(M)$ is an $OA(R/q, \nu -1, q, t -1)$;
	\item \label{it-OA-DS+} if $M$ is an $OA(R, \nu, q, t)$, then $B(M)$ is a $D_t(R, \nu +1, q)$;
	\item \label{it-DS-OA_E} $M$ is a $D_t(R, \nu, q)$ {\em if and only if} $C(M)$ is an $OA(q\times R, \nu, q, t)$.
\end{enumerate}
\end{property}
For example, on both sides of \cref{tab-OA-DS-ex-q=2}, the array $N$ is the map by $B$ of the array $M$. In \cref{tab-OA-DS-ex-q=3}, the array $N$ is the map by $C$ of the array $M$.
\begin{table}\footnotesize{\begin{center} 
\setlength\arraycolsep{3pt}	
$\begin{array}{cp{5pt}|p{5pt}cp{3pt}c}
	\begin{array}{c}
		\rho(3, 2, 2) =1/4\\[3pt]
		\begin{array}{ccc}
			M^1		&M^2	&M^3	\\\hline
			0		&0		&0		\\
			0		&1		&1		\\
			1		&0		&1		\\
			1		&1		&0		\\
		\end{array}\\\\\\[-3pt]
		\rho_E(4, 2, 2) =\rho_E(4, 2, 3) =1/4\\[3pt]
		\begin{array}{cccc}
			N^1		&N^2	&N^3	&N^4	\\\hline
			0		&0		&0		&0		\\
			0		&1		&1		&0		\\
			1		&0		&1		&0		\\
			1		&1		&0		&0		\\
		\end{array}\\
	\end{array}
	&&&\begin{array}{c}
		\rho(4, 2, 2) =1/6\\[3pt]	
		\begin{array}{cccc}
			M^1		&M^2	&M^3	&M^4	\\\hline
			0		&0		&0		&0		\\
			0		&0		&0		&0		\\
			0		&0		&1		&1		\\
			0		&1		&0		&1		\\
			0		&1		&1		&0		\\
			0		&1		&1		&1		\\
			1		&0		&0		&1		\\
			1		&0		&1		&0		\\
			1		&0		&1		&1		\\
			1		&1		&0		&0		\\
			1		&1		&0		&1		\\
			1		&1		&1		&0		\\
		\end{array}
	\end{array}
	&&\begin{array}{c}
		\rho_E(5, 2, 2) =\rho_E(5, 2, 3) =1/6\\[3pt]	
		\begin{array}{ccccc}
			N^1		&N^2	&N^3	&N^4	&N^5\\\hline
			0		&0		&0		&0		&0	\\
			0		&0		&0		&0		&0	\\
			0		&0		&1		&1		&0	\\
			0		&1		&0		&1		&0	\\
			0		&1		&1		&0		&0	\\
			0		&1		&1		&1		&0	\\
			1		&0		&0		&1		&0	\\
			1		&0		&1		&0		&0	\\
			1		&0		&1		&1		&0	\\
			1		&1		&0		&0		&0	\\
			1		&1		&0		&1		&0	\\
			1		&1		&1		&0		&0	\\
		\end{array}
	\end{array}
\end{array}$
\caption{\label{tab-OA-DS-ex-q=2}
Two pairs $(M, N)$ of an OA $M$ and a DS $N$, where $N$ results from applying the transformation \cref{eq-OA-DS+} to $M$, illustrating the equality $\rho_E(\nu, 2, 2t) =\rho(\nu -1, 2, 2t)$ of relation \cref{eq-OA-E-F-bin} with $t =1$ for $\nu =4$ (on the left) and $\nu =5$ (on the right). 
In both cases, \cref{thm-E-2CSPs} implies that $M$ realizes $\rho(\nu -1, 2, 2)$, which with  \cref{eq-OA-E-F-bin} implies that $N$ realizes $\rho_E(\nu, 2, 2)$. 
Furthermore, if we applied the transformation \cref{eq-DS-OA_E} to the array $N$, we would obtain an OA that realizes $\rho(\nu, 2, 3)$, which according to \cref{eq-OA-E-F-bin} is equal to $\rho_E(\nu, 2, 2)/2$.}
\end{center}}\end{table}

\begin{property}[see {\em e.g.} \cite{HSS99}]\label{pty-DS-2}
If $M$ is a difference scheme with Boolean entries of even strength $t$, then it is a difference scheme of strength $t +1$. 
\end{property}

\Cref{pty-OA-DS,pty-DS-2} obviously imply the following relations between the minimum number of rows on the one hand, and the maximum frequencies on the other hand, in orthogonal arrays and difference schemes of different orders:

\begin{proposition}\label{prop-OA-E-F}
For three integers $q\geq 2$, $t\geq 1$, and $\nu >t$, we have:
\begin{align}\label{eq-OA-E-F}
&\left\{\begin{array}{rlll}
	E(\nu, q, t)			&\leq F(\nu -1, q, t)		&\leq 	1/q\times F(\nu, q, t +1) 	&\leq E(\nu, q, t +1)\\
	\rho_E(\nu, q, t)	&\geq \rho(\nu -1, q, t)		&\geq 	q\times \rho(\nu, q, t +1)	&\geq \rho_E(\nu, q, t +1)
\end{array}\right.
\end{align}

When $q =2$ and $t$ is even, we even have the equalities:
\begin{align}\label{eq-OA-E-F-bin}
\left\{\begin{array}{rlll}
	E(\nu, 2, t)		&= F(\nu -1, 2, t)		&= F(\nu, 2, t +1)/2		&= E(\nu, 2, t +1)\\
	\rho_E(\nu, 2, t)	&= \rho(\nu -1, 2, t)	&= 2 \rho(\nu, 2, t +1)		&= \rho_E(\nu, 2, t +1)
\end{array}\right.
\end{align}
\end{proposition}
The arrays of \cref{tab-OA-DS-ex-q=2,tab-OA-DS-ex-q=3} provide some illustration of relations \cref{eq-OA-E-F,eq-OA-E-F-bin}.
 %
We now review known upper bounds in the literature for $F(q, \nu, t)$, $E(q, \nu, t)$, $\rho(q, \nu, t)$, and $\rho_E(q, \nu, t)$, starting with small values of $\nu$. For $\nu\in\set{t, t +1}$, we have $F(t +1, q, t) =F(t, q, t) =q^t$ and $E(t, q, t) =q^{t -1}$. In \cite{B52}, Bush exhibits other triples $(\nu, q, t)$ for which $F(\nu, q, t)$ equals $q^t$:
\begin{theorem}[\cite{B52}]\label{thm-Bush}
Let $q\geq 2$, $t \geq 2$, and $\nu\geq t$ be three integers. 
Then $F(\nu, q, t) =q^t$ if 
		$\nu\leq t +1$, 	
	or 	$q$ is a prime power greater than $t$ and $\nu\leq q +1$,
	or 	$t =3$, $q$ is a power of 2 greater than 3, and $\nu\leq q +2$.
\end{theorem}

For larger integers $\nu$ and $t =2$, Colbourn et al. in \cite{CSV19} explicitly study orthogonal arrays that maximize their maximum frequency. They notably establish the following result:
\begin{theorem}[\cite{CSV19}]\label{thm-E-2CSPs}
Let $q \geq 2$ and $\nu \geq q$ be two integers such that $\nu$ is 1 or 0 modulo $q$. 
Then $1/\rho(\nu, q, 2)$ is equal to: 
	$$\left\{\begin{array}{ll}
		\nu (q -1) +1	&\text{if $\nu\equiv 1\bmod{q}$}\\
		\nu (q -1) +q	&\text{if $\nu\equiv 0\bmod{q}$}
	\end{array}\right.$$
\end{theorem}

When $q =2$, \cref{thm-E-2CSPs} provides the exact value of $\rho(\nu, 2, 2)$, which is $1/(\nu +1)$ if $\nu$ is odd and $1/(\nu +2)$ otherwise. Combined with equalities \cref{eq-OA-E-F-bin}, it also gives the exact value of $\rho(\nu, 2, 3)$ and $\rho_E(\nu, 2, 3)$:

\begin{corollary}\label{cor-E-3CSP2}
For all integers $\nu\geq 3$, 
	$1/\rho(\nu, 2, 3) =2\nu$ if $\nu$ is even and $2(\nu +1)$ otherwise.
Equivalently,
	$1/\rho_E(\nu, 2, 3) =\nu$ if $\nu$ is even and $\nu +1$ otherwise.
\end{corollary}

For values of $t$ greater than $2$, upper bounds for $F(\nu, q, t)$ and $E(\nu, q, t)$ are obtained from infinite families of {\em linear codes}.

\begin{definition}
Let $L$, $r$, and $q$ be three positive integers, with $q$ a prime power.
A $q$-ary {\em linear code} $\mathcal{C}$ of {\em length} $L$ and {\em dimension} $r$ is an $r$-dimensional subspace of $\mathbb{F}_q^L$. 
The {\em distance} of $\mathcal{C}$ is the minimum Hamming distance between two of its codewords. 
The {\em dual code} of $\mathcal{C}$ is the set of vectors $v\in\mathbb{F}_q^L$ such that $\sum_{j =1}^L v_j c_j =0$, for every codeword $c\in\mathcal{C}$.
\end{definition}%

First, for binary alphabets (i.e., when $q =2$), we apply Delsarte's Theorem \cite{D73} to binary {\em BCH codes} (named after Bose, Ray-Chaudhuri and Hocquenghem). 
For two positive integers $\kappa\geq 3$ and $t$ such that $2^\kappa -1\geq 2t +1$, the {\em primitive binary BCH code} of length $2^\kappa -1$ and {\em design distance} $2t +1$ is a binary linear code of dimension at least $2^\kappa -1 -t \kappa$ and distance at least $2t +1$ (see, {\em e.g.}, \cite{MWS77}). 
Delsarte's Theorem \cite{D73} states that if $\mathcal{C}$ is a linear code of length $L$, dimension $r$, and distance $d$ over $\mathbb{F}_q$, then the codewords of its dual form a simple $OA(q^{L -r}, L, q, d -1)$. Considering $q =2$, $L =2^\kappa -1$, $r\geq 2^\kappa -1 -t\kappa$, and $d\geq 2t +1$, there thus exists an $OA(R, 2^\kappa -1, 2, 2t)$ with $R\leq 2^{t\kappa}$. 
The following upper bound is therefore valid for $F(\nu, 2, 2t)$:

\begin{theorem}[\cite{H59,BC60,D73}]\label{thm-BCH}
For all integers $t\geq 1$ and $\nu\geq\max\set{2t +1, 7}$ such that $\nu +1$ is a power of 2, we have $F(\nu, 2, 2t)\leq (\nu +1)^t$.
\end{theorem}%

Using the equalities \cref{eq-OA-E-F-bin}, equivalently, for all positive integers $t\geq 1$ and $\nu\geq \max\set{2t +2, 8}$ such that $\nu$ is a power of 2, we have $E(\nu, 2, 2t +1)\leq\nu^t$ and $F(\nu, 2, 2t +1)\leq 2\nu^t$.

For larger prime powers $q$, Bierbrauer constructs in \cite{Bi96} orthogonal arrays that are trace-codes of Reed-Solomon codes. Let $\nu$ and $t$ be two integers such that $q^\nu\geq t\geq 2$, and let $\phi$ be a $\mathbb{F}_q$-linear surjective map from $\mathbb{F}_q^\nu$ to $\mathbb{F}_q$. Then consider the array $B$ defined by:
$$\begin{array}{rll}
B^c_{(a, z)}
	&=\phi\left(\sum_{j =1}^{t -1} a_j c^j\right) +z, 	
	&a =(a_1 ,\ldots, a_{t -1})\in(\mathbb{F}_q^\nu)^{t -1},\ z\in\mathbb{F}_q,\ c\in\mathbb{F}_q^\nu
\end{array}$$

Bierbrauer shows that $B$ is an $OA(q\times q^{\nu(t -1)}, q^\nu, q, t)$, and that the rows of $B$ all have the same multiplicity. He further shows that for some $\phi$ this multiplicity is at least $\lambda\nu$, where $\lambda$ is the largest integer such that $t >q^\lambda$. So there exists an $OA(q\times q^{\nu(t -1 -\lambda)}, q^\nu, q, t)$.
Moreover, we observe that for $z\in\mathbb{F}_q$, the rows $B_{(a, z)}$, $a\in(\mathbb{F}_q^\nu)^{t -1}$ in $B$ are the shift by $\mathbf{z}$ of the rows $B_{(a, 0)}$, $a\in(\mathbb{F}_q^\nu)^{t -1}$. 
If $q$ is prime, the field $\mathbb{F}_q$ is isomorphic to $\mathbb{Z}_q$. Using \cref{it-DS-OA_E} of \cref{pty-OA-DS}, we deduce that in this case, the rows $B_{(a, z)}$ of $B$ with (e.g.) $z =0$ form a $D_t(q^{\nu(t -1)}, q^\nu, q)$. 
These observations imply the following upper bounds on $F(\nu, q, k)$ and $E(\nu, q, k)$:

\begin{theorem}[\cite{Bi96}]\label{thm-Bierb}
For all integers $q\geq 2$, $t \geq 2$, and $\nu\geq t$ where $q$ is a prime power and $\nu$ is a power of $q$, we have $F(\nu, q, t) \leq q\times\nu^{t -\lceil\log_q t\rceil}$ and, if $q$ is a prime, $E(\nu, q, t) \leq \nu^{t -\lceil\log_q t\rceil}$.
\end{theorem}%

\subsection{Derived approximation guarantees for $\mathsf{k\,CSP\!-\!q}$, %
	$\mathsf{CSP(\mathcal{E}_q)}$ and $\mathsf{CSP(\mathcal{I}_q^t)}$}\label{sec-E-dapx}

The arrays discussed in \cref{sec-rho+E} allow us to derive lower bounds on the average differential ratio from \cref{thm-E-OA}.
 %
Although these arrays are typically constructed for values of $q$ that are prime powers, we can extend the  lower bounds they induce on the average differential ratio to arbitrary values of $q$ by reducing to the case where the alphabet size is a prime power. 
\begin{theorem}\label{thm-reduc}
Let $q\geq 2$, $k\geq 2$, and $d >q$ be three integers. Suppose that for $\mathsf{k\,CSP\!-\!d}$, the average differential ratio is lower bounded by some quantity $\rho$, which may depend on the structure of the instance's primary hypergraph; then the same lower bound holds for $\mathsf{k\,CSP\!-\!q}$.
\end{theorem}

\begin{proof}
Consider an instance $I$ of $\mathsf{k\,CSP\!-\!q}$. 
We denote by $\mathcal{M}$ the set of all surjective maps from $\Sigma_d$ to $\Sigma_q$.
Given a sequence $\pi =(\pi_1 ,\ldots, \pi_n)$ of $n$ maps from $\mathcal{M}$, 
we interpret $I$ as the instance $f_\pi(I)$ of $\mathsf{CSP\!-\!d}$, where:
\begin{enumerate}
	\item for each $j\in[n]$, $f_\pi(I)$ contains a variable $z_j$ with domain $\Sigma_d$;
	\item for each $i\in[m]$, $f_\pi(I)$ contains a constraint $P_i (\pi_{i_1}(z_{i_1}) ,\ldots, \pi_{i_{k_i}}(z_{i_{k_i}}))$ with the same weight $w_i$ as $C_i$ in $I$.
\end{enumerate}

By construction, the two instances $f_\pi(I)$ and $I$ share the same primary hypergraph.
To retrieve solutions of $I$ from those of $f_\pi(I)$, we define $g_\pi(I, .)$ by: 
$$\begin{array}{rll}
	g_\pi(I, z) &=(\pi_1(z_1) ,\ldots, \pi_n(z_n)),	&z\in\Sigma_d^n
\end{array}$$ 

This function is surjective and satisfies $v(f_\pi(I), z) =v(I, g_\pi(I, z))$ for all $z\in\Sigma_d^n$.
The extremal solution values on $I$ and $f_\pi(I)$ therefore satisfy:
\begin{align}\label{eq-E-reduc-ext}
\begin{array}{rlrl}
	\mathrm{opt}(f_\pi(I)) &=\mathrm{opt}(I),	&\mathrm{wor}(f_\pi(I)) &=\mathrm{wor}(I)
\end{array}
\end{align}

In contrast, $\mathbb{E}_Z[v(f_\pi(I), Z)]$ may differ from $\mathbb{E}_X[v(I, X)]$ due to the fact that $g_\pi(I, .)$ can associate two vectors $x\neq x'$ of $\Sigma_q^n$ with a different number of vectors from $\Sigma_d^n$. Therefore, instead of $n$ specific maps $\pi_1 ,\ldots, \pi_n$, we consider a sequence $\Pi =(\Pi_1 ,\ldots, \Pi_n)$ of random maps that are independently and uniformly distributed over $\mathcal{M}$. We establish the equality:
\begin{align}\label{eq-tout_q-E}
	\mathbb{E}_\Pi \left[\mathbb{E}_Z[v(f_\Pi(I), Z)]\right]	&= \mathbb{E}_X[v(I, X)]
\end{align}

Let $b\in\Sigma_q$ and $b'\in\Sigma_q\backslash\set{b}$. 
We define the function $\sigma:\mathcal{M}\mapsto\mathcal{M}$ that maps each $\tau\in\mathcal{M}$ to a new function $\sigma(\tau)\in\mathcal{M}$ given by:
$$\begin{array}{rl}
	\sigma(\tau)(c)	&=\left\{\begin{array}{ll}
						b'		&\text{if $\tau(c) =b$}\\
						b		&\text{if $\tau(c) =b'$}\\
						\tau(c)	&\text{otherwise}\\
					\end{array}\right.			
\end{array}$$

The function $\sigma$ is a bijection on $\mathcal{M}$ since it permutes the preimages of two values while leaving all others unchanged. Given any $a\in\Sigma_d$, we have:
$$\begin{array}{rll}
\card{\set{\tau\in\mathcal{M}\,|\,\tau(a) =b}}
	&=\card{\set{\tau\in\mathcal{M}\,|\,\sigma(\tau)(a) =b}}	&\text{since $\sigma$ is a bijective}\\
	&=\card{\set{\tau\in\mathcal{M}\,|\,\tau(a) =b'}}			&\text{by definition of $\sigma$}
\end{array}$$

It follows that for all $a\in\Sigma_d$, the cardinalities $\card{\set{\tau\in\mathcal{M}\,|\,\tau(a) =b}}$, $b\in\Sigma_q$ are all equal to $\card{\mathcal{M}}/q$. Given $z\in\Sigma_d^n$, we successively deduce:
$$\begin{array}{rll}
\mathrm{P}_\Pi[g_\Pi(I, z) =x]	
	&=\prod_{j =1}^n \mathrm{P}_{\Pi_j}[\Pi_j(z_j) =x_j]					\\
	&=\prod_{j =1}^n \left(\card{\set{\tau\in\mathcal{M}\,|\,\tau(z_j) =x_j}}/\card{\mathcal{M}}\right)
	&=1/q^n																	\\[4pt]
\mathbb{E}_\Pi[v(I, g_\Pi(I, z))]	
	&=\sum_{x\in\Sigma_q^n} v(I, x) \times\mathrm{P}_\Pi[g_\Pi(I, z) =x]	\\
	&=\sum_{x\in\Sigma_q^n} v(I, x) \times1/q^n
	&=\mathbb{E}_X[v(I, X)]
\end{array}$$

Thus, the expected objective value of $v(I, g_\Pi(I, z))$ coincides with $\mathbb{E}_X[v(I, X)]$ for every $z\in\Sigma_d^n$, implying that $\mathbb{E}_Z\left[\mathbb{E}_\Pi[v(I, g_\Pi(I, Z))]\right]$ coincides with $\mathbb{E}_X[v(I, X)]$. Equality \cref{eq-tout_q-E} follows, considering:
$$\begin{array}{rlll}
	\mathbb{E}_\Pi\left[\mathbb{E}_Z[v(f_\Pi(I), Z)]\right]	
		&= \mathbb{E}_\Pi\left[\mathbb{E}_Z[v(I, g_\Pi(I, Z))]\right]
		&= \mathbb{E}_Z\left[\mathbb{E}_\Pi[v(I, g_\Pi(I, Z))]\right]
\end{array}$$

By \cref{eq-tout_q-E}, there exists a vector $\pi_*\in\mathcal{M}^n$ such that $\mathbb{E}_X[v(I, X)]\geq\mathbb{E}_Z[v(f_{\pi_*}(I), Z)]$, while given such a $\pi_*$, by \cref{eq-E-reduc-ext} we have:
$$\begin{array}{rl}
\displaystyle	\frac{\mathbb{E}_Z[v(f_{\pi_*}(I), Z)] -\mathrm{wor}(I)} 
					 {\mathrm{opt}(I) -\mathrm{wor}(I)}
&\displaystyle=	\frac{\mathbb{E}_Z[v(f_{\pi_*}(I), Z)] -\mathrm{wor}(f_{\pi_*}(I))}
					 {\mathrm{opt}(f_{\pi_*}(I)) -\mathrm{wor}(f_{\pi_*}(I))}
\end{array}$$

We conclude that the average differential ratio on $I$ is at least the average differential ratio on $f_{\pi_*}(I)$, completing the proof.
\end{proof}

Note that in the general case, the transformation $f_\pi$ does not necessarily map an instance of $\mathsf{CSP(\mathcal{E}_q)}$ or $\mathsf{CSP(\mathcal{I}_q^t)}$ to an instance of $\mathsf{k\,CSP(\mathcal{E}_d)}$ or $\mathsf{k\,CSP(\mathcal{I}_d^t)}$. 
For example, consider the case where $d =3$, $q =2$, and the maps $\pi_1$ and $\pi_2$ both assign the value $a\bmod 2$ to each $a\in\Sigma_3$. 
The function $XNOR^2$ on $\Sigma_2^2$ belongs to $\mathcal{E}_2$ (and thus, to $\mathcal{I}_2^1$), while the function $P$ on $\Sigma_3^2$ that assigns the value $XNOR^2(a \bmod 2, b \bmod 2)$ to each $(a, b)\in\Sigma_3^2$ does not belong to $\mathcal{I}_3^1$ (and thus, not to $\mathcal{E}_3$). 
For instance, we have
	$P(0, 0) +P(0, 1) +P(0, 2) =2\times XNOR^2(0, 0) +XNOR^2(0, 1) =2$, 
	and $P(1, 0) +P(1, 1) +P(1, 2) =2\times XNOR^2(1, 0) +XNOR^2(1, 1) =1$, 
which prevents $P$ from satisfying condition \cref{eq-I_q^t} at rank $t =1$.

\smallskip
Consider an instance $I$ of $\mathsf{k\,CSP\!-\!q}$. We argue that for all integers $s\geq\nu$ and $d\geq q$, the quantity $\rho(s, d, \min\set{s, k})$ provides a valid lower bound for the average differential ratio on $I$. First, we deduce from \cref{thm-E-OA,thm-reduc} that $\rho(\nu, d, \min\set{\nu, k})$ is a proper lower bound on this ratio. Then we observe that the first $\nu$ columns of an $OA(R, s, d, t)$ with $R^*$ repeated rows form an $OA(R, \nu, d, t)$ with at least $R^*$ repeated rows, implying the inequality  $\rho(\nu, d, \min\set{\nu, k})\geq\rho(s, d, \min\set{\nu, k})$. Similarly, an $OA(R, s, d, t)$ is also an $OA(R, s, d, t')$ for all integers $t'\in[t]$, implying $\rho(s, d, \min\set{\nu, k})\geq\rho(s, d, \min\set{s, k})$.

Based on these observations, we derive the following estimates of the average differential ratio on $\mathsf{k\,CSP\!-\!q}$ instances from \cref{thm-Bush,thm-E-2CSPs,cor-E-3CSP2,thm-BCH,thm-Bierb} (a complete proof of \cref{cor-E-k_partite,cor-E-2CSPs,cor-E-3CSPs,cor-E-kCSP-2,cor-E-kCSP-q,cor-E-kCSP(I_q^t)-Bush,cor-E-kCSP(E_q)} can be found in \cref{sec-ap-avd}):

\begin{corollary}[Consequence of \cref{thm-E-OA,thm-reduc,thm-Bush}]\label{cor-E-k_partite}
Let $q\geq 2$ and $k\geq 2$ be two integers, and $I$ be an instance of $\mathsf{k\,CSP\!-\!q}$ with a strong chromatic number $\nu\in\set{k ,\ldots, \max\set{k +1, 2q}}$. We denote by $p^\kappa$ the smallest prime power greater than or equal to $q$. Then on $I$, the average differential ratio is bounded below by: 
\begin{enumerate}
	\item\label{it-B-nu} $1/q^k$ if $\nu\leq k +1$, 
				or $q >k$ is a prime power and $\nu\leq q +1$, 
				or $k =3$, $q >3$ is a power of 2, and $\nu\leq q +2$;
	\item\label{it-B-qpp} $1/p^{\kappa k}$ --- and thus, by $1/(2q -2)^k$ --- if $p^\kappa >k$ and $\nu\leq p^\kappa +1$;
	\item\label{it-B-qp2} $1/2^{3\lceil\log_2 q\rceil}$ --- and thus, by $1/(2q -2)^3$ --- if $k =3$, $q \geq 3$, and $\nu\leq 2^{\lceil\log_2 q\rceil} +2$.
\end{enumerate}
\end{corollary}

\begin{corollary}[Consequence of \cref{thm-E-OA,thm-E-2CSPs}]\label{cor-E-2CSPs}
For all integers $q\geq 2$, on any instance of $\mathsf{2\,CSP\!-\!q}$ with a strong chromatic number $\nu\geq 2$, the average differential ratio is bounded below by:
$$\begin{array}{rll}
	\displaystyle\frac{1}{q\lceil(\nu -1)/q\rceil(q -1) +q}
		&\displaystyle	\sim \frac{1}{(q -1)\nu}	
\end{array}$$
\end{corollary}
In particular for $q =2$, this ratio is at least $1/(\nu +1)$ if $\nu$ is odd and $1/(\nu +2)$ otherwise.

\begin{corollary}[Consequence of \cref{thm-E-OA,cor-E-3CSP2}]\label{cor-E-3CSPs}
On any instance of $\mathsf{3\,CSP\!-\!2}$ with a strong chromatic number $\nu\geq 3$, the average differential ratio is at least $1/\left(4\lceil\nu/2\rceil\right)\sim 1/(2 \nu)$. 
\end{corollary}

\begin{corollary}[Consequence of \cref{thm-E-OA,prop-OA-E-F,thm-BCH}]\label{cor-E-kCSP-2}
For all integers $k\geq 4$, on any instance of $\mathsf{k\,CSP\!-\!2}$ with a strong chromatic number $\nu\geq k$, the average differential ratio is at least:
\begin{center}
$\begin{array}{rll}
1/2^{\lceil\log_2(\nu +1)\rceil k/2}
	&\geq 1/(2 \nu)^{k/2}							&\text{if $k$ is even;}\\ 
1/2\times 1/2^{\lceil\log_2 \nu\rceil(k -1)/2}
	&\geq 1/2\times 1/(2\nu)^{(k -1)/2}				&\text{if $k$ is odd.}
\end{array}$
\end{center}
\end{corollary}

\begin{corollary}[Consequence of \cref{thm-E-OA,thm-reduc,thm-Bierb}]\label{cor-E-kCSP-q}
Let $q\geq 3$ and $k\geq 2$ be two integers, and $I$ be an instance of $\mathsf{k\,CSP\!-\!q}$ with a strong chromatic number $\nu\geq k$. We denote by $p^\kappa$ the smallest prime power such that $p^\kappa\geq q$ (thus $p^\kappa =q$ if $q$ is a prime power and $p^\kappa\leq 2(q -1)$ otherwise).
Then the average differential ratio on $I$ is bounded below by: 
\begin{center}
$\begin{array}{rll}
\displaystyle\frac{1}{p^\kappa}\!\times\!
 	\frac{1}{p^{\kappa\lceil\log_{p^\kappa}\nu\rceil(k -\lceil\log_{p^\kappa} k\rceil)}}
&\geq\displaystyle
	\frac{1}{p^\kappa}\!\times\!\frac{1}{(p^\kappa\nu)^{k -\lceil\log_{p^\kappa} k\rceil}}
\end{array}$
\end{center}
\end{corollary}

If there is some integer $t >0$ for which the constraints each involve a function of $\mathcal{I}_q^t$, then, in accordance with \cref{thm-E-OA}, we can refine the lower bound on the average differential ratio to $\rho(\nu -t, q, \min\set{\nu -t, k})$. This leads to an improvement in our estimate of this ratio for instances of $\mathsf{k\,CSP(\mathcal{I}_q^t)}$ with a small strong chromatic number:

\begin{corollary}[Consequence of \cref{thm-E-OA,thm-Bush}]\label{cor-E-kCSP(I_q^t)-Bush}
Let $q\geq 2$, $k\geq 2$, and $t\in[k -1]$ be three integers where $q$ is a prime power, and $I$ be an instance of $\mathsf{k\,CSP(\mathcal{I}_q^t)}$ with a strong chromatic number $\nu\in\set{k ,\ldots,\max\set{k +t +1, q +t +2}}$. 
Then:

$\bullet$ If $\nu -t\leq k$, the average differential ratio on $I$ is bounded below by 
	$1/q^{\nu -t}$.

$\bullet$ If $\nu -t\leq k +1$, 
			or $q >k$ and $\nu -t\leq q +1$, 
			or $k =3$, $q >3$, $q$ is a power of 2, and $\nu -t\leq q +2$,
this ratio is at least $1/q^k$.
\end{corollary}

Now consider an instance $I$ of $\mathsf{CSP(\mathcal{E}_q)}$, along with an integer $s\geq\nu$. Similar to the numbers $\rho(\nu, q, t)$, we have the inequality 
$\rho_E\left(\nu, q, \min\set{t, \nu}\right) 
	\geq\rho_E\left(s, q, \min\set{t, s}\right)$.
It therefore follows from \cref{thm-E-OA} that $\rho_E\left(s, q, \min\set{k, s}\right)$ is a proper lower bound on the average differential ratio on $I$.  
Given this fact, for $\mathsf{k\,CSP(\mathcal{E}_q)}$ in the case where $q$ is prime and either $k$ is odd or $q\geq 3$, we obtain an improvement in our estimate of this ratio by a multiplicative factor of $q$ over the estimates previously obtained for $\mathsf{k\,CSP\!-\!q}$ and $\mathsf{k\,CSP(\mathcal{I}_q^1)}$:

\begin{corollary}[Consequence of \cref{thm-E-OA,thm-Bush,thm-BCH,thm-Bierb}, \cref{prop-OA-E-F}, and \cref{cor-E-3CSP2}]\label{cor-E-kCSP(E_q)}
Let $q$ and $k\geq 3$ be two integers where $q$ is prime, and $I$ be an instance of $\mathsf{k\,CSP(\mathcal{E}_q)}$ with a strong chromatic number $\nu\geq k$. Then:

$\bullet$ If $q =2$ and $k$ is odd, the average differential ratio on $I$ is at least:
\begin{center}
$\begin{array}{rll}
1/2^{k -1}							&										&\text{if $\nu \leq k +1$};\\
1/\left(2\lceil\nu/2\rceil\right) 	&\sim 1/\nu								&\text{if $k =3$};\\
1/2^{\lceil\log_2\nu\rceil(k -1)/2}	&\geq 1/(2\nu)^{\lfloor k/2\rfloor}	&\text{if $k \geq 5$.}
\end{array}$
\end{center}

$\bullet$ If $q\geq3$, this ratio is bounded below by:
\begin{center}
$\begin{array}{rlll}
 	1/q^{\lceil\log_q\nu\rceil(k -\lceil\log_q k\rceil)} 
 		&\geq 1/(q\nu)^{(k -\lceil\log_q k\rceil)}.
\end{array}$
\end{center}
\end{corollary}

\subsection{Concluding remarks}
\label{sec-E-tight}
\begin{table}[t]\small{
\setlength\arraycolsep{0.15cm}
{\em Approximability bounds in $k$-partite instances of $\mathsf{Ek\,CSP(\mathcal{I}_q^t)}$:}\\
$\begin{array}{lll|ll|l}
\multicolumn{1}{c}{k}	&\multicolumn{1}{c}{q}	&\multicolumn{1}{c|}{t}
	&\multicolumn{1}{c}{\text{gapx det.}}	&\multicolumn{1}{c|}{\text{dapx det.}}			&\multicolumn{1}{c}{\text{avd}}
\\\hline
\multirow{2}{*}{=2}		&\multirow{2}{*}{=2}		&\multirow{2}{*}{=1}
	&2\ln(1 +\sqrt{2})/\pi >0.561	&1/2 +\ln(1 +\sqrt{2})/\pi>0.78\text{ \cite{AN04}}^*
	&\multirow{2}{*}{=1/2}
\\
		&					&
	&\neg\ 3/4 +\varepsilon\text{ \cite{H97}}	&\neg\ 7/8 +\varepsilon\text{ \cite{H97}}	&
\\
\geq 3	&\geq 2				&=k -1
	&\neg\ \varepsilon\text{ \cite{C13}}	&\neg\ 1/q +\varepsilon\text{ \cite{C13}}			&\geq 1/q
\\
\geq 3	&\geq 2,\ \leq k	&=2	
	&\neg\ \varepsilon\text{ \cite{C13}}	&\neg\ O(k/q^{k -1}) +\varepsilon\text{ \cite{C13}}	&\geq 1/q^{k -2}
\\
\geq 3	&\geq k				&=2
	&\neg\ \varepsilon\text{ \cite{C13}}	&\neg\ O(k/q^{k -2}) +\varepsilon\text{ \cite{C13}}	&\geq 1/q^{k -2}
\\\multicolumn{6}{c}{}\\[-3pt]
\end{array}$

{\em Gain approximability bounds for $\mathsf{E3\,CSP(\mathcal{I}_q^2)}$ and $\mathsf{k\,CSP\!-\!q}$:}\\
$\begin{array}{l|ll|l}	
\multicolumn{1}{c|}{\text{Restriction}}	&\multicolumn{1}{c}{\text{gapx det.}}
								&\multicolumn{1}{c|}{\text{gapx exp.}}	&\multicolumn{1}{c}{\text{avd}}
\\\hline
\mathsf{E3\,CSP(\mathcal{I}_2^2)}
	&\Omega(1/m)\text{ \cite{HV04}}^*		&\Omega(\sqrt{\ln n/n})\text{ \cite{KN08}}^*	&=1/2
\\\hline
\mathsf{2\,CSP\!-\!2}			
	&\Omega(1/\ln n)\text{ \cite{NRT99}}	&											&\Omega(1/\nu)
\\
\mathsf{k\,CSP\!-\!2},\ k\geq 3	&\Omega(1/m)\text{ \cite{HV04}}		&\Omega(1/\sqrt{m})\text{ \cite{HV04}}
																		&\Omega(1/\nu^{\lfloor k/2\rfloor})
\\
\mathsf{k\,CSP\!-\!2^\kappa},\ k\geq 2,\,\kappa\geq 2
					&\Omega(1/m)\text{ \cite{HV04,CT18-E}}	&\Omega(1/\sqrt{m})\text{ \cite{HV04,CT18-E}}	
														&\Omega(1/\nu^{k -\lceil\log_{2^\kappa} k\rceil})
\\\multicolumn{4}{c}{}\\[-3pt]
\end{array}$ 

{\em Other differential approximability bounds for $\mathsf{k\,CSP\!-\!q}$:}\\
$\begin{array}{l|ll|l}
\multicolumn{1}{c|}{\text{Conditions on $k$ and $q$}}
	&\multicolumn{1}{c}{\text{dapx det.}}	&\multicolumn{1}{c|}{\text{dapx exp.}}	&\multicolumn{1}{c}{\text{avd}}\\\hline
\text{$k =2$ or $(k, q) =(3, 2)$}	
	&\Omega(1)		\text{ \cite{N98, CT18}}		&	
	&\Omega(1/\nu)\\ 
\text{$k \geq 3$ and $q \geq 3$}	
	&\Omega(1/m)	\text{ \cite{HV04,CT18-E}}		&\Omega(1/\sqrt{m})	\text{ \cite{HV04,CT18-E}}
	&\Omega(1/\nu^{k -\lceil\log_{\Theta(q)} k\rceil})
\end{array}$
\caption{\label{tab-E_cmp} 
	Differential (dapx) and gain (gapx) approximability bounds for $\mathsf{k\,CSP\!-\!q}$  
		that are achievable by either deterministic (det.) or randomized (exp.) algorithms, and their comparison with known lower bounds for the average differential ratio (avd).
The given inapproximability bounds hold for all constant $\varepsilon >0$, provided that $\mathbf{P}\neq\mathbf{NP}$. 
The bounds marked with $^*$ are discussed in \cref{sec-ap-bib}.}}
\end{table}

We consider the following questions: is the average solution value a good approximation of the optimum value? How tight are our estimates of the average differential approximation ratio? How accurate are our bounds on numbers $\rho(\nu, q, t)$ and $\rho_E(\nu, q, t)$? Addressing these questions, we identify potential areas of improvement and directions for future research.

\smallskip
{\bf{\em Quality of the obtained approximability bounds.}}
We evaluate our estimate of the average differential ratio in light of the gain and differential approximability bounds in the literature. We summarize the approximability bounds known to us in \cref{tab-E_cmp}. 
On the one hand, we compare the differential approximation guarantee implied by $\mathbb{E}_X[v(I, X)]$ with those obtained using dedicated algorithms.
On the other hand, we examine how well $\mathbb{E}_X[v(I, X)]$ approximates the advantage over $\mathrm{wor}(I)$ compared to how well the advantage over $\mathbb{E}_X[v(I, X)]$ can be approximated.

For such symptomatic CSPs as the restriction of $\mathsf{Ek\,CSP(\mathcal{I}_q^{k -1})}$ to $k$-partite instances, for any $k\geq 3$, $\mathbb{E}_X[v(I, X)]$ trivially provides a differential approximation guarantee of $1/q$. According to \cite{H97,C13}, this is the best constant factor achievable under the assumption $\mathbf{P}\neq\mathbf{NP}$. For this particular CSP, approximating the optimal advantage over a random assignment within any positive constant is $\mathbf{NP\!-\!hard}$. 
The same facts hold for $\mathsf{CSP(\mathcal{O}_q)}$.

In contrast, for $\mathsf{2\,CSP\!-\!q}$, $\mathbb{E}_X[v(I, X)]$ provides a rather weak approximation, since we next establish that $\Omega(1/n)$ is a tight lower bound on the average differential ratio, while $\mathsf{2\,CSP\!-\!q}$ is approximable within some constant differential factor (see \cref{sec-reduc}). 
In the Boolean case, $\mathsf{2\,CSP\!-\!2}$ is even approximable within a gain approximation ratio of $\Omega(1/\ln n)$ \cite{NRT99,NWY00,M01}. 

For $\mathsf{k\,CSP\!-\!2}$ with $k\geq 3$, $\abs{\mathbb{E}_X[v(I, X)] -\mathrm{wor}(I)}$ approximates the instance diameter within a factor of $\Omega(1/n^{\lfloor k/2\rfloor})$. In dense instances, this is comparable to the expected gain approximability bound of $\Omega(1/\sqrt{m})$ established in \cite{HV04}.

\smallskip
{\bf{\em Tightness of the analysis.}}
Given two integers $k\geq 3$ and $q\geq 2$, for $k$-partite instances of $\mathsf{CSP(\mathcal{I}_q^2)}$, the inapproximability bound we can derive from \cite{C13} is a factor $O(k)$ of the lower bound we obtain for the average differential ratio. For $k$-partite instances of $\mathsf{CSP(\mathcal{I}_q^{k -1})}$, we have already pointed out that the guarantee obtained is basically optimal.
 %
Moreover, we observe that our analysis for $\mathsf{2\,CSP(\mathcal{E}_q)}$ and $\mathsf{3\,CSP(\mathcal{E}_2)}$ is either tight or asymptotically tight. Given three positive integers $q$, $k$, and $n \geq k$, we denote by $I^{q, k}_n$ the instance of $\mathsf{CSP(\set{AllEqual^{k, q}})}$ that considers all the $k$-ary constraints that can be formed over a set of $qn$ variables. This instance is trivially satisfiable, and has a strong chromatic number of $q n$. 
Furthermore, a worst solution on $I^{q, k}_n$ assigns the value $a$ to exactly $n$ variables for each $a\in\Sigma_q$. 
Thus, considering that 
	$\mathrm{opt}(I^{q, k}_n) =\binom{q n}{k}$, 
	$\mathrm{wor}(I^{q, k}_n) =q\times \binom{n}{k}$,
	and $\mathbb{E}_X[v(I^{q, k}_n, X)] =\binom{q n}{k}/q^{k -1}$, 
the average differential ratio on $I^{q, k}_n$ is equal to:
\begin{align}\nonumber
\frac{\binom{q n}{k}/q^{k -1} -q\times\binom{n}{k}}{\binom{q n}{k} -q\times\binom{n}{k}}
\end{align}

The above fraction is $1/(q n)$ when either $k =2$, or $k =3$ and $q =2$. 
Hence, for $k =2$, the average differential ratio on $I^{q, k}_n$ is asymptotically a factor $(q -1)$ of the lower bound \cref{cor-E-2CSPs} provides for this ratio. For $(k, q) =(3, 2)$, this ratio matches the bound given in \cref{cor-E-kCSP(E_q)}.

\begin{table}
\footnotesize
\setlength\arraycolsep{3.8pt}
$$\begin{array}{|c|c|cccccccccccc|}
\multicolumn{14}{l}{\text{\begin{minipage}{0.9\textwidth}
The ratio $R^*/R$ in DS's that minimize $R$ among those realizing $\rho_E(\nu, q, t)$
(for $q =2$ and $t$ odd, we recall that $\rho_E(\nu, 2, t) =\rho_E(\nu, 2, t -1)$ \cref{eq-OA-E-F-bin}):
\end{minipage}}}\\
\multicolumn{2}{c|}{} &\multicolumn{12}{c|}{\nu}\\
\cline{3-14}
q 	&t	&2		&3			&4			&5					&6			&7			&8			&9			&10		&11			&12			&13		\\
\hline \multirow{3}{*}{2}	
	&2	&\gr{1/2}	&\gr{1/4}	&\gr{1/4}	&2/12*			&2/12*		&\gr{1/8}	&\gr{1/8}	&2/20*		&2/20*	&\gr{1/12}	&\gr{1/12}	&2/28*	\\	
	&4	&-			&-			&\gr{1/8}	&\gr{1/16}		&\gr{1/16}	&3/80		&4/144		&6/240		&6/336	&6/336		&			&		\\
	&6	&-			&-			&-			&-				&\gr{1/32}	&\gr{1/64}	&\gr{1/64}	&4/448		&6/960	&25/5184	&			&		\\
\hline \multirow{3}{*}{3}	
	&2	&\gr{1/3}	&\gr{1/3}	&3/15		&\gr{1/6}		&\gr{1/6}	&3/24		&\gr{1/9}	&\gr{1/9}	&3/33	&			&			&		\\
	&3	&-			&\gr{1/9}	&\gr{1/9}	&2/27			&\gr{2/27}	&8/162		&8/162		&			&		&			&			&		\\
	&4	&-			&-			&\gr{1/27}	&\gr{1/27}		&7/297		&5/243		&			&			&		&			&			&		\\
	&5	&-			&-			&-			&\gr{1/81}		&\gr{1/81}	&27/3240	&			&			&		&			&			&		\\
\hline \multirow{2}{*}{4}	
	&2	&\gr{1/4}	&\gr{2/8}	&\gr{2/8}	&4/24			&14/104		&\gr{2/16}	&\gr{2/16}	&			&		&			&			&		\\	
	&3	&-			&\gr{1/16}	&\gr{1/16}	&\gr{2/32}		&\gr{2/32}	&			&			&			&		&			&			&		\\
	&4	&-			&-			&\gr{1/64}	&\gr{2/128}		&			&			&			&			&		&			&			&		\\
\hline
\multicolumn{14}{c}{}\\
\multicolumn{14}{l}{
\text{\begin{minipage}{0.9\textwidth}
The ratio $R^*/R$ in DS's that maximize $R^*$ among those realizing $E(\nu, q, t)$ 
(for $q =2$ and $t$ odd, we recall that $E(\nu, 2, t) =E(\nu, 2, t -1)$ \cref{eq-OA-E-F-bin}):	
\end{minipage}}}\\
\multicolumn{2}{c|}{} &\multicolumn{12}{c|}{\nu}\\
\cline{3-14}
q 	&t	&2			&3			&4			&5				&6			&7			&8			&9			&10		&11			&12			&13		\\
\hline \multirow{3}{*}{2}	
	&2	&\gr{1/2}	&\gr{1/4}	&\gr{1/4}	&1/8			&1/8		&\gr{1/8}	&\gr{1/8}	&1/12		&1/12	&\gr{1/12}	&\gr{1/12}	&1/16	\\
	&4	&-			&-			&\gr{1/8}	&\gr{1/16}		&\gr{1/16}	&1/32		&1/64		&1/64		&		&			&			&		\\
	&6	&-			&-			&-			&-				&\gr{1/32}	&\gr{1/64}	&\gr{1/64}	&1/128		&1/256	&			&			&		\\
\hline \multirow{3}{*}{3}	
	&2	&\gr{1/3}	&\gr{1/3}	&1/6		&\gr{1/6}		&\gr{1/6}	&1/9		&\gr{1/9}	&\gr{1/9}	&1/12	&			&			&		\\
	&3	&-			&\gr{1/9}	&\gr{1/9}	&1/18			&\gr{2/27}	&1/27		&1/27		&			&		&			&			&		\\
	&4	&-			&-			&\gr{1/27}	&\gr{1/27}		&1/81		&1/81		&			&			&		&			&			&		\\
	&5	&-			&-			&-			&\gr{1/81}		&\gr{1/81}	&1/243		&			&			&		&			&			&		\\
\hline \multirow{1}{*}{4}	
	&2	&\gr{1/4}	&\gr{2/8}	&\gr{2/8}	&2/16			&2/16		&\gr{2/16}	&\gr{2/16}	&			&		&			&			&		\\
\hline
\multicolumn{14}{c}{}\\
\multicolumn{14}{l}{
\text{\begin{minipage}{0.9\textwidth}
The ratio $R^*/R$ in OAs that minimize $R$ among those realizing $\rho(\nu, q, t)$ 
(for $q =2$, we recall that $\rho(\nu, 2, t) =\rho_E(\nu +1, 2, t)$ if $t$ is even and $\rho_E(\nu, 2, t -1)/2$ otherwise \cref{eq-OA-E-F-bin}):
\end{minipage}}}\\
\multicolumn{2}{c|}{} &\multicolumn{7}{c|}{\nu}\\
\cline{3-9}
	q 	&t	&2			&3			&4				&5				&6			&7			&8			&\multicolumn{5}{|c}{}\\
\cline{1-9} \multirow{2}{*}{3}	
		&2	&\gr{1/9}	&\gr{1/9}	&\gr{1/9}		&2/27			&3/45*		&3/45*		&7/135		&\multicolumn{5}{|c}{}\\	
		&3	&-			&\gr{1/27}	&\gr{1/27}		&\gr{2/54}		&\gr{2/81}	&			&			&\multicolumn{5}{|c}{}\\	
		&4	&-			&-			&\gr{1/81}		&\gr{1/81}		&4/324		&			&			&\multicolumn{5}{|c}{}\\
\cline{1-9} \multirow{1}{*}{4}	
		&2	&\gr{1/16}	&\gr{1/16}	&\gr{1/16}		&\gr{1/16}		&7/160		&			&			&\multicolumn{5}{|c}{}\\
\cline{1-9}
\multicolumn{14}{c}{}\\
\multicolumn{14}{l}{
\text{\begin{minipage}{0.9\textwidth}
The ratio $R^*/R$ in OAs that maximize $R^*$ among those realizing $F(\nu, q, t)$ 
(for $q =2$, we recall that $F(\nu, 2, t) =E(\nu +1, 2, t)$ if $t$ is even and $E(\nu, 2, t -1)/2$ otherwise \cref{eq-OA-E-F-bin}):
\end{minipage}}}\\
\multicolumn{2}{c|}{} &\multicolumn{6}{c|}{\nu}\\
\cline{3-8}
	q 	&t	&2			&3			&4				&5				&6			&7			&\multicolumn{6}{|c}{}\\
\cline{1-8} \multirow{2}{*}{3}	
		&2	&\gr{1/9}	&\gr{1/9}	&\gr{1/9}		&1/18			&1/18		&1/18		&\multicolumn{6}{|c}{}\\
		&4	&-			&-			&\gr{1/81}		&\gr{1/81}		&2/243		&			&\multicolumn{6}{|c}{}\\
\cline{1-8} \multirow{1}{*}{4}	
		&2	&\gr{1/16}	&\gr{1/16}	&\gr{1/16}		&\gr{1/16}		&1/32		&			&\multicolumn{6}{|c}{}\\
\cline{1-8}
\end{array}$$
\caption{The ratio of the maximum multiplicity $R^*$ of a row to the total number $R$ of rows in OAs and DS's that verify certain optimality conditions. 
Such arrays can be computed by computer solving linear programs (see \cref{sec-cd} for more details). 
The $*$ mark indicates the cases where the corresponding arrays and the proof of their optimality follow from \cite{CSV19}.
We use gray color to indicate cases where a same array realizes both numbers either $E(\nu, q, t)$ and $\rho_E(\nu, q, t)$, or $F(\nu, q, t)$ and $\rho(\nu, q, t)$. 
}\label{tab-rho} 
\end{table}

As noted by Stinson in \cite{S19}, Mukerjee, Qian and Wu in \cite{MQW08} provide an upper bound for $\rho(\nu, q, t)$ for all integers $q\geq 2$, $t\geq 2$, and $\nu\geq t$. In their work, an $OA(R, \nu, q, t)$ is termed {\em nested} if it contains an $OA(R', \nu, q', t)$ as a subarray for some positive integers $R' <R$ and $q' \leq q$. The authors provide a lower bound on the ratio $R/R'$, which generalizes the Rao bound for $R$ in an $OA(R, \nu, q, t)$ \cite{R47}. 
Interpreting $R^*$ identical rows of an $OA(R, \nu, q, t)$ as an $OA(R^*, \nu, 1, t)$, the bound of \cite{MQW08} with $q' =1$ provides an upper bound for $\rho(\nu, q, t)$ \cite{S19}. 
It precisely follows from \cite{MQW08} that $1/\rho(\nu, q, t)$ is at least:
\begin{align}\nonumber 
\begin{array}{rll}
	\sum_{j =0}^{t/2} (q -1)^j \binom{\nu}{j}	
		&\sim \left((q -1)\nu\right)^{\lfloor t/2\rfloor}/\lfloor t/2\rfloor!
			&\text{if $t$ is even;}\\
	\sum_{j =0}^{\lfloor t/2\rfloor} (q -1)^j \binom{\nu}{j} +(q -1)^{\lceil t/2\rceil} \binom{\nu -1}{\lfloor t/2\rfloor}
		&\sim q\times \left((q -1)\nu\right)^{\lfloor t/2\rfloor}/\lfloor t/2\rfloor!
			&\text{if $t$ is odd.}
\end{array}
\end{align}

Thus for $\mathsf{2\,CSP\!-\!q}$ and $\mathsf{3\,CSP\!-\!2}$, the best asymptotic lower bounds we can derive from our approach for the average differential ratio are, respectively, $1/\left((q -1)\nu\right)$ and $1/(2\nu)$.
When $q =2$ and $k\geq 4$, our estimate for $\rho(\nu, 2, k)$ is asymptotically a factor of $1/\left(\lfloor k/2\rfloor!\times 2^{\lceil k/2\rceil}\right)$ of the bound given by \cite{MQW08}. 
 
\smallskip
{\bf{\em Directions for further research.}}
 %
Our estimates of the numbers $\rho(\nu, q, t)$ and $\rho_E(\nu, q, t)$ could potentially be improved for most triples $(\nu, q, t)$. 
First, except for $t =2$, the lower bounds we considered for $\rho(\nu, q, t)$ are derived from simple arrays. By definition of $\rho(\nu, q, t)$ and $F(\nu, q, t)$, the inequality $\rho(\nu, q, t)\geq 1/F(\nu, q, t)$ holds for all sets $(\nu, q, t)$ of parameters. This raises the question of how much better $\rho(\nu, q, t)$ can be compared to $1/F(\nu, q, t)$. \Cref{tab-rho} presents the exact value of $F(\nu, q, t)$ and $\rho(\nu, q, t)$ for some triples $(\nu, q, t)$. For a fair comparison, we report the minimum number of rows in arrays achieving $\rho(\nu, q, t)$ and the maximum number of zero-rows in arrays achieving $F(\nu, q, t)$.
The same table also provides similar information for difference schemes.

Second, we found few results on difference schemes in the literature. Our analysis suggests searching for DS's that maximize the total frequency of the words $\mathbf{a}$, $a\in\Sigma_q$. 
In particular, similar to what Bush did in \cite{B52} for OAs, given two positive integers $q$ and $t$, one should search for the largest $\nu\geq t$ such that $E(\nu, q, t) =q^{t -1}$.
For example, consider the two equations below:
\begin{align}
	y_1 +\ldots +y_\nu 		-y_{\nu +1} -\ldots -y_{2\nu} 					&\equiv 0\bmod{q}	\label{eq-Eq-1}\\
	y_1 +\ldots +y_{\nu -1} +2 y_\nu 	-y_{\nu +1} -\ldots -y_{2\nu +1} 	&\equiv 0\bmod{q}	\label{eq-Eq-2}
\end{align}

Let $P$ be the predicate that accepts solutions to equation \cref{eq-Eq-1}. If we fix in this equation the value of any $2\nu -1$ variables, there is exactly one possible assignment for the remaining variable that  satisfies the equation. Thus $P\in\mathcal{I}_q^{2\nu -1}$, and the accepting entries of $P$ form the rows of an $OA(q^{2\nu -1}, 2\nu, q, 2\nu -1)$.
Furthermore, a vector $y\in\Sigma_q^{2\nu}$ satisfies \cref{eq-Eq-1} {\em if and only if} all vectors of the form $y +\mathbf{a}$ satisfy \cref{eq-Eq-1}. 
We deduce from \cref{it-DS-OA_E} of \cref{pty-OA-DS} that solutions $y$ to equation \cref{eq-Eq-1} that additionally satisfy, for instance, $y_1 =0$ form a $D_{2\nu -1}(q^{2\nu -2}, 2\nu, q)$. 
By a similar argument, the predicate whose accepting entries are the solutions of equation \cref{eq-Eq-2} belongs to $\mathcal{E}_q$, and to $\mathcal{I}_q^{2\nu}$ if $q$ is odd. 
Thus, assuming that $q$ is odd, the solutions $y$ of equation \cref{eq-Eq-2} that additionally satisfy $y_1 =0$ form the rows of a $D_{2\nu}(q^{2\nu -1}, 2\nu +1, q)$. 
Therefore, we have:
\begin{align}
	E(t +1, q, t) &=q^{t -1},	&q, t\in\mathbb{N}\backslash\set{0},\ \text{$t$ or $q$ is odd}
\end{align}

We conclude that if either $k$ or $q$ is odd, the average differential ratio on $(k +1)$-partite instances of $\mathsf{k\,CSP(\mathcal{E}_q)}$ is at least $1/q^{k -1}$.

\section{Reducing the alphabet size}\label{sec-reduc}
In the general case, CSPs become harder as the size of the alphabet increases.
On the one hand, we can reduce to smaller alphabets by increasing the arity of the constraints. Let $p \geq 2$, $q \geq p$, and $k \geq 1$ be three integers, and $\kappa =\lceil\log_p q\rceil$. Then a function $P$ of $k$ variables $x_1 ,\ldots, x_k\in\Sigma_q$ can be interpreted as a function of $k$ strings $y_1 ,\ldots, y_k\in\Sigma_p^\kappa$ where for each $j\in[k]$, $y_j =(y_{j, 1} ,\ldots, y_{j, \kappa})$ is the base $p$ encoding of $x_j$. Thus, an instance of $\mathsf{k\,CSP\!-\!q}$ can be encoded as an instance of $\mathsf{(\lceil \log_p q\rceil k)CSP\!-\!p}$. 
On the other hand, in the most general case, CSPs become harder as the constraint arity increases. In fact, a function of $h$ variables can be interpreted as a function of $k >h$ variables, whose value depends only on its $h$ first inputs (the $k -h$ last inputs are simply ignored).

We investigate whether it is possible to reduce the alphabet size without increasing the constraint arity, possibly at the cost of a reduced approximation guarantee. More formally, given three positive integers $k$, $q$, and $p$ with $p <q$, we seek {\em a reduction} from $\mathsf{k\,CSP\!-\!q}$ to $\mathsf{k\,CSP\!-\!p}$  that {\em preserves the differential approximation ratio}, up to a possible multiplicative factor. 
 %
\label{sec-def-reduc} 
Consider two optimization CSPs $\mathsf{\Pi}$ and $\mathsf{\Pi}'$. A {\em reduction} from $\mathsf{\Pi}$ to $\mathsf{\Pi}'$ can be seen as a polynomial-time algorithm $\mathcal{A}$ for $\mathsf{\Pi}$ that makes call to a hypothetical algorithm $\mathcal{A}'$ for $\mathsf{\Pi}'$ as a subroutine. Such an algorithm is a {\em differential approximation-preserving reduction} ($D$-reduction for short) if there exists $\gamma >0$ such that for any $\rho$-differential approximation algorithm $\mathcal{A}'$ for $\mathsf{\Pi}'$, the induced algorithm $\mathcal{A}$ is a $(\gamma\times\rho)$-differential approximation algorithm for $\mathsf{\Pi}$.  
The quantity $\gamma$ is called {\em the expansion factor} of the reduction. For example, $\mathsf{3\,CSP\!-\!2}$ $D$-reduces to $\mathsf{E2\,Lin\!-\!2}$ with $\gamma =1/2$ \cite{CT18}, while $\mathsf{k\,CSP\!-\!q}$ $D$-reduces to $\mathsf{k\lceil\log_2 q\rceil\,Lin\!-\!2}$ with no loss on the approximation guarantee (i.e., with $\gamma =1$) \cite{CT18-E}.

Most often, the algorithm $\mathcal{A}$ relies on a pair $(f, g)$ of polynomial-time algorithms, where $f$ maps each instance $I$ of $\mathsf{\Pi}$ into an instance $f(I)$ of $\mathsf{\Pi}'$, and $g$ associates with each instance $I$ of $\mathsf{\Pi}$ and each solution $x$ of $f(I)$ a solution $g(I, x)$ of $I$. 
$\mathcal{A}$ then proceeds as follows: starting with $I$, it first computes $f(I)$, then solves $f(I)$ using $\mathcal{A}'$ to obtain a solution $x$, and finally applies $g$ to derive a solution for $I$.
Such a pair $(f, g)$ defines a $D$-reduction with expansion $\gamma$ if, for every instance $I$ of $\mathsf{\Pi}$ and every solution $x$ of $f(I)$, we have the inequality:
$$\begin{array}{rll}
\displaystyle	\frac{v(I, g(I, x)) -\mathrm{wor}(I)}{\mathrm{opt}(I) -\mathrm{wor}(I)}
&\displaystyle	\geq\gamma\times 
				\frac{v(f(I), x) -\mathrm{wor}(f(I))}{\mathrm{opt}(f(I)) -\mathrm{wor}(f(I))}
\end{array}$$

In this section, we present a reduction where the transformation $f$ {does not map an instance $I$ of $\mathsf{\Pi}$ into a single instance of $\mathsf{\Pi}'$, but instead, associates it with multiple instances $f_1(I), \ldots, f_R(I)$. In this case, the solution returned by $\mathcal{A}$ is the best solution it can derives from the solutions $x_1 ,\ldots, x_R$ that algorithm $\mathcal{A}'$ returns on instances $f_1(I) ,\ldots, f_R(I)$.

\subsection{Previous related works and preliminary remarks}\label{sec-reduc-intro}

Consider three natural numbers $q$, $p\in[q]$, and $k\geq 2$, along with an instance $I$ of $\mathsf{k\,CSP\!-\!q}$. Our goal is to derive an approximate solution of $I$ from approximate solutions of instances of $\mathsf{k\,CSP\!-\!p}$. 
Hereafter, we refer to the set of $p$-element subsets of $\Sigma_q$ as $\mathcal{P}_p(\Sigma_q)$. 
For any $S =(S_1 ,\ldots, S_n)\in\mathcal{P}_p(\Sigma_q)^n$, the restriction of $I$ to the solution set $S$, denote by $I(S)$, can be interpreted as an instance of $\mathsf{k\,CSP\!-\!p}$. Specifically, given a family of bijections ($\pi_j: \Sigma_p\rightarrow S_j \,|\, j\in[n])$, we define the instance $f_S(I)$ of $\mathsf{k\,CSP\!-\!p}$ as follows:

\begin{enumerate}
	\item for each variable $x_j$ of $I$, there is in $f_S(I)$ a variable $z_j$ with domain $\Sigma_p$;
	\item for each constraint $C_i =P_i(x_{i_1} ,\ldots, x_{i_{k_i}})$ of $I$, there is in $f_S(I)$ a constraint
			\\$P_i(\pi_{i_1}(z_{i_1}) ,\ldots, \pi_{i_{k_i}}(z_{i_{k_i}}))$ with the same associated weight $w_i$ as $C_i$.
\end{enumerate} 

By construction, for any $z\in\Sigma_p^n$, the solution $g_S(I, z):= (\pi_1(z_1) ,\ldots, \pi_n(z_n))$
of $I(S)$ achieves on $I(S)$ --- and thus, on $I$ --- the same objective value as the solution $z$ on $f_S(I)$. This implies that the worst solution value on $f_S(I)$ is at least as good as the worst solution value on $I$. 

The consideration of sub-instances $I(S)$ seems natural having in view a reduction from $\mathsf{k\,CSP\!-\!q}$ to $\mathsf{k\,CSP\!-\!p}$. In fact, if the terms $w_i\times P_i$ of the objective function are non-negative, there exists a standard approximation-preserving reduction from $\mathsf{Max\,Ek\,CSP\!-\!q}$ to $\mathsf{Max\,Ek\,CSP\!-\!p}$ that consists precisely in choosing a solution subset $S$ by randomly selecting $n$ subsets $S_1 ,\ldots, S_n$ independently and uniformly over $\mathcal{P}_p(\Sigma_q)$ \cite{CMM09}. The argument is based on the fact that, in expectation, the optimal value of $I(S)$ remains a constant fraction of the optimal value. Specifically, we have \cite{CMM09}:
\begin{align}\label{eq-opt-class}\textstyle
	\mathbb{E}_S[\mathrm{opt}(I(S))]	&\geq p^k/q^k \times \mathrm{opt}(I)
\end{align}

Suppose that for each $S\in\mathcal{P}_p(\Sigma_q)^n$, we can compute in polynomial time a solution $x(S)$ that is $\rho$-standard approximate on $I(S)$. In expectation, these solutions are $(p^k/q^k \times\rho)$-standard approximate on $I$, as shown below:
$$\begin{array}{rll}
\mathbb{E}_S[v(I, x(S))] &\geq \mathbb{E}_S[\rho\times\mathrm{opt}(I(S))] 
							&\text{by assumption on $x(S)$, $S\in\mathcal{P}_p(\Sigma_q)^n$}\\
	&\geq\rho\times p^k/q^k \times \mathrm{opt}(I)	&\text{by \cref{eq-opt-class}}
\end{array}$$
 %
The reduction can then be derandomized using an alternative distribution over $\mathcal{P}_p(\Sigma_q)^n$, at the cost of a multiplicative factor of $(1 -\varepsilon)$ on the approximation guarantee \cite{CMM09}. 

To establish \cref{eq-opt-class}, the authors of \cite{CMM09} associate a presumed optimal solution $x^*$ with a family $(x^*(S)\,|\,S\in\mathcal{P}_p(\Sigma_q)^n)$ of solutions where, for each $S\in\mathcal{P}_p(\Sigma_q)^n$, $x^*(S)$ can be any solution that coincides with $x^*$ at its coordinates indexed by $j$ such that $x^*_j\in S_j$. 
They observe that a constraint $P_i(x_{J_i})$ evaluates the same on $x^*(S)$ as on $x^*$ if, for every $j\in J_i$, $x^*_j\in S_j$, which occurs with probability: 
$$\begin{array}{rl}\prod_{j =1}^k \binom{q -1}{p -1}/\binom{q}{p} &=(p/q)^k\end{array}$$

They deduce that the expected value of $v(I, x^*(S))$ over $\mathcal{P}_p(\Sigma_q)^n$ satisfies:
\begin{align}\label{eq-x*(S)}
&\mathbb{E}_S[v(I, x^*(S))]\textstyle\, 
	=\sum_{i =1}^m w_i\times\mathbb{E}_S\left[P_i(x^*(S)_{J_i})\right]\nonumber\\
&\textstyle\, 
	=\sum_{i =1}^m w_i\times\left(
		\mathbb{P}_S\left[x^*_{J_i}\in S_{J_i}\right]		\times P_i(x^*_{J_i})
		+\mathbb{P}_S\left[x^*_{J_i}\notin S_{J_i}\right]	\times 
			\mathbb{E}_S\left[P_i(x^*(S)_{J_i})\,|\,x^*_{J_i}\notin S_{J_i}\right]
	\right)\nonumber\\
&\textstyle\ 
	\geq p^k/q^k\times\mathrm{opt}(I) +(1 -p^k/q^k)\sum_{i =1}^m w_i\times
					\mathbb{E}_S\left[P_i(x^*(S)_{J_i})\,|\,x^*_{J_i}\notin S_{J_i}\right]
\end{align}

Since $\mathrm{opt}(I(S))\geq v(I, x^*(S))$, $S\in\mathcal{P}_p(\Sigma_q)^n$, inequality  \cref{eq-opt-class} follows directly from \cref{eq-x*(S)} under the assumption that for all $i\in[m]$, $w_i\times P_i$ is non-negative. In contrast, we do not know how to compare the quantity 
$\sum_{i =1}^m w_i\times\mathbb{E}_S\left[P_i(x^*(S)_{J_i})\,|\,x^*_{J_i}\notin S_{J_i}\right]$ to a solution value. Thus, we cannot derive from inequality \cref{eq-x*(S)} a lower bound on the expected differential ratio reached at $\mathrm{opt}(I(S))$.} 

Therefore, we seek an alternative connection between the optimal values on sub-instances $I(S)$ and the optimal value on $I$, which would allow us to derive differential approximate solutions on $I$ from differential approximate solutions of sub-instances $I(S)$.
To simplify our analysis, we restrict our focus to solution sets of the form $T^n$, where $T$ is a subset of $\Sigma_q$ with cardinality $p$. 
Identifying $T^n$ with $T$, we henceforth denote by $I(T)$ the restriction of $I$ to the solution set $T^n$. Note that choosing a best solution among a set of hypothetical approximate solutions $x(T)$ of the subinstances $I(T)$ requires comparing only $\binom{q}{p}$ solution values. While this number can be large, it remains constant. 

Similar to our approach for estimating the average differential ratio, we associate with each solution $x$ of a given instance $I$ two multisets of solutions, $\mathcal{X}(I, x)$ and $\mathcal{Y}(I, x)$, both of equal size $R$. Here, $\mathcal{X}(I, x)$ is a subset of $\set{T^n\,|\,T\in\mathcal{P}_p(\Sigma_q)}$, 
	$\mathcal{Y}(I, x)$ contains a certain number $R^* >0$ of occurrences of $x$, 
	and the sum of the solution values over $\mathcal{X}(I, x)$ is the same as over $\mathcal{Y}(I, x)$, i.e.:
$$\begin{array}{rl}
	\sum_{y\in\mathcal{X}(I, x)} v(I, y)	&=\sum_{y\in\mathcal{Y}(I, x)} v(I, y)
\end{array}$$

Taking $\mathcal{X}(I, .)$ and $\mathcal{Y}(I, .)$ at an optimal solution $x^*$, we deduce that a best solution over $\mathcal{X}(I, x^*)$, and hence a best solution over $\{T^n\,|\,T\in\mathcal{P}_p(\Sigma_q)\}$, achieves a differential ratio of at least $R^*/R$. 

\subsection{Partition-based solution multisets}\label{sec-reduc-pbsm}

We describe how we construct our solution multisets $\mathcal{X}(I, .)$ and $\mathcal{Y}(I, .)$ for an instance $I$ of $\mathsf{k\,CSP\!-\!q}$.

\smallskip
$\bullet$ {\bf Solution multisets association.} 
Each solution $x$ induces a partition of the index set $[n]$ into $q$ subsets, based on the values taken by its coordinates. Given a solution $x\in\Sigma_q^n$ and a vector $u\in\Sigma_q^q$, we define the solution $y(x, u)$ obtained from $x$ by assigning, for each $c\in\Sigma_q$, the value $u_c$ to the coordinates equal to $c$ in $x$. Formally, for each $j\in[n]$, $y(x, u)_j =u_{x_j}$. 
In particular, if $u =(0, 1 ,\ldots, q -1)$, then $y(x, u) =x$.
Furthermore, the components of a solution $y(x, u)$ take as many different values as those of $u$.

We then consider two arrays $\Psi$ and $\Phi$ on $\Sigma_q$, each containing $q$ columns and having the same number of rows, denoted by $R$. For $x\in\Sigma_q^n$, we associate with $(\Psi, \Phi)$ the solution multisets $\mathcal{X}(I, x)$ and $\mathcal{Y}(I, x)$ defined by:
\begin{align}\nonumber 
	\begin{array}{rlcrlc}
		  \mathcal{X}(I, x)	&=\left(y(x, \Psi_r)\,|\,r\in[R]\right),
		&&\mathcal{Y}(I, x)	&=\left(y(x, \Phi_r)\,|\,r\in[R]\right)
	\end{array}
\end{align}

$\bullet$ {\bf Conditions.}
We denote by $R^*$ the number of rows of the form $(0, 1 ,\ldots, q -1)$ in $\Phi$. 
The arrays $\Psi$ and $\Phi$ must satisfy certain conditions. 
First, each row of $\Psi$ must contain at most $p$ distinct values. 
This ensures that the solutions in $\mathcal{X}(I, x)$ are each restricted to the domain of some sub-instance $I(T)$. 
Then, we require that $R^*>0$ to ensure that the solution values $v(I, y)$, $y\in \mathcal{Y}(I, x)$ cover $\mathrm{opt}(I)$ when $x$ is optimal.
Finally, since our ultimate goal is to connect $\mathrm{opt}(I)$ to solution values on the $I(T)$ sub-instances, the sum of the solution values over $\mathcal{X}(I, x)$ and $\mathcal{Y}(I, x)$ must be the same, which means that $\Psi$ and $\Phi$ should verify:
\begin{align}\label{eq-vq-p--P=Q}
	\textstyle\sum_{r =1}^R v(I, y(x, \Psi_r))	&\textstyle=\sum_{r =1}^R v(I, y(x, \Phi_r)),	&x\in\Sigma_q^n
\end{align}

Suppose that $\Psi$ and $\Phi$ fulfill the above requirements, and let $x^*$ be an optimal solution on $I$. We assume {\em w.l.o.g.} that $I$ is a maximization instance. 
Furthermore, we denote by $\mathcal{T}$ a subset of $\mathcal{P}_p(\Sigma_q)$ such that each row of $\Psi$ contains only values from some $T\in\mathcal{T}$. Then, we have:
\begin{align}
\textstyle\nonumber
\max_{T\in\mathcal{T}} \left\{\mathrm{opt}\left(I(T)\right)\right\}
	&\textstyle\nonumber \geq	\sum_{r =1}^R v\left(I, y(x^*,\Psi_r)\right)/R
									\qquad\text{by assumption on $\mathcal{T}$}\\
	&\textstyle\nonumber =	 	\sum_{r =1}^R v\left(I, y(x^*,\Phi_r)\right)/R	
									\qquad\text{by \cref{eq-vq-p--P=Q}}\\
	&\textstyle\label{eq-opt-Psi}	\geq	\left(R^*\times\mathrm{opt}(I) + (R -R^*)\mathrm{wor}(I)\right)/R
\end{align}

Thus, the best solutions over $\cup_{T\in\mathcal{T}} T^n$ are $R^*/R$-differential approximate on $I$. Now suppose that for each $T\in\mathcal{T}$, we can compute a $\rho$-differential approximate solution $x(T)$ on $I(T)$. Then we have the following bound on the best objective value performed by these solutions:
\begin{align}
\textstyle\nonumber
\max_{T\in\mathcal{T}} \left\{v\left(I, x(I(T))\right)\right\}
	&\textstyle\nonumber
	\geq	\max_{T\in\mathcal{T}} \left\{\rho\,\mathrm{opt}(I(T)) +(1 -\rho)\mathrm{wor}(I(T))\right\}\\
	&\textstyle\label{eq-opt-T}
	\geq	\rho\,\max_{T\in\mathcal{T}} \left\{\mathrm{opt}(I(T))\right\} +(1 -\rho)\mathrm{wor}(I)
\end{align}

Combining \cref{eq-opt-Psi} and \cref{eq-opt-T}, we get:
\begin{align}
\textstyle\nonumber
\max_{T\in\mathcal{T}} \left\{v\left(I, x(I(T))\right)\right\}	&\textstyle\geq\rho\left(
		 R^*/R\times\mathrm{opt}(I) +(1 -R^*/R)\mathrm{wor}(I)
	\right) +(1 -\rho)\mathrm{wor}(I)\\	\Rightarrow\ 
\textstyle\label{eq-apx-T}
\max_{T\in\mathcal{T}} \left\{v\left(I, x(I(T))\right)\right\}	&\textstyle\geq
	\rho\,R^*/R \times \mathrm{opt}(I) +(1 -\rho\,R^*/R)\mathrm{wor}(I)
\end{align}

Thus, the solutions $x(T)$ that achieve the best objective value over $\mathcal{T}$ are $\rho\times R^*/R$-differential approximate on $I$.

\subsection{Connection to combinatorial designs}\label{sec-reduc-CD}  

We need to identify the conditions under which $\Psi$ and $\Phi$ verify \cref{eq-vq-p--P=Q}. 
Clearly, a sufficient condition for the mean of the solution values to be the same over $(y(x, \Psi_r))\,|\,r\in[R])$ and $(y(x, \Phi_r))\,|\,r\in[R])$ is that each constraint of $I$ takes the same value on average over these two solution multisets. 
In other words, $\Psi$ and $\Phi$ verify \cref{eq-vq-p--P=Q} if they satisfy:
\begin{align}\label{eq-vq-p--P=Q-i}\textstyle
	\sum_{r =1}^R P_i(y(x, \Psi_r)_{J_i})	&\textstyle=\sum_{r =1}^R P_i(y(x, \Phi_r)_{J_i}), &i\in[m],\ x\in\Sigma_q^n
\end{align}
	
Consider a solution $x\in\Sigma_q^n$ and a constraint $C_i =P_i(x_{J_i})$ of $I$. 
The coordinates $x_j$, $j\in J_i$ of $x$ take at most $\min\set{q, k}$ pairwise distinct values. 
Let $H =(c_1 ,\ldots, c_t)$ be the subsequence of $(0, 1 ,\ldots, q -1)$ induced by these values. A sufficient condition for $(\Psi, \Phi)$ to satisfy \cref{eq-vq-p--P=Q-i} at $(i, x)$ is that the function $P_i$ is evaluated on the same multisets of entries over both $\left(y(x, \Psi_r)\,|\,r\in[R]\right)$ and $\left(y(x, \Phi_r)\,|\,r\in[R]\right)$.
By the definition of solutions $y(x, u)$, this occurs {\em if and only if} $(\Psi^H_r\,|\,r\in[R])$ and $(\Phi^H_r\,|\,r\in[R])$ reduce to the same multisubset of $\Sigma_q^t$. We observe that this cannot happen unless $|H|\leq p$, since $\Phi^H_r =(c_1 ,\ldots, c_t)$ must hold for {\em at least one} row of $\Phi$, while $\card{\set{\Psi^{c_1}_r ,\ldots, \Psi^{c_t}_r}}\leq p$ must hold for {\em all} rows of $\Psi$. Considering that $\card{H}\leq\card{J_i}\leq k$, we conclude that $(\Psi,\Phi)$ satisfies \cref{eq-vq-p--P=Q-i} in case where $k\leq p\leq q$ and the arrays $\Psi$ and $\Phi$ verify for all $k$-length subsequences $H$ of $\Sigma_q$ the property that their subarrays $\Psi^H$ and $\Phi^H$ coincide (up to row ordering).

From now on, we assume that $k\leq p\leq q$. 
Consider the function $\mu^\Psi -\mu^\Phi:\Sigma_q^q\rightarrow[-1, 1]$. Its mean is clearly zero. Furthermore, for a $t$-length subsequence $H$ of $\Sigma_q$, the subarrays $\Psi^H$ and $\Phi^H$ consist of the same collection of words {\em if and only if}, for all $v\in\Sigma_q^t$, the total frequency of words $u =(u_0, u_1 ,\ldots, u_{q -1})\in\Sigma_q^q$ such that $u_H =v$ is the same in $\Psi$ as in $\Phi$. Therefore, the subarrays $\Psi^H$ and $\Phi^H$ coincide for any subsequence $H$ of length at most $k$ of $\Sigma_q$ {\em if and only if} $\mu^\Psi -\mu^\Phi$ is balanced $k$-wise independent. 
These considerations lead us to introduce the following family of array pairs:

\begin{definition}\label{def-reduc-CD}
Let $k >0$, $p\geq k$, and $q\geq p$ be three integers. 
Two arrays $\Psi$ and $\Phi$ with $q$ columns and entries from $\Sigma_q$ form a {\em $(q, p)$-alphabet reduction pair of arrays} (for short, a $(q, p)$-ARPA) {\em of strength $k$} if they satisfy the conditions below:
\begin{enumerate}
	\item $\Phi$ contains at least 1 row of the form $(0, 1, \ldots, q -1)$;
	\item \label{it-thm-q-p-P} the components of each row of $\Psi$ take at most $p$ different values;
	\item \label{it-thm-q-p-P=Q} the function $\mu^\Psi -\mu^\Phi$ is balanced $k$-wise independent.
\end{enumerate} 

For two integers $R^* >0$ and $R\geq R^*$, we denote by $\Gamma(R, R^*, q, p, k)$ the (possibly empty) set of $(q, p)$-ARPAs of strength $k$ in which the row $(0, 1, \ldots, q -1)$ has multiplicity $R^*$ and the arrays have $R$ rows each.

Furthermore, we define $\gamma(q, p, k)$ as the largest number $\gamma$ for which $\Gamma(R, R^*, q, p, k)\neq\emptyset$ holds for at least two integers $R^* >0$ and $R\geq R^*$ such that $R^*/R =\gamma$.
\end{definition}

\Cref{tab-Gamma-ex} shows pairs of arrays that realize $\gamma(q, p, k)$. The preceding discussion establishes the following connection between ARPAs and the reducibility of $\mathsf{k\,CSP\!-\!q}$ to $\mathsf{k\,CSP\!-\!p}$:

\begin{table}[t]\begin{center}{\footnotesize
\setlength\arraycolsep{3pt}
$\begin{array}{cp{8pt}c}
\begin{array}{cp{4pt}c}
	\multicolumn{3}{c}{\gamma(4, 3, 2) =1/3}\\[4pt]	
	\begin{array}{cccc}
		\Psi^0	&\Psi^1	&\Psi^2	&\Psi^3\\\hline
		0		&0		&2		&3	\\
		0		&1		&0		&3	\\
		0		&1		&2		&2	\\
		0		&1		&2		&2	\\
		3		&0		&0		&2	\\
		3		&1		&2		&3	\\
	\end{array}&&\begin{array}{cccc}
		\Phi^0	&\Phi^1	&\Phi^2	&\Phi^3\\\hline
		0		&0		&2		&2	\\
		0		&1		&0		&2	\\
		0		&1		&2		&3	\\
		0		&1		&2		&3	\\
		3		&0		&0		&3	\\
		3		&1		&2		&2	\\
	\end{array}\\\\[-4pt]
	\multicolumn{3}{c}{\gamma(5, 3, 2) = 1/6}\\[4pt]
	\begin{array}{ccccc}
		\Psi^0	&\Psi^1	&\Psi^2	&\Psi^3	&\Psi^4\\\hline
		0		&1		&3		&3		&3	\\
		0		&2		&2		&2		&4	\\
		1		&1		&2		&1		&4	\\
		1		&2		&3		&1		&3	\\
		3		&3		&2		&3		&4	\\
		3		&3		&3		&2		&3	\\
	\end{array}&&\begin{array}{cccccl}
		\Phi^0	&\Phi^1	&\Phi^2	&\Phi^3	&\Phi^4\\\hline
		0		&1		&2		&3		&4	\\
		0		&2		&3		&2		&3	\\
		1		&1		&3		&1		&3	\\
		1		&2		&2		&1		&4	\\
		3		&3		&2		&2		&4	\\
		3		&3		&3		&3		&3	\\
	\end{array}
\end{array}&&
\begin{array}{cp{4pt}c}
	\multicolumn{3}{c}{\gamma(5, 4, 3) =1/5}\\[4pt]	
	\begin{array}{ccccc}
		\Psi^0	&\Psi^1	&\Psi^2	&\Psi^3	&\Psi^4\\\hline
		0	&0	&1	&0	&3	\\
		0	&0	&2	&3	&4	\\
		0	&0	&2	&3	&4	\\
		0	&1	&1	&3	&4	\\
		0	&1	&1	&3	&4	\\
		0	&1	&2	&0	&4	\\
		0	&1	&2	&0	&4	\\
		0	&1	&2	&3	&3	\\
		0	&1	&2	&3	&3	\\
		4	&0	&1	&0	&4	\\
		4	&0	&1	&3	&3	\\
		4	&0	&2	&0	&3	\\
		4	&1	&1	&0	&3	\\
		4	&1	&2	&3	&4	\\
		4	&1	&2	&3	&4	\\
	\end{array}&&\begin{array}{cccccl}
		\Phi^0	&\Phi^1	&\Phi^2	&\Phi^3	&\Phi^4	\\\hline
		0	&0	&1	&3	&4	\\
		0	&0	&2	&0	&4	\\
		0	&0	&2	&3	&3	\\
		0	&1	&1	&0	&4	\\
		0	&1	&1	&3	&3	\\
		0	&1	&2	&0	&3	\\
		0	&1	&2	&3	&4	\\
		0	&1	&2	&3	&4	\\
		0	&1	&2	&3	&4	\\
		4	&0	&1	&0	&3	\\
		4	&0	&1	&0	&3	\\
		4	&0	&2	&3	&4	\\
		4	&1	&1	&3	&4	\\
		4	&1	&2	&0	&4	\\
		4	&1	&2	&3	&3	\\
	\end{array}
\end{array}
\end{array}$
\caption{\label{tab-Gamma-ex}
ARPAs realizing $\gamma(4, 3, 2)$, $\gamma(5, 3, 2)$, and $\gamma(5, 4, 3)$.
These array pairs have been computed by computer (see \cref{sec-cd}).
}}\end{center}\end{table}

\begin{theorem}\label{thm-reduc-CD}
For all constant integers $k >0$, $p\geq k$, and $q\geq p$, $\mathsf{k\,CSP\!-\!q}$ $D$-reduces to $\mathsf{k\,CSP\!-\!p}$ with an expansion factor of $\gamma(q, p, k)$ on the approximation guarantee by a reduction that involves solving $\binom{q}{p}$ instances of $\mathsf{k\,CSP\!-\!p}$.

Moreover, for all integers $R^* >0$ and $R\geq R^*$ such that an ARPA in $\Gamma(R, R^*, q, p, k)$ is known, $\mathsf{k\,CSP\!-\!q}$ $D$-reduces to $\mathsf{k\,CSP\!-\!p}$ with an expansion factor of $R/R^*$ on the approximation guarantee. 
The reduction requires solving at most $\min\set{R, \binom{q}{p}}$ instances of $\mathsf{k\,CSP\!-\!p}$.
\end{theorem}

\begin{proof}
We consider an instance $I$ of $\mathsf{k\,CSP\!-\!q}$ together with an ARPA $(\Psi,\Phi)$ of $\Gamma(R, R^*, q, p, k)$, where $R^*>0$ and $R\geq R^*$. We denote by $\mathcal{T}$ a subset of $\mathcal{P}_p(\Sigma_q)$ such that, for every $r\in[R]$, the set $\set{\Psi_r^0 ,\ldots, \Psi_r^{q -1}}$ is contained in some $T\in\mathcal{T}$. 

We have already argued that $(\Psi,\Phi)$ verifies \cref{eq-vq-p--P=Q-i} and hence \cref{eq-vq-p--P=Q}. Now assume that there exists a $\rho$-differential approximation algorithm $\mathcal{A}$ for $\mathsf{k\,CSP\!-\!p}$. We can use $\mathcal{A}$ to compute, for each $T\in\mathcal{T}$, a $\rho$-differential approximate solution $x(T)$ of $I(T)$. If the goal is to maximize on $I$, then we have also argued that these solutions and $(\Psi,\Phi)$ verify \cref{eq-apx-T} (and otherwise, they verify the reverse inequality).

If $(\Psi,\Phi)$ is explicitly known, then we can derive from $\Psi$ a subset $\mathcal{T}$ of $\mathcal{P}_p(\Sigma_q)$ of cardinality at most $R$. Besides, regardless of whether a concrete ARPA, either in $\Gamma(R, R^*, q, p, k)$ or realizing $\gamma(q,p,k)$, is known, we can always set $\mathcal{T}$ to $\mathcal{P}_p(\Sigma_q)$.
\end{proof}

In addition, the previous discussion establishes the following relationship between the optimal values on the sub-instances $I(T)$ and the optimal value on $I$: 
\begin{theorem}\label{thm-gap}
For all fixed integers $k\geq 2$, $p\geq k$, and $q\geq p$, and all instances $I$ of $\mathsf{k\,CSP\!-\!q}$, the best solutions among those whose components take at most $p$ different values are $\gamma(q, p, k)$-approximate. Formally, let $\mathrm{opt}(I\,|\,\mathcal{P}_p(\Sigma_q))$ refer to the quantity:
$$\begin{array}{ll}
	\max_{T\subseteq\Sigma_q: |T| =p} \left\{\max_{x\in T^n} v(I, x)\right\}
		&\text{if the goal on $I$ is to maximize;}\\
	\min_{T\subseteq\Sigma_q: |T| =p} \left\{\min_{x\in T^n} v(I, x)\right\}
		&\text{if the goal on $I$ is to minimize.}
\end{array}$$

Then we have the following bound on the differential ratio reached at   $\mathrm{opt}(I\,|\,\mathcal{P}_p(\Sigma_q))$:
\begin{align}\label{eq-reduc-opt}
\frac{\mathrm{opt}(I\,|\,\mathcal{P}_p(\Sigma_q)) -\mathrm{wor}(I)}{\mathrm{opt}(I) -\mathrm{wor}(I)}
	&\geq\gamma(q, p, k)
\end{align}
\end{theorem}
\begin{proof}
Assuming {\em w.l.o.g.} that the goal on $I$ is to maximize, the inequality \cref{eq-opt-Psi} holds for all $(q, p)$-ARPAs of strength $k$.
\end{proof}

\subsection{A lower bound for $\gamma(q, p, k)$ and derived approximation results}
\label{sec-reduc-gamma}

Given three integers $k >0$, $p\geq k$ and $q\geq p$, our goal is to construct $(q, p)$-ARPAs of strength $k$ that maximize the ratio $R^*/R$.
If $p =q$, this ratio is trivially 1. This follows from considering for both $\Psi$ and $\Phi$ a same array consisting of identical rows of the form $(0, 1 ,\ldots, q -1)$.
\begin{property}\label{pty-gamma-base}	
For all positive integers $k$, $q\geq k$, and $R^*$, $\Gamma(R^*, R^*, q, q, k)\neq\emptyset$.
\end{property}

If $p >k$, then we observe that extending each row of a $(q -p +k, k)$-ARPA of strength $k$ by 
$(q -p +k, q -p +k +1, \ldots, q -1)$ yields a $(q, p)$-ARPA of the same strength $k$. 
\begin{property}\label{pty-gamma-p>k}	
For all positive integers $k$, $p >k$, $q \geq p$, $R^*$, and $R\geq R^*$, 
	if $\Gamma(R, R^*, q -p +k, k, k)\neq\emptyset$, then $\Gamma(R, R^*, q, p, k)\neq\emptyset$.
\end{property}

Now assume $p =k <q$. Before we go any further, we introduce a quantity that is involved in subsequent arguments.
\begin{property}\label{pty-T}
Given two natural numbers $a$ and $b$ such that $a >b$, we define:
\begin{align}\label{eq-T-def}
T(a, b)	&\textstyle :=\sum_{r =0}^b	\binom{a}{r} \binom{a -1 -r}{b -r}
\end{align}

These numbers satisfy the following relation:
\begin{align}
\label{eq-T_base}		
	T(b +1, b) 		&= \textstyle 2^{b +1} -1,							&b\in\mathbb{N}							\end{align}

They also satisfy the following recurrence relations:
\begin{align}
\label{eq-T_a+b+}		
	T(a, b) 		&= \textstyle 2^b \binom{a -1}{b} +T(a -1, b -1),	&a, b\in\mathbb{N},\ a >b >0	\\
\label{eq-T_b+}			
	T(a, b) 		&= \textstyle 2^b \binom{a}{b} -T(a, b -1),			&a, b\in\mathbb{N},\ a >b >0	\\			\label{eq-T_a+}			
	T(a, b)			&= \textstyle 2^{b +1} \binom{a -1}{b} -T(a -1, b),	&a, b\in\mathbb{N},\ a >b +1 	\\
\label{eq-T_a+b+bis}	
	T(a, b) 		&= \textstyle 2 T(a -1, b -1) +T(a -1, b),			&a, b\in\mathbb{N},\ a >b +1,\,b >0		\end{align}
\end{property}

\begin{proof}[Proof (sketch)]
For \cref{eq-T_base}: according to \cref{eq-T-def}, $T(b +1, b)$ is equal to $\sum_{r =0}^b\binom{b +1}{r}$. 
The recursions \cref{eq-T_a+b+,eq-T_b+} are obtained by applying Pascal's rule to coefficients of the form $\binom{a}{r}$ and $\binom{a -1 -r}{b -r}$, respectively. 
We then derive the identity \cref{eq-T_a+} by subtracting \cref{eq-T_a+b+} from \cref{eq-T_b+}, both evaluated at $(a, b +1)$.
Finally, the identity \cref{eq-T_a+b+bis} is $2\times\cref{eq-T_a+b+} -\cref{eq-T_a+}$.
\end{proof}

\begin{algorithm}
\begin{algorithmic}[1]
\State\label{alg-Colq}Duplicate the first column of each array $\Psi$ and $\Phi$ into a $q$-th column
\State\label{alg-R*}Substitute for each row of the form $(0, 1 ,\ldots, q -2, 0)$ in $\Phi$ the row $(0, 1 ,\ldots, q -2, q -1)$
\For{$h =k -1$ {\bf down to} 0}\label{alg-For_h}
	\ForAll{$J\subseteq\Sigma_{q -1}$ with $|J| =h$}\label{alg-For_J}
		\Comment{$\alpha(J)$ is the word of $\Sigma_q^{q -1}$ defined by \cref{eq-def-alpha}}
		\If{$h\equiv k -1\bmod{2}$}
			\State Insert $\binom{q -h -2}{k -h -1}\times R^*$ rows 
					of the form $\left(\alpha(J), q -1\right)$ in $\Psi$
			\State Insert $\binom{q -h -2}{k -h -1}\times R^*$ rows 
					of the form $\left(\alpha(J), 0\right)$ in $\Phi$
		\Else
			\State Insert $\binom{q -h -2}{k -h -1}\times R^*$ rows 
					of the form $\left(\alpha(J), q -1\right)$ in $\Phi$
			\State Insert $\binom{q -h -2}{k -h -1}\times R^*$ rows 
					of the form $\left(\alpha(J), 0\right)$ in $\Psi$
		\EndIf
	\EndFor
\EndFor
\end{algorithmic}
\caption{Transforming an ARPA $(\Psi,\Phi)\in\Gamma(R, R^*, q -1, k, k)$ into an ARPA of $\Gamma\left(R', R^*, q, k, k\right)$ for $R' =R +R^* \times T(q -1, k -1)$.
\label{algo-Gamma_rec}}
\end{algorithm}

We now show how to derive $(q, k)$-ARPAs of strength $k$ from $(q -1, k)$-ARPAs of the same strength $k$.  

\begin{lemma}\label{lem-gamma-p=k}\label{def-T}
Let $k$, $q >k$, $R^*$, and $R \geq R^*$ be four positive integers, and $R' :=R +R^*\times T(q -1, k -1)$, where $T(q -1, k -1)$ is defined by \cref{eq-T-def}. Then, if there exists an ARPA in $\Gamma(R, R^*, q -1, k, k)$, then there also exists a corresponding ARPA in $\Gamma\left(R', R^*, q, k, k\right)$.
\end{lemma}

\begin{proof}
Assume there exists $(\Psi,\Phi)\in\Gamma(R, R^*, q -1, k, k)$. We use \cref{algo-Gamma_rec} to transform $(\Psi,\Phi)$ into an ARPA of $\Gamma(R', R^*, q, k, k)$. 
\Cref{tab-Gamma-rec} shows the construction starting from the basic family $\Gamma(1, 1, k, k, k)$ when either $k =2$ and $q \in\{3, 4, 5, 6\}$, or $k =3$ and $q\in\{4, 5\}$.

In this algorithm, for $J\subseteq\Sigma_{q -1}$, $\alpha(J) =(\alpha(J)_0 ,\ldots, \alpha(J)_{q -2})$ refers to the word of $\Sigma_q^{q -1}$ defined by:
\begin{align}\label{eq-def-alpha}
	\alpha(J)_j	&:=\left\{\begin{array}{rl}
									j		&\text{if $j\in J$}\\
									q -1	&\text{otherwise}
								\end{array}\right.,		&j\in\Sigma_{q -1}
\end{align}

We argue why $(\Psi, \Phi)$, at the end of \cref{algo-Gamma_rec}, is an element of $\Gamma(R', R^*, q, k, k)$. First, for a natural number $h$, the number of $h$-cardinality subsets of $\Sigma_{q -1}$ is equal to $\binom{q -1}{h}$. 
The construction therefore inserts into each of the two arrays a number of new rows equal to:
$$\begin{array}{rll}
R^*\times\sum_{h =0}^{k -1}\binom{q -1}{h}\times\binom{q -h -2}{k -h -1}
	&=R^*\times T(q -1, k -1) 
	&=R' -R 
\end{array}$$

Second, it is clear that the array $\Phi$ contains exactly $R^*$ rows of the form $(0, 1 ,\ldots, q -1)$. 
Moreover, each row $\Psi_r$ (where $r\in[R']$) of $\Psi$ contains at most $k$ distinct values. 
For rows of index $r\leq R$, this follows from the initial definition of $\Psi$ and the way its $q$-th column is initialized. 
Rows of larger index are of the form either $(\alpha(q, J), q -1)$ with $|J|\leq k -1$ or $(\alpha(q, J), 0)$ with $|J|\leq k -2$, implying exactly $|J| +1\leq k$ and at most $|J| +2\leq k$ different values, respectively. 

It remains to show that the difference $\mu_\Psi -\mu_\Phi$ in frequency of words appearing in $\Psi$ and $\Phi$ is balanced $k$-wise independent. 
Formally, if $R(\Psi,\Phi)$ denotes the number of rows in $\Psi$ and $\Phi$, then we must prove that, at the end of \cref{algo-Gamma_rec}, $\Psi$ and $\Phi$ satisfy:
\begin{align}\label{eq-qkk}
\begin{array}{rl}
\card{\set{r\in[R(\Psi,\Phi)]\,|\,\Psi_r^J =v}} 
	-\card{\set{r\in[R(\Psi,\Phi)]\,|\,\Phi_r^J =v}}	&=0,
\ \ \ \ \\\multicolumn{2}{r}{
	J =(j_1 ,\ldots, j_k)\in\Sigma_q^k,\ j_1 <\ldots< j_k,\ v\in\Sigma_q^k}
\end{array}
\end{align}

The proof is given in detail in \cref{sec-qkk}. Here we just state its principles. 
Since $(\Psi,\Phi)$ is initially a $(q -1, k)$-ARPA of strength $k$, it satisfies \cref{eq-qkk} (with $R(\Psi,\Phi) =R$) before proceeding to \cref{alg-R*} of the algorithm. Executing this step causes $(\Psi, \Phi)$ to temporarily violate \cref{eq-qkk} at pairs $(J, v)$ such that $j_k =q -1$ and $v$ is the word either $(j_1 ,\ldots, j_{k -1}, 0)$ or $(j_1 ,\ldots, j_{k -1}, q -1)$. 
The first iteration of the outer {\bf for} loop corrects these violations of \cref{eq-qkk}. 
However, it also introduces new violations of \cref{eq-qkk} at pairs $(J, v)$ such that $j_k =q -1$, $v_k\in\set{0, q -1}$, 
	$(v_1 ,\ldots, v_{k -1})\in\set{j_1, q -1}\times\ldots\times\set{j_{k -1}, q -1}$, 
and $v_s =j_s$ holds for at most $k -2$ integers $s\in[k -1]$. 
Subsequent iterations of the outer {\bf for} loop, where $h <k -1$, iteratively correct these violations of \cref{eq-qkk} at those of these pairs that satisfy $v_s =j_s$ for {\em exactly} $h$ integers $s\in[k -1]$. 
\end{proof}

Repeatedly applying the construction of \cref{lem-gamma-p=k} starting from the basic ARPA (where $\Psi =\Phi =\left\{(0, 1 ,\ldots, k -1)\right\}$) yields an ARPA of $\Gamma(R, 1, q, k, k)$, for $R =1 +T(k, k -1) +\ldots+ T(q -1, k -1)$. We formalize this construction in \cref{algo-Gamma}, and show the resulting array pairs for $(k, q)\in\set{(2, 6), (3, 5)}$ in \cref{tab-Gamma-rec}. 
Thus, we obtain a lower bound on the number $\gamma(q, p, k)$ we are interested in for all triples $(q, p, k)$ of integers with $q\geq p\geq k>0$.

\begin{table}\begin{center}{\footnotesize
\setlength\arraycolsep{2pt}
$\begin{array}{cp{4pt}c}
	R^*/R =2/(T(6, 2) +1) =1/25		&&R^*/R =2/(T(5, 3) +1) =1/25\\[6pt]
	\begin{array}{cp{2pt}c}
		\begin{array}{cc|c|c|c|c|}
			\Psi^0	&\Psi^1	&\Psi^2	&\Psi^3	&\Psi^4&\Psi^5\\
	\hline
			0	&1						&0						&0						&0						&0\\
	\cline{1-2}
			0	&\multicolumn{1}{c}{2}	&2						&0						&0						&0\\
			2	&\multicolumn{1}{c}{1}	&2						&2						&2						&2\\
			2	&\multicolumn{1}{c}{2}	&0						&2						&2						&2\\
	\cline{1-3}
			0	&\multicolumn{1}{c}{3}	&\multicolumn{1}{c}{3}	&3						&0						&0\\
			3	&\multicolumn{1}{c}{1}	&\multicolumn{1}{c}{3}	&3						&3						&3\\
			3	&\multicolumn{1}{c}{3}	&\multicolumn{1}{c}{2}	&3						&3						&3\\
			3	&\multicolumn{1}{c}{3}	&\multicolumn{1}{c}{3}	&0						&3						&3\\
			3	&\multicolumn{1}{c}{3}	&\multicolumn{1}{c}{3}	&0						&3						&3\\
	\cline{1-4}
			0	&\multicolumn{1}{c}{4}	&\multicolumn{1}{c}{4}	&\multicolumn{1}{c}{4}	&4						&0\\
			4	&\multicolumn{1}{c}{1}	&\multicolumn{1}{c}{4}	&\multicolumn{1}{c}{4}	&4						&4\\
			4	&\multicolumn{1}{c}{4}	&\multicolumn{1}{c}{2}	&\multicolumn{1}{c}{4}	&4						&4\\
			4	&\multicolumn{1}{c}{4}	&\multicolumn{1}{c}{4}	&\multicolumn{1}{c}{3}	&4						&4\\
			4	&\multicolumn{1}{c}{4}	&\multicolumn{1}{c}{4}	&\multicolumn{1}{c}{4}	&0						&4\\
			4	&\multicolumn{1}{c}{4}	&\multicolumn{1}{c}{4}	&\multicolumn{1}{c}{4}	&0						&4\\
			4	&\multicolumn{1}{c}{4}	&\multicolumn{1}{c}{4}	&\multicolumn{1}{c}{4}	&0						&4\\
	\cline{1-5}
			0	&\multicolumn{1}{c}{5}	&\multicolumn{1}{c}{5}	&\multicolumn{1}{c}{5}	&\multicolumn{1}{c}{5}	&5\\
			5	&\multicolumn{1}{c}{1}	&\multicolumn{1}{c}{5}	&\multicolumn{1}{c}{5}	&\multicolumn{1}{c}{5}	&5\\
			5	&\multicolumn{1}{c}{5}	&\multicolumn{1}{c}{2}	&\multicolumn{1}{c}{5}	&\multicolumn{1}{c}{5}	&5\\
			5	&\multicolumn{1}{c}{5}	&\multicolumn{1}{c}{5}	&\multicolumn{1}{c}{3}	&\multicolumn{1}{c}{5}	&5\\
			5	&\multicolumn{1}{c}{5}	&\multicolumn{1}{c}{5}	&\multicolumn{1}{c}{5}	&\multicolumn{1}{c}{4}	&5\\
			5	&\multicolumn{1}{c}{5}	&\multicolumn{1}{c}{5}	&\multicolumn{1}{c}{5}	&\multicolumn{1}{c}{5}	&0\\
			5	&\multicolumn{1}{c}{5}	&\multicolumn{1}{c}{5}	&\multicolumn{1}{c}{5}	&\multicolumn{1}{c}{5}	&0\\
			5	&\multicolumn{1}{c}{5}	&\multicolumn{1}{c}{5}	&\multicolumn{1}{c}{5}	&\multicolumn{1}{c}{5}	&0\\
			5	&\multicolumn{1}{c}{5}	&\multicolumn{1}{c}{5}	&\multicolumn{1}{c}{5}	&\multicolumn{1}{c}{5}	&0\\
	\hline
		\end{array}&&
		\begin{array}{cc|c|c|c|c|}
			\Phi^0	&\Phi^1	&\Phi^2	&\Phi^3	&\Phi^4&\Phi^5\\
	\hline
			0	&1						&2						&3						&4						&5\\
	\cline{1-2}
			0	&\multicolumn{1}{c}{2}	&0						&0						&0						&0\\
			2	&\multicolumn{1}{c}{1}	&0						&2						&2						&2\\
			2	&\multicolumn{1}{c}{2}	&2						&2						&2						&2\\
	\cline{1-3}
			0	&\multicolumn{1}{c}{3}	&\multicolumn{1}{c}{3}	&0						&0						&0\\
			3	&\multicolumn{1}{c}{1}	&\multicolumn{1}{c}{3}	&0						&3						&3\\
			3	&\multicolumn{1}{c}{3}	&\multicolumn{1}{c}{2}	&0						&3						&3\\
			3	&\multicolumn{1}{c}{3}	&\multicolumn{1}{c}{3}	&3						&3						&3\\
			3	&\multicolumn{1}{c}{3}	&\multicolumn{1}{c}{3}	&3						&3						&3\\
	\cline{1-4}
			0	&\multicolumn{1}{c}{4}	&\multicolumn{1}{c}{4}	&\multicolumn{1}{c}{4}	&0						&0\\
			4	&\multicolumn{1}{c}{1}	&\multicolumn{1}{c}{4}	&\multicolumn{1}{c}{4}	&0						&4\\
			4	&\multicolumn{1}{c}{4}	&\multicolumn{1}{c}{2}	&\multicolumn{1}{c}{4}	&0						&4\\
			4	&\multicolumn{1}{c}{4}	&\multicolumn{1}{c}{4}	&\multicolumn{1}{c}{3}	&0						&4\\
			4	&\multicolumn{1}{c}{4}	&\multicolumn{1}{c}{4}	&\multicolumn{1}{c}{4}	&4						&4\\
			4	&\multicolumn{1}{c}{4}	&\multicolumn{1}{c}{4}	&\multicolumn{1}{c}{4}	&4						&4\\
			4	&\multicolumn{1}{c}{4}	&\multicolumn{1}{c}{4}	&\multicolumn{1}{c}{4}	&4						&4\\
	\cline{1-5}
			0	&\multicolumn{1}{c}{5}	&\multicolumn{1}{c}{5}	&\multicolumn{1}{c}{5}	&\multicolumn{1}{c}{5}	&0\\
			5	&\multicolumn{1}{c}{1}	&\multicolumn{1}{c}{5}	&\multicolumn{1}{c}{5}	&\multicolumn{1}{c}{5}	&0\\
			5	&\multicolumn{1}{c}{5}	&\multicolumn{1}{c}{2}	&\multicolumn{1}{c}{5}	&\multicolumn{1}{c}{5}	&0\\
			5	&\multicolumn{1}{c}{5}	&\multicolumn{1}{c}{5}	&\multicolumn{1}{c}{3}	&\multicolumn{1}{c}{5}	&0\\
			5	&\multicolumn{1}{c}{5}	&\multicolumn{1}{c}{5}	&\multicolumn{1}{c}{5}	&\multicolumn{1}{c}{4}	&0\\
			5	&\multicolumn{1}{c}{5}	&\multicolumn{1}{c}{5}	&\multicolumn{1}{c}{5}	&\multicolumn{1}{c}{5}	&5\\
			5	&\multicolumn{1}{c}{5}	&\multicolumn{1}{c}{5}	&\multicolumn{1}{c}{5}	&\multicolumn{1}{c}{5}	&5\\
			5	&\multicolumn{1}{c}{5}	&\multicolumn{1}{c}{5}	&\multicolumn{1}{c}{5}	&\multicolumn{1}{c}{5}	&5\\
			5	&\multicolumn{1}{c}{5}	&\multicolumn{1}{c}{5}	&\multicolumn{1}{c}{5}	&\multicolumn{1}{c}{5}	&5\\
	\hline
		\end{array}
	\end{array}&&
	\begin{array}{cp{2pt}c}
		\begin{array}{ccc|c|c|}
			\Psi^0	&\Psi^1	&\Psi^2						&\Psi^3						&\Psi^4\\
\hline
			0		&1		&2							&0							&0\\
\cline{1-3}
			0		&1		&\multicolumn{1}{c}{3}		&3							&0\\
			0		&3		&\multicolumn{1}{c}{2}		&3							&0\\
			3		&1		&\multicolumn{1}{c}{2}		&3							&3\\
			0		&3		&\multicolumn{1}{c}{3}		&0							&0\\
			3		&1		&\multicolumn{1}{c}{3}		&0							&3\\
			3		&3		&\multicolumn{1}{c}{2}		&0							&3\\
			3		&3		&\multicolumn{1}{c}{3}		&3							&3\\
\cline{1-4}
			0		&1		&\multicolumn{1}{c}{4}		&\multicolumn{1}{c}{4}		&4\\
			0		&4		&\multicolumn{1}{c}{2}		&\multicolumn{1}{c}{4}		&4\\
			0		&4		&\multicolumn{1}{c}{4}		&\multicolumn{1}{c}{3}		&4\\
			4		&1		&\multicolumn{1}{c}{2}		&\multicolumn{1}{c}{4}		&4\\
			4		&1		&\multicolumn{1}{c}{4}		&\multicolumn{1}{c}{3}		&4\\
			4		&4		&\multicolumn{1}{c}{2}		&\multicolumn{1}{c}{3}		&4\\

			0		&4		&\multicolumn{1}{c}{4}		&\multicolumn{1}{c}{4}		&0\\
			0		&4		&\multicolumn{1}{c}{4}		&\multicolumn{1}{c}{4}		&0\\
			4		&1		&\multicolumn{1}{c}{4}		&\multicolumn{1}{c}{4}		&0\\
			4		&1		&\multicolumn{1}{c}{4}		&\multicolumn{1}{c}{4}		&0\\
			4		&4		&\multicolumn{1}{c}{2}		&\multicolumn{1}{c}{4}		&0\\
			4		&4		&\multicolumn{1}{c}{2}		&\multicolumn{1}{c}{4}		&0\\
			4		&4		&\multicolumn{1}{c}{4}		&\multicolumn{1}{c}{3}		&0\\
			4		&4		&\multicolumn{1}{c}{4}		&\multicolumn{1}{c}{3}		&0\\

			4		&4		&\multicolumn{1}{c}{4}		&\multicolumn{1}{c}{4}		&4\\
			4		&4		&\multicolumn{1}{c}{4}		&\multicolumn{1}{c}{4}		&4\\
			4		&4		&\multicolumn{1}{c}{4}		&\multicolumn{1}{c}{4}		&4\\
\hline
		\end{array}&&
		\begin{array}{ccc|c|c|}
			\Phi^0	&\Phi^1	&\Phi^2						&\Phi^3						&\Phi^4\\
\hline
			0		&1		&2							&3							&4\\
\cline{1-3}
			0		&1		&\multicolumn{1}{c}{3}		&0							&0\\
			0		&3		&\multicolumn{1}{c}{2}		&0							&0\\
			3		&1		&\multicolumn{1}{c}{2}		&0							&3\\
			0		&3		&\multicolumn{1}{c}{3}		&3							&0\\
			3		&1		&\multicolumn{1}{c}{3}		&3							&3\\
			3		&3		&\multicolumn{1}{c}{2}		&3							&3\\
			3		&3		&\multicolumn{1}{c}{3}		&0							&3\\
\cline{1-4}
			0		&1		&\multicolumn{1}{c}{4}		&\multicolumn{1}{c}{4}		&0\\
			0		&4		&\multicolumn{1}{c}{2}		&\multicolumn{1}{c}{4}		&0\\
			0		&4		&\multicolumn{1}{c}{4}		&\multicolumn{1}{c}{3}		&0\\
			4		&1		&\multicolumn{1}{c}{2}		&\multicolumn{1}{c}{4}		&0\\
			4		&1		&\multicolumn{1}{c}{4}		&\multicolumn{1}{c}{3}		&0\\
			4		&4		&\multicolumn{1}{c}{2}		&\multicolumn{1}{c}{3}		&0\\
			0		&4		&\multicolumn{1}{c}{4}		&\multicolumn{1}{c}{4}		&4\\
			0		&4		&\multicolumn{1}{c}{4}		&\multicolumn{1}{c}{4}		&4\\
			4		&1		&\multicolumn{1}{c}{4}		&\multicolumn{1}{c}{4}		&4\\
			4		&1		&\multicolumn{1}{c}{4}		&\multicolumn{1}{c}{4}		&4\\
			4		&4		&\multicolumn{1}{c}{2}		&\multicolumn{1}{c}{4}		&4\\
			4		&4		&\multicolumn{1}{c}{2}		&\multicolumn{1}{c}{4}		&4\\
			4		&4		&\multicolumn{1}{c}{4}		&\multicolumn{1}{c}{3}		&4\\
			4		&4		&\multicolumn{1}{c}{4}		&\multicolumn{1}{c}{3}		&4\\
			4		&4		&\multicolumn{1}{c}{4}		&\multicolumn{1}{c}{4}		&0\\
			4		&4		&\multicolumn{1}{c}{4}		&\multicolumn{1}{c}{4}		&0\\
			4		&4		&\multicolumn{1}{c}{4}		&\multicolumn{1}{c}{4}		&0\\
\hline
		\end{array}
	\end{array}
\end{array}$
\caption{\label{tab-Gamma-rec}Inductive construction for $\Gamma\left((T(q, k) +1)/2, 1, q, k, k\right)$: illustration when $(k, q)\in \{(2, 6), (3, 5)\}$.}
}\end{center}\end{table}

\begin{algorithm}
\begin{algorithmic}[1]
	\State$\Psi,\Phi\gets\left\{(0, 1 ,\ldots, k -1)\right\}$
		\Comment{$(\Psi,\Phi)$ is $(k, k)$-ARPA of strength $k$}
 	\For{$i =k +1$ {\bf to} $q$}
		\State Extend $(\Psi, \Phi)$ to an $(i, k)$-ARPA of strength $k$ using \cref{algo-Gamma_rec}
	\EndFor
	\State\Return $(\Psi, \Phi)$ 
\end{algorithmic}
\caption{Given two positive integers $k$ and $q\geq k$, constructing an ARPA of $\Gamma(R, 1, q, k, k)$ for $R =(T(q, k) +1)/2$.
\label{algo-Gamma}}
\end{algorithm}

\begin{theorem}\label{thm-gamma_qpk}
Let $k >0$, $p \geq k$, and $q \geq p$ be three integers. 
If $p =q$, then $\gamma(q, q, k) =1$. 
If $p >k$, then $\gamma(q, p, k)\geq \gamma(q -p +k, k, k)$. 
Otherwise, we have:
\begin{align}\label{eq-gamma_qkk-LB}
	\gamma(q, k, k) &\textstyle\geq 2/\left(\sum_{r =0}^k \binom{q}{r}\binom{q -1 -r}{k -r} +1\right)
\end{align}%
\end{theorem}

\begin{proof}
The first two cases are trivial (see \cref{pty-gamma-base,pty-gamma-p>k}).
So let us assume that $k =p <q$. According to the previous discussion, it is sufficient to show the equality:
\begin{align}\label{eq-T-sum}
\textstyle 	1 +\sum_{i =k}^{q -1} T(i, k -1)	&=(T(q, k) +1)/2,
				&q\in\mathbb{N}\backslash\set{0 ,\ldots, k}
\end{align}	

This equality is satisfied at rank $q =k +1$ since by \cref{eq-T_base}, we have:
$$\begin{array}{rlll}
	1 +T(k, k -1)	&=2^k &=2^{k +1}/2 &=(T(k +1, k) +1)/2
\end{array}$$

Now suppose $q >k +1$, and the equality \cref{eq-T-sum} of rank $q -1$ holds. 
We observe successively:
$$\begin{array}{rll}
	1 +\sum_{i =k}^{q -1} T(i, k -1)
		&\!=(T(q -1, k) +1)/2 +T(q -1, k -1)	&\!\text{by induction hypothesis}\\
		&\!=(T(q, k) +1)/2						&\!\text{by \cref{eq-T_a+b+bis}}
\end{array}$$ 

So \cref{eq-T-sum} is verified, which completes the proof.
\end{proof}

By applying \cref{thm-gamma_qpk,thm-reduc-CD}, we obtain that for any three integers $k\geq 2$, $p\geq k$, and $q >p$, $\mathsf{k\,CSP\!-\!q}$ $D$-reduces to $\mathsf{k\,CSP\!-\!p}$ with an expansion factor on the approximation guarantee of:
$$\begin{array}{rll}	
	\gamma(q, p, k)	&\geq \gamma(q -p +k, k, k) &\geq 2/(T(q -p +k, k) +1)
\end{array}$$

It is not too hard to see that the following upper bound holds for the quantity $(T(a, b) +1)/2$ (see \cref{sec-T-UB} for a detailed proof):
\begin{align}\label{eq-T-UB}
	(T(a, b) +1)/2		&\leq (2a -b)^b/(2\times b!),	&a, b\in\mathbb{N},\ a >b
\end{align}
\Cref{thm-gamma_qpk,thm-reduc-CD} thus imply the following results:

\begin{corollary}\label{cor-reduc-dapx}
Let $k\geq 2$, $p\geq k$, and $q >p$ be three integers.
Then on an instance $I$ of $\mathsf{k\,CSP\!-\!q}$, the differential ratio of the best solutions among those whose coordinates take at most $p$ distinct values is least:
$$\begin{array}{rll}
	\gamma 	&:=2/\left(\sum_{r =0}^k \binom{q -p +k}{r}\binom{q -p +k -1 -r}{k -r} +1\right)	
			&\geq 2(k!)/(2q -2p +k)^k
\end{array}$$

Furthermore, if $\mathsf{k\,CSP\!-\!p}$ is approximable within some differential factor $\rho$, then $\mathsf{k\,CSP\!-\!q}$ is approximable within differential factor $\gamma\times\rho$.
\end{corollary}%

In particular, for all $q\geq 3$, $\mathsf{2\,CSP\!-\!q}$ $D$-reduces to $\mathsf{2\,CSP\!-\!2}$ with an expansion of $1/(q -1)^2$ on the approximation guarantee. 
Therefore, the result of \cite{N98} implies that for all $q\geq 2$, $\mathsf{2\,CSP\!-\!q}$ is differentially approximable within a constant factor using semidefinite programming along with derandomization techniques.

\begin{corollary}\label{cor-dapx-2CSP-q}
For all constant integers $q\geq 2$, $\mathsf{2\,CSP\!-\!q}$ is approximable within differential approximation ratio $(2 -\pi/2)/(q -1)^2\geq 0.429/(q -1)^2$. 
\end{corollary}%

\subsection{A refined analysis for the special case of $\mathsf{CSP(\mathcal{E}_q)}$}

As in the previous section, when restricting to $\mathsf{k\,CSP(\mathcal{E}_q)}$, the requirements on the manipulated arrays can be slightly relaxed. On the one hand, in an instance of $\mathsf{CSP(\mathcal{E}_q)}$, the constraints all evaluate the same on any two solutions $x$ and $x +\mathbf{a}$. 
On the other hand, for any $u\in\Sigma_q^q$ and any $a\in\Sigma_q$, the solution $y(x, u +\mathbf{a})$ consists of $y(x, u) +\mathbf{a}$. Thus, for $\mathsf{CSP(\mathcal{E}_q)}$, two solutions $y(x, u)$ and $y(x, u +\mathbf{a})$ realize the same objective value. 
This suggests that when reducing $\mathsf{k\,CSP(\mathcal{E}_q)}$ to $\mathsf{k\,CSP\!-\!p}$, the following slight relaxation of ARPAs should be considered:

\begin{definition}\label{def-reduc-CD_E}
Let $k >0$, $p\geq k$, and $q\geq p$ be three integers. 
Two arrays $\Psi$ and $\Phi$ with $q$ columns and entries in $\Sigma_q$ form a {\em relaxed $(q, p)$-alphabet reduction pair of arrays of strength $k$} if they satisfy the conditions below:
\begin{enumerate}
\item $\Phi$ contains at least 1 row of the form $(a, 1 +a ,\ldots, q -1 +a)$, for some $a\in\Sigma_q$;
\item \label{it-thm-q-p-P_E} the components of each row of $\Psi$ take at most $p$ different values;
\item \label{it-thm-q-p-P=Q_E} the function $\mu^\Psi_E -\mu^\Phi_E$ is balanced $k$-wise independent.
\end{enumerate} 

For two integers $R^* >0$ and $R\geq R^*$, we denote by $\Gamma_E(R, R^*, q, p, k)$ the (possibly empty) set of relaxed $(q, p)$-ARPAs of strength $k$ in which the total frequency of the words $(a, a +1, \ldots, a +q -1)$, $a\in\Sigma_q$ is $R^*$ and the arrays have $R$ rows each.

Furthermore, we define $\gamma_E(q, p, k)$ as the largest number $\gamma$ for which there exist two integers $R^* >0$ and $R\geq R^*$ such that $R^*/R =\gamma$ and $\Gamma_E(R, R^*, q, p, k)\neq\emptyset$.
\end{definition}

\Cref{tab-Gamma_E-ex} shows pairs of arrays that realize the value $\gamma_E(q, p, k)$. 
Similar to $\gamma(q, p, k)$, $\gamma_E(q, p, k)$ provides a lower bound on the expansion factor of the reduction from $\mathsf{k\,CSP\!-\!q}$ to $\mathsf{k\,CSP\!-\!p}$, when applied to instances of $\mathsf{k\,CSP(\mathcal{E}_q)}$. As in the general case, this quantity is also a lower bound on the differential ratio of a best solution among those whose coordinates take at most $p$ different values.

\begin{table}[t]\begin{center}{\footnotesize
\setlength{\arraycolsep}{3pt}
$\begin{array}{cp{8pt}c}
	\begin{array}{cp{4pt}c}
		\multicolumn{3}{c}{\gamma_E(4, 3, 2) =1/2}\\[6pt]	
		\begin{array}{cccc}
			\Psi^0	&\Psi^1	&\Psi^2	&\Psi^3\\\hline
			0		&0		&1		&2	\\
			0		&0		&1		&2	\\
			0		&1		&0		&3	\\
			0		&1		&2		&0	\\
			0		&1		&2		&0	\\
			0		&1		&2		&1	\\
			0		&1		&3		&3	\\
			0		&1		&3		&3	\\
			0		&2		&2		&3	\\
			0		&2		&2		&3	\\
			0		&3		&0		&1	\\
			0		&3		&2		&3	\\
		\end{array}&&\begin{array}{cccc}
			\Phi^0	&\Phi^1	&\Phi^2	&\Phi^3	\\\hline
			0		&0		&0		&0	\\
			0		&0		&3		&1	\\
			0		&1		&2		&3	\\
			0		&1		&2		&3	\\
			0		&1		&2		&3	\\
			0		&1		&2		&3	\\
			0		&1		&2		&3	\\
			0		&1		&2		&3	\\
			0		&2		&0		&2	\\
			0		&2		&1		&1	\\
			0		&3		&1		&0	\\
			0		&3		&3		&2	\\
		\end{array}
	\end{array}&&
	\begin{array}{cp{4pt}c}
	\multicolumn{3}{c}{\gamma_E(5, 4, 3) =1/3}\\[6pt]	
		\begin{array}{ccccc}
			\Psi^0	&\Psi^1	&\Psi^2	&\Psi^3	&\Psi^4\\\hline
			0	&0	&1	&2	&3	\\
			0	&0	&2	&3	&4	\\
			0	&0	&2	&3	&4	\\
			0	&1	&1	&3	&4	\\
			0	&1	&2	&2	&4	\\
			0	&1	&2	&3	&0	\\
			0	&1	&2	&3	&3	\\
			0	&1	&2	&4	&4	\\
			0	&1	&3	&3	&4	\\
			0	&2	&2	&3	&4	\\
			0	&2	&2	&3	&4	\\
			0	&2	&3	&4	&0	\\
		\end{array}&&\begin{array}{ccccc}
			\Phi^0	&\Phi^1	&\Phi^2	&\Phi^3	&\Phi^4	\\\hline
			0	&0	&1	&3	&4	\\
			0	&0	&2	&2	&4	\\
			0	&0	&2	&3	&3	\\
			0	&1	&1	&2	&3	\\
			0	&1	&2	&3	&4	\\
			0	&1	&2	&3	&4	\\
			0	&1	&2	&3	&4	\\
			0	&1	&2	&3	&4	\\
			0	&1	&3	&4	&0	\\
			0	&2	&2	&3	&0	\\
			0	&2	&2	&4	&4	\\
			0	&2	&3	&3	&4	\\
		\end{array}
	\end{array}
\end{array}$
\caption{\label{tab-Gamma_E-ex}
Pairs of arrays realizing $\gamma_E(4, 3, 2)$ and $\gamma_E(5, 4, 3)$.
These pairs have been computed by computer (see \cref{sec-cd}).
}}\end{center}\end{table}

\begin{theorem}\label{thm-reduc-CD_E}
For all constant integers $k\geq 2$, $p\geq k$, and $q\geq p$, $\mathsf{k\,CSP(\mathcal{E}_q)}$ $D$-reduces to $\mathsf{k\,CSP\!-\!p}$ with an expansion of $\gamma_E(q, p, k)$ on the approximation guarantee.

Furthermore, on an instance $I$ of $\mathsf{k\,CSP(\mathcal{E}_q)}$, the best solutions among those whose components take at most $p$ different values are $\gamma_E(q, p, k)$-approximate.
\end{theorem}
	
\begin{proof}
Suppose there exists $(\Psi,\Phi)\in\Gamma_E(R, R^*, q, p, k)$, where $R^* >0$ and $R\geq R^*$, and consider an instance $I$ of $\mathsf{CSP(\mathcal{E}_q)}$. We assume {\em w.l.o.g.} that the goal on $I$ is to maximize. We also may assume {\em w.l.o.g.} that the first column of $\Phi$ contains only zeros. If this is not the case, we can replace each row $\Phi_r$ of $\Phi$ with $\Phi_r^0\neq 0$ by $\Phi_r -\mathbf{\Phi_r^0}$, yielding a new element of $\Gamma_E(R, R^*, q, p, k)$ that satisfies this condition. 
Under this assumption, $R^*$ corresponds exactly to the number of rows of the form $(0, 1 ,\ldots, q -1)$ in $\Phi$. We can therefore use exactly the same argument as for the general case, which is based on the equality \cref{eq-vq-p--P=Q}. 

Let $i\in[m]$ and $x\in\Sigma_q^n$. Our goal is to show that $(\Psi,\Phi)$ satisifies the equality \cref{eq-vq-p--P=Q-i} of rank $(i, x)$. 
Since $I$ is an instance of $\mathsf{CSP(\mathcal{E}_q)}$, the function $P_i$ is invariant under a uniform shift of all its inputs. Therefore, a sufficient condition for $(\Psi, \Phi)$ to satisfy \cref{eq-vq-p--P=Q-i} at $(i, x)$ is that, for each $v\in\Sigma_q^{k_i}$, there are as many rows $\Psi_r$ of $\Psi$ satisfying $y(x, \Psi_r)_{J_i}\in\{v, v +\mathbf{1} ,\ldots, v +\mathbf{q -1}\}$ as there are in $\Phi$. 
Let $H =(c_1 ,\ldots, c_t)$ be the subsequence of $(0, 1 ,\ldots, q -1)$ induced by the coordinates $x_j$, $j\in J_i$ of $x$. 
Then, by definition of solutions $y(x, u)$, this occurs {\em if and only if} for each $v\in\Sigma_q^t$, there are as many rows $\Psi_r$ of $\Psi$ satisfying $\Psi_r^H\in\{v, v +\mathbf{1} ,\ldots, v +\mathbf{q -1}\}$ as there are in $\Phi$.
Equivalently, for all $v\in\Sigma_q^t$, the total frequency of words $u =(u_0, u_1 ,\ldots, u_{q -1})\in\Sigma_q^q$ such that $u_H\in\{v, v +\mathbf{1} ,\ldots, v +\mathbf{q -1}\}$ must be the same in $\Psi$ as in $\Phi$. 
Considering that $\card{H}\leq\card{J_i}\leq k$, we conclude that the arrays $\Psi$ and $\Phi$ do indeed satisfy \cref{eq-vq-p--P=Q-i} and hence \cref{eq-vq-p--P=Q}, provided that $\mu^\Psi_E -\mu^\Phi_E$ is balanced $k$-wise independent.
\end{proof}

In particular, for all $q\in\{3, 4, 5, 7, 8\}$, we have $\gamma_E(q, 2, 2) =1/q$ (see \cref{tab-gamma_E,sec-gEq22}). For these values of $q$, we deduce from \cref{thm-reduc-CD_E} and \cite{N98} that $\mathsf{2\,CSP(\mathcal{E}_q)}$ is approximable within a differential factor of $0.429/q$ (whereas for $\mathsf{2\,CSP\!-\!q}$ we only obtained a guarantee of $0.429/(q -1)^2$): 
\newcommand{\apx}[1]{}		
\newcommand{\pvo}[1]{#1}	
\begin{table}\footnotesize
$$\begin{array}{cc}
\begin{array}{|c|c|cccccc|}
\multicolumn{2}{c|}{} &\multicolumn{6}{c|}{q}\\\cline{3-8}
k&p	&3		&4		&5 		&6						&7						&8		\\
\hline\multirow{6}{*}{2}
&2	&1/3	&1/4	&1/5	&9/59					&1/7					&\pvo{1/8}	\\
&3	&-		&1/2	&2/5	&4/13					&2/7					&93/404	\\
&4 	&-		&-		&3/5	&7/15					&3/7					&3/8	\\
&5 	&-		&-		&-		&2/3					&11/21					&13/28	\\
&6 	&-		&-		&-		&-						&5/7					&4/7	\\
&7 	&-		&-		&-		&-						&-						&3/4	\\
\hline\multirow{5}{*}{3}
&3	&-		&1/4	&1/11	&38425/701342			&\apx{3676/107221}		&		\\
&4	&-		&-		&1/3	&1/6					&\pvo{5/52}				&		\\
&5	&-		&-		&-		&4/9					&\pvo{2/9}				&		\\
&6	&-		&-		&-		&-						&1/2					&		\\
&7	&-		&-		&-		&-						&-						&9/16	\\
\hline
\end{array}
&\begin{array}{|c|c|ccc|}	\multicolumn{5}{c}{}\\
\multicolumn{2}{c|}{} &\multicolumn{3}{c|}{q}\\\cline{3-5}
k&p	&5 		&6					&7						\\
\hline\multirow{2}{*}{4} 
&4	&1/11	&\apx{557/17632} 	&\apx{0.013964734}		\\
&5	&-		&1/6				&\apx{0.058898}			\\
\hline\multirow{2}{*}{5} 
&5	&-		&1/16				&\apx{0.01281777623}	\\
&6	&-		&-					&1/10					\\
\comment{\hline\multirow{1}{*}{6}
&6	&-		&-					&\apx{1/42}				&		\\}
\hline
\end{array}
\end{array}$$
\caption{Numbers $\gamma_E(q, p, k)$ for some triples $(q, p, k)$. 
These values can be calculated by computer by solving linear programs (see \cref{sec-cd} for more details).
}\label{tab-gamma_E}
\end{table}

\begin{corollary}\label{cor-reduc-2E_q}
For $q\in\{3, 4, 5, 7, 8\}$, $\mathsf{2\,CSP(\mathcal{E}_q)}$ is approximable within differential approximation ratio $0.429/q$.
\end{corollary}%

\subsection{Concluding remarks}

From \cref{cor-reduc-dapx}, we know that $\mathsf{k\,CSP\!-\!q}$ reduces to $\mathsf{k\,CSP\!-\!k}$ with a constant-order expansion on the differential approximation guarantee. Together with the result of \cite{N98}, this allows us to obtain new constant approximation results for 2-CSPs. 
The question of whether $\mathsf{k\,CSP\!-\!q}$ is approximable to within some constant differential factor remains open for pairs of integers $k$ and $q$ such that either $\min\set{k, q}\geq 3$, or $k\geq 4$ and $q =2$. Nevertheless, \cref{cor-reduc-dapx} allows us to restrict this question to the case where $k\geq q$.

Recently, we showed that the ARPAs constructed by \cref{algo-Gamma} realize the numbers $\gamma(q, k, k)$\footnote{This result is available in the arXiv preprint {\em J.-F. Culus, S. Toulouse, Optimizing alphabet reduction pairs of arrays, arXiv:2406.10930 (2024)}.\label{CT24-reg0}}. 
Therefore, the obtained estimate of the expansion factor of the reduction from $\mathsf{k\,CSP\!-\!q}$ to $\mathsf{k\,CSP\!-\!k}$ --- and of the differential ratio of optimal solutions over $\{T^n\,|\,T\subseteq\Sigma_q: |T| =k\}$ --- is the best possible within our framework. 
\begin{table}\footnotesize
$$\begin{array}{cc}
\begin{array}{|c|c|ccccc|}
\multicolumn{2}{c|}{} 			&\multicolumn{5}{c|}{q}
\\\cline{3-7}		k		&p	&3		&4		&5 		&6		&7		
\\\hline\multirow{5}{*}{2}	&2	&1/4	&1/9	&1/16	&1/25	&1/36	
\\							&3	&-		&1/3	&1/6	&1/10	&1/15	
\\							&4	&-		&-		&4/9	&1/4	&4/25	
\\							&5	&-		&-		&-		&1/2	&3/10	
\\							&6	&-		&-		&-		&-		&9/16	
\\\hline\multirow{4}{*}{3}	&3	&-		&1/8	&1/25	&1/56	&1/105	
\\							&4	&-		&-		&1/5	&2/27	&1/28	
\\							&5	&-		&-		&-		&1/4	&5/49	
\\							&6	&-		&-		&-		&-		&2/7	
\\\hline
\end{array}
&\begin{array}{|c|c|ccc|}	\multicolumn{5}{c}{}\\
\multicolumn{2}{c|}{} 			&\multicolumn{3}{c|}{q}
\\\cline{3-5}		k		&p	&5 		&6		&7		
\\\hline\multirow{3}{*}{4}	&4	&1/16	&1/65	&1/176	
\\							&5	&-		&1/10	&1/36	
\\							&6	&-		&-		&5/33	
\\\hline\multirow{2}{*}{5}	&5	&-		&1/32	&1/161	
\\							&6	&-		&-		&2/35	
\\\hline\multirow{1}{*}{6}	&6	&-		&-		&1/64	
\comment{
\\							&7	&-		&-		&-		&1/35	&1/190
\\							&8	&-		&-		&-		&-		&4/85
\\\hline\multirow{2}{*}{7}	&7	&-		&-		&-		&1/128	&1/897
\\							&8	&-		&-		&-		&-		&5/315
\\\hline\multirow{1}{*}{8}	&8	&-		&-		&-		&-		&1/256
}
\\\hline
\end{array}
\end{array}$$
\caption{Numbers $\gamma(q, p, k)$ for some triples $(q, p, k)$. 
These values can be calculated by computer by solving linear programs (see \cref{sec-cd} for more details).
}\label{tab-gamma}
\end{table}

 %
For the case where $p >k$, our estimate relies on the inequality $\gamma(q, p, k)\geq\gamma(q -p +k, k, k)$. However, given three integers $k >0$, $p >k$, and $q >p$, $\gamma(q, p, k)$ is greater than $\gamma(q -1, p -1, k)$\textsuperscript{\ref{CT24-reg0}}. 
For instance, we have (see \cref{tab-gamma}):
$$\begin{array}{rlll}
	\gamma(6, 4, 2)  =1/4 	&>\gamma(5, 3, 2) =1/6	&>\gamma(4, 2, 2) =1/9	&\\
	\gamma(6, 5, 3)  =1/4 	&>\gamma(5, 4, 3) =1/5	&>\gamma(4, 3, 3) =1/8	&\\
\end{array}$$

Thus, a closer study of $(q, p)$-ARPAs of strength $k <p$ could provide a finer estimate of the expansion factor of the reduction when reducing to $\mathsf{k\,CSP\!-\!p}$ for some $p\in\set{k +1 ,\ldots, q -1}$, as well as a tighter bound on the highest differential ratio achieved by a solution whose coordinates take at most $p$ different values.
 %
Similarly, investigating relaxed ARPAs has the potential to improve the results obtained when reducing from $\mathsf{k\,CSP(\mathcal{E}_q)}$, especially for the case where $p =k =2$, since an approximation algorithm is known for $\mathsf{2\,CSP(\mathcal{E}_2)}$. 
Although for all integers $k\geq 2$, $p \geq k$ and $q >p$, $\mathsf{k\,CSP(\mathcal{E}_q)}$ $D$-reduces to $\mathsf{k\,CSP\!-\!p}$ with a multiplicative factor of $\gamma(q, p, k)$ --- for which we know a lower bound --- on the approximation guarantee, it is very likely that for three such integers,  $\gamma_E(q, p, k) >\gamma(q, p, k)$. This is notably true in all the few cases we have computed. For instance (see \cref{tab-gamma,tab-gamma_E}), we have:
$$\begin{array}{rlcrl}
	\gamma_E(5, 2, 2)/\gamma(5, 2, 2) &=16/5	&&\gamma_E(5, 3, 2)/\gamma(5, 3, 2) &=12/5\\
	\gamma_E(5, 3, 3)/\gamma(5, 3, 3) &=25/11	&&\gamma_E(5, 4, 3)/\gamma(5, 4, 3) &=5/3
\end{array}$$

Designing such pairs of arrays appears to be more challenging compared to ARPAs. 
For these latter, we observe that two vectors $u, v\in\mathbb{Z}_q^k$ satisfying $u =v$ can always be viewed as two vectors of $\mathbb{Z}_{q +1}^k$ that satisfy $u =v$. 
Thus, $(q, p)$-ARPAs can be interpreted as partial $(q +1, p)$-ARPAs. 
In contrast, two vectors $u, v\in\mathbb{Z}_q^k$ may satisfy the condition $u_j -u_1\equiv v_j -v_1\bmod{q}$, $j\in\{2 ,\ldots, k\}$, but $u_j -u_1\not\equiv v_j -v_1\bmod{(q +1)}$ for some $j\in\set{2 ,\ldots, k}$. Consequently, in the most general case, a relaxed $(q, p)$-ARPAs cannot be interpreted as a pair of subarrays of a relaxed $(q +1, p)$-ARPAs, and vice versa.

\section{At the neighborhood of any solution}\label{sec-vois}
\label{sec-vois-prev} 
In \cref{sec-average}, we investigate whether the average solution value provides some differential approximation guarantee. In this section, we address a similar question: we analyze the differential ratio of optimal solutions on Hamming balls of a given radius $d\geq k$, or on the union of the shifts of such a ball by vectors of the form $\mathbf{a}$. 
For example, we observed in \cref{sec-average} that for $\mathsf{CSP(\mathcal{O}_q)}$, a best solution over the neighborhood $\tilde{B}^0(x)$ of an arbitrary solution $x$ is $1/q$-differential approximate.

\subsection{From the average solution value to solutions with optimal value on Hamming balls of radius 1}
\label{sec-B1} 

A common heuristic approach for CSPs is to fix some radius $d$ and then compute a {\em local optimum with respect to} (w.r.t.) $B^d$, i.e., a solution $x$ such that $v(I, x)\geq v(I, y)$ for all $y\in B^d(x)$ when maximizing, with value $v(I, x)\leq v(I, y)$ for all $y\in B^d(x)$ when minimizing. 
Finding these solutions takes only polynomial time for instances where the objective function takes integer values and the diameter is polynomially bounded. Otherwise, finding local optima w.r.t. $B^d$ is $\mathbf{PLS}$-complete, even for $\mathsf{Max\,Cut}$ and the $B^1$ neighborhood function \cite{SY91}. 

Many articles address the question of whether such solutions provide an approximation guarantee \cite{KMSV99,MPT03A,AP95A,A97}. 
In particular, Khanna et al. showed that for $\mathsf{Max\,2\,CCSP}$, there is no constant integer $d$ for which local search with respect to $B^d$ guarantees any constant-factor standard approximation \cite{KMSV99}. In the differential approximation paradigm, this result extends by reduction to $\mathsf{Max\,Ek\,Sat}$ for all $k\geq 2$, \cite{MPT03A}. 

We identify a few cases for which the results of \cref{sec-average} imply a differential approximation guarantee at local optima with respect to the neighborhood function either $B^1$ or $\tilde{B}^1$, as well as at solutions with optimal value in the neighborhood $B^1(x)$ or $\tilde{B}^1(x)$ of any solution $x$.

\begin{property}\label{pty-I_q^(k-1)-B^1}
On an instance $I$ of $\mathsf{Ek\,CSP(\mathcal{I}_q^{k -1})}$, where $q\geq 2$ and $k\geq 2$, the differential ratio of local optima with respect to $B^1$ is at least the average differential ratio. Furthermore, the differential ratio of solutions with optimal value over Hamming balls of radius 1 is at least $1/n \times kq/(q -1)$ times this ratio.
\end{property}

\begin{proof}
Let $x\in\Sigma_q^n$. 
We want to evaluate the sum of the solution values over $B^1(x)\backslash\set{x}$. 
Consider a constraint $C_i =P_i(x_{i_1} ,\ldots, x_{i_k})$ of $I$.
Over $B^1(x)\backslash\set{x}$, $P_i$ is evaluated as follows:
\begin{itemize}
	\item for all $s\in[k]$, once at each 
			$(x_{i_1} ,\ldots, x_{i_{s -1}}, x_{i_s} +a, x_{i_{s +1}} ,\ldots, x_{i_k})$, $a\in[q -1]$; 
	\item $(q -1) \times (n -k)$ times at $(x_{i_1} ,\ldots, x_{i_k})$.
\end{itemize}
Since $P_i\in\mathcal{I}_q^{k -1}$, by \cref{eq-I_q^t}, for all $s\in[k]$ we have:
$$\begin{array}{rl}
\sum_{a =1}^{q -1}P_i(x_{i_1} ,\ldots, x_{i_{s -1}}, x_{i_s} +a, x_{i_{s +1}} ,\ldots, x_{i_k})
	&=q\times r_{P_i} -P_i(x_{i_1} ,\ldots, x_{i_k}) 
\end{array}$$
Thus, the sum of the evaluations of $C_i$ over $B^1(x)\setminus\set{x}$ gives the following expression:
\begin{align}\label{eq-Iqk-1}
\textstyle\sum_{y\in B^1(x): y\neq x} P_i(y_{i_1} ,\ldots, y_{i_k})
	&=k\times q\,r_{P_i} +\left((q -1) n -qk\right)\times P_i(x_{i_1} ,\ldots, x_{i_k})
\end{align}

Noting that $(q -1)n$ is the cardinality of $\card{B^1(x)\backslash\set{x}}$, the sum of the solution values over $B^1(x)\backslash\set{x}$ is therefore:
\begin{align}\label{eq-sum_B^1}
\textstyle\sum_{y\in B^1(x): y\neq x} v(I, y)	&\textstyle 
	=qk\times\mathbb{E}_X[v(I, X)] +\left(\card{B^1(x)\backslash\set{x}} -qk\right)\times v(I, x)
\end{align}

Since $v(I,x)$ cannot be worse than $\mathrm{wor}(I)$, it follows from \cref{eq-sum_B^1} that the average differential ratio over $B^1(x)\backslash\set{x}$ is at least a fraction $qk/\card{B^1(x)\backslash\set{x}}$ of the average differential ratio reached at $\mathbb{E}_X[v(I, X)]$.  
Furthermore, if $x$ is a local optimum with respect to $B^1$, then assuming {\em w.l.o.g.} that the goal on $I$ is to maximize, we have:
\begin{align}\label{eq-avB1}
	\card{B^1(x)\backslash\set{x}} \times v(I, x)	&\textstyle\geq\sum_{y\in B^1(x): y\neq x} v(I, y)
\end{align}

Combining \cref{eq-avB1} with \cref{eq-sum_B^1}, we get the inequality $v(I, x)\geq\mathbb{E}_X[v(I, X)]$, meaning that the differential ratio achieved by $x$ on $I$ is at least the average differential ratio.
\end{proof}

It follows from \cref{pty-I_q^(k-1)-B^1} that the differential guarantees established in \cref{sec-average} for $\mathsf{Ek\,CSP(\mathcal{I}_q^{k -1})}$, which hold at the average solution value, also hold for local optima with respect to $B^1$.
Furthermore, these guarantees extend to solutions with optimal value over Hamming balls of radius 1, although up to a multiplicative factor of $O(1/n)$ on the approximation guarantee. 
Note that $\mathsf{Ek\,CSP(\mathcal{I}_q^{k -1})}$ covers the restriction of $\mathsf{Lin\!-\!q}$ to equations of the form $x_{i_1} +\ldots+ x_{i_k}\equiv\alpha_i\bmod{q}$. In particular, for $q =2$, according to \cref{cor-E-2CSPs,cor-E-kCSP-2}, local optima with respect to $B^1$ and solutions with optimal value over Hamming balls of radius 1 provide a differential approximation guarantee of $\Omega(1/n^k)$ and $\Omega(1/n^{k +1})$, respectively, for $\mathsf{E(2k)Lin\!-\!2}$. 
For $\mathsf{E2\,Lin\!-\!2}$, we obtain precisely the ratios $1/(2\lceil\nu/2\rceil)$ and $2/(\lceil \nu/2\rceil\times n)$.  

Now we see that for $\mathsf{3\,CSP(\mathcal{E}_2)}$ similar conclusions can be drawn as for $\mathsf{Ek\,CSP(\mathcal{I}_q^{k -1})}$:

\begin{property}\label{pty-3CSP(E_2)-B^1}
For $\mathsf{3\,CSP(\mathcal{E}_2)}$, the differential ratio of local optima with respect to $B^1$ is at least $1/(2\lceil\nu/2\rceil)$, and the differential ratio of solutions with optimal value over Hamming balls of radius 1 is at least $4/n$ times this bound.
\end{property}

\begin{proof}
Consider a constraint $C_i =P_i(x_{J_i})$ of $I$. 
If $C_i$ is of the form $P_i(x_{i_1}, x_{i_2})$, then since $\mathcal{E}_2\subseteq\mathcal{I}_2^1$, $C_i$ satisfies \cref{eq-Iqk-1} with $q =k =2$. 
If it is of the form $P_i(x_{i_1}, x_{i_2}, x_{i_3})$, then over $B^1(x)\backslash\set{x}$, $C_i$ is evaluated once at $(\bar x_{i_1}, x_{i_2}, x_{i_3})$, $(x_{i_1}, \bar x_{i_2}, x_{i_3})$, and $(x_{i_1}, x_{i_2}, \bar x_{i_3})$, and $n -3$ times at $(x_{i_1}, x_{i_2}, x_{i_3})$. 
Since $P_i\in\mathcal{E}_2$, $P_i(\bar x_{i_1}, x_{i_2}, x_{i_3}) +P_i(x_{i_1}, \bar x_{i_2}, x_{i_3}) +P_i(x_{i_1}, x_{i_2}, \bar x_{i_3}) +P_i(x_{i_1}, x_{i_2}, x_{i_3})$ is one half of the sum of the values taken by $P_i$ over $\set{0, 1}^3$, and hence is equal to $4 \times r_{P_i}$.
We deduce that, also in this case, the sum of the evaluations of $C_i$ over $B^1(x)\backslash\set{x}$ is equal to $4 \times r_{P_i} +(n -4) P_i(x_{J_i})$.  
So the solution values over $B^1(x)\backslash\set{x}$ verify \cref{eq-sum_B^1} with $qk/\card{B^1(x)\backslash\set{x}} =4/n$.
The rest of the proof is the same as for \cref{pty-I_q^(k-1)-B^1}, with the addition of the lower bound of $1/(2\lceil\nu/2\rceil)$ given in \cref{cor-E-kCSP(E_q)} for the differential ratio reached at $\mathbb{E}_X[v(I, X)]$.
\end{proof}

Finally, the approximation guarantees of \cref{pty-3CSP(E_2)-B^1} somehow extend, by reduction, to $\mathsf{2\,CSP\!-\!2}$:

\begin{property}\label{pty-2CSP-B^1}
For $\mathsf{2\,CSP\!-\!2}$, the differential ratio of local optima with respect to $\tilde{B}^1$ is  at least $1/(2\lceil(\nu +1)/2\rceil)$, and the differential ratio of a best solution in the neighborhood $\tilde{B}^1(x)$ of any solution $x$ is at least $4/(n +1)$ times this bound.
\end{property}

\begin{proof}
Consider an instance $I$ of $\mathsf{2\,CSP\!-\!2}$. 
From $I$, we construct an instance $J$ of $\mathsf{3\,CSP(\mathcal{E}_2)}$ as follows:
	first, we introduce a new binary variable $z_0$;
	then, we replace each constraint $P_i(x_{i_1})$ or $P_i(x_{i_1}, x_{i_2})$ of $I$ with the new constraint $P_i(x_{i_1} +z_0)$ or $P_i(x_{i_1} +z_0, x_{i_2} +z_0)$. 
(Note that such a transformation is quite common, see {\em e.g.} \cite{GW95}.)
Since the construction adds the same new variable to the support of each constraint, the strong chromatic number of $J$ is $\nu +1$. Furthermore, the objective functions of $I$ and $J$ satisfy:
$$\begin{array}{rlrll}
	v(J, (x, z_0))	&=v(I, (x_1 +z_0 ,\ldots, x_n +z_0)),	&x\in\set{0, 1}^n,\ z_0\in\set{0, 1}
\end{array}$$

Consequently, for any solution $(x, z_0)$ of $J$, this solution and its corresponding solution $x +\mathbf{z_0}$ of $I$ achieve the same differential ratio on their respective instances. 
Consider then a solution $\tilde{x}$ of $I$ that is optimal over the neighborhood $\tilde{B}^1(x)$ of some other solution $x$ of $I$. We can assume {\em w.l.o.g.} that $\tilde{x}\in B^1(x)$, as otherwise we replace $x$ by $\bar x$. Suppose {\em w.l.o.g.} that the goal on $I$ (and hence on $J$) is to maximize. We observe:
\begin{itemize}
	\item $v(I, \tilde{x})\geq v(I, y)$, $y\in B^1(x)$ 
			$\Leftrightarrow$ $v(J, (\tilde{x}, 0))\geq v(J, (y, 0))$, $y\in B^1(x)$;
	\item $v(I, \tilde{x})\geq v(I, \bar x)$ 
			$\Leftrightarrow$ $v\left(J, (\tilde{x}, 0)\right)\geq v\left(J, (x, 1)\right)$;
\end{itemize}

Thus, the assumption that on $I$, $\tilde{x}$ is optimal over $\tilde{B}^1(x)$, and hence over $B^1(x)\cup\set{\bar x}$, implies that on $J$, $(\tilde{x}, 0)$ is optimal over $B^1(x)\times\set{0} \cup\set{(x, 1)}$, which is exactly $B^1((x, 0))$.  
In particular, if $x$ on $I$ is a local optimum with respect to $\tilde{B}^1$, then $(x, 0)$ on $J$ is a local optimum with respect to $B^1$. 
We conclude by observing that \cref{pty-3CSP(E_2)-B^1} applies to instance $J$, which is $(\nu +1)$-partite and manipulates $n +1$ variables. 
\end{proof}

We summarize the approximation guarantees induced by 	\cref{pty-I_q^(k-1)-B^1,pty-3CSP(E_2)-B^1,pty-2CSP-B^1} in \cref{tab-B1}. 
\begin{table}\footnotesize{\begin{center}
\begin{tabular}{l|l|l|l}
Restriction	&Conditions on $q$, $k$ and $\nu$	&$\!$Loc. opt.	&$\!$Opt. sol.\\
\hline&&&\\[-8pt]
\multirow{6}{*}{$\mathsf{Ek\,CSP(\mathcal{I}_q^{k -1})}$}
	&$\nu\leq 2k -1$	&$\!1/q^{\nu -k +1}$			&\multirow{4}{*}{$\!\Omega(1/n)$}
\\[1pt]\cline{2-3}&&&\\[-8pt]
 	&$\nu\leq 2k$		&\multirow{2}{*}{$\!1/q^k$}		&
\\[1pt]\cline{2-2}&&&\\[-8pt]
	&$q$ prime power $>k$ and $\nu\leq q +k$	&		&
\\[1pt]\cline{2-3}&&&\\[-8pt]
	&$2^{\lceil\log_2q\rceil} >k$ and $\nu\leq 2^{\lceil\log_2q\rceil} +1$	&$\!1/(2(q -1))^k$	&
\\[1pt]\cline{2-4}&&&\\[-8pt]
	&$q\geq 3$ and $\nu\geq k$	&$\!1/O(\nu^{k -\lceil\log_{\Theta(q)} k\rceil})$	
									&$\!1/O(\nu^{k -\lceil\log_{\Theta(q)} k\rceil}\times n)$
\\[1pt]\hline&&&\\[-8pt]
$\mathsf{Ek\,CSP(\mathcal{I}_2^{k -1})}$
	&$k\geq 3$ and $\nu\geq k$		&$\!1/O(\nu^{\lfloor k/2\rfloor})$
												&$\!1/O(\nu^{\lfloor k/2\rfloor}\times n)$
\\[1pt]\hline&&&\\[-8pt]
\multirow{2}{*}{$\mathsf{E3\,CSP(\mathcal{I}_q^2)}$}
	&$q$ power of 2 $>3$ and $\nu\leq q +4$			&$\!1/q^3$
													&\multirow{2}{*}{$\!\Omega(1/n)$}
\\[1pt]\cline{2-3}&&&\\[-8pt]
	&$2^{\lceil\log_2q\rceil} >3$ and $\nu\leq 2^{\lceil\log_2q\rceil} +2$
													&$\!1/(2(q -1))^3$		&
\\[1pt]\hline&&&\\[-8pt]
$\mathsf{E2\,CSP(\mathcal{I}_q^1)}$ &$\nu\geq 2$	&\multirow{3}{*}{$\!1/O(\nu)$}
													&\multirow{3}{*}{$\!1/O(\nu\times n)$}
\\[1pt]\cline{1-2}&&&\\[-8pt]
$\mathsf{3\,CSP(\mathcal{E}_2)}$ 	&$\nu\geq 3$	&			&
\\[1pt]\cline{1-2}&&&\\[-8pt]
$\mathsf{2\,CSP\!-\!2}$ 			&$\nu\geq 2$	&			&
\end{tabular}

\caption{\label{tab-B1}
Lower bounds on the differential ratio of local optima (the ``Loc. opt.'' column), or of solutions with optimal value in the neighborhood of any solution (the ``Opt. sol.'' column), with respect to the neighborhood functions $B^1$ for $\mathsf{Ek\,CSP(\mathcal{I}_q^{k -1})}$ (where $q, k\geq 2$) and 	$\mathsf{3\,CSP(\mathcal{E}_2)}$, and $\tilde{B}^1$ for $\mathsf{2\,CSP\!-\!2}$.} 
\end{center}}\end{table}

\subsection{Hamming balls with radius at least $k$}

We investigate whether, for $\mathsf{k\,CSP\!-\!q}$, extremal solutions over Hamming balls of fixed radius $d\geq k$ provide any differential approximation guarantees.
Consider an instance $I$ of $\mathsf{k\,CSP\!-\!q}$, and let $x^*$ and $x$ be two solutions of $I$ where $x^*$ is optimal. We denote by $\kappa$ the Hamming distance between $x^*$ and $x$. 
Assuming $\kappa\geq k$, consider an integer $d\in\set{k ,\dots, \kappa}$. 
We are interested in the vectors from $\set{x^*_1, x_1}\times\ldots\times\set{x^*_n, x_n}$ that are at Hamming distance $d$ from $x$, and denote by $N^d(x^*, x)$ the set of such vectors. 
The average solution value over $N^d(x^*, x)$ can be expressed as:
$$\begin{array}{l}
	\sum_{y\in N^d(x^*, x)} v(I, y)/\card{N^d(x^*, x)}\\
	\ \ \ \ =\sum_{i =1}^m w_i\times\sum_{y\in N^d(x^*, x)} P_i(y_{J_i})/\binom{\kappa}{d}\\
	\ \ \ \ =\sum_{i =1}^m w_i 	\left(
									\sum_{y\in N^d(x^*, x): y_{J_i} =x^*_{J_i}} P_i(x^*_{J_i}) 
									+\sum_{y\in N^d(x^*, x): y_{J_i} \neq x^*_{J_i}} P_i(y_{J_i})
								\right)/\binom{\kappa}{d}
\end{array}$$ 

For each constraint $C_i$, let $\kappa_i$ be the number of indices $j\in J_i$ for which $x_j\neq x^*_j$. Then the number of vectors $y\in N^d(x^*, x)$ satisfying $y_{J_i} =x^*_{J_i}$ is $\binom{\kappa -\kappa_i}{d -\kappa_i}$. If for all $i\in[m]$, $w_i P_i$ is non-negative, then $w_i P_i(y_{J_i})\geq 0$ for all $y\in N^d(x^*, x)$. If this is the case, and the goal is to maximize on $I$, then we deduce that the average solution value over $N^d(x^*, x)$ is bounded below by:
$$\begin{array}{rll}
	\displaystyle	 		\frac{\min_{i =1}^m \binom{\kappa -\kappa_i}{d -\kappa_i}}{\binom{\kappa}{d}} 	\times	v(I, x^*)
	&\geq	\displaystyle	\frac{\binom{\kappa -k}{d -k}}{\binom{\kappa}{d}}							\times	\mathrm{opt}(I)
	&=		\displaystyle	\frac{d(d -1) \ldots (d -k +1)}{\kappa(\kappa -1) \ldots (\kappa -k +1)}		\times	\mathrm{opt}(I)
\end{array}$$

Considering that $N^d(x^*, x)\subseteq B^d(x)$ while $\kappa\leq n$, we conclude that the highest standard ratio achieved over $B^d(x)$ is at least $k!\binom{d}{k}/n^k$.
In contrast, deriving a similar conclusion for the differential ratio appears to be more challenging. In particular, in the general case, there is no straightforward way to compare the two quantities:
$$\textstyle\sum_{y\in N^d(x^*, x)} v(I, y) -\binom{\kappa -k}{d -k}\times\mathrm{opt}(I)
\text{ and }\left(\binom{\kappa}{d} -\binom{\kappa -k}{d -k}\right)\times\mathrm{wor}(I)$$

To evaluate the highest differential ratio of a solution in the neighborhood $B^d(x)$ of an arbitrary solution $x$, where $d\geq k$, we introduce an approach based on solution multisets. We associate with each pair $(x^*, x)\in\Sigma_q^n\times\Sigma_q^n$ a pair $\left(\mathcal{X}(I, x^*, x), \mathcal{Y}(I, x^*, x)\right)$ of solution multisets. These multisets must be of the same size $R$, and satisfy that
	$\mathcal{X}(I, x^*, x)$ is a subset of $B^d(x)$, 
	$\mathcal{Y}(I, x^*, x)$ contains $x^*$ a positive number (denoted $R^*$) of times, 
	and the sum of solution values is the same over both solution multisets.
Provided that $x^*$ is optimal, the best solutions over $\mathcal{X}(I, x^*, x)$, and thus the best solutions over $B^d(x)$, realize a differential ratio of at least $R^*/R$.

\subsection{Partition-based solution multisets}

Since the case where $B^d(x)$ contains an optimal solution is trivial, we consider pairs $(x^*, x)$ of solutions that are within a Hamming distance of at least $d +1$ from each other. Furthermore, we restrict the solution multisets $\mathcal{X}(I, x^*, x)$ and $\mathcal{Y}(I, x^*, x)$ to solutions belonging to the set $\set{x^*_1, x_1}\times,\ldots,\times\set{x^*_n, x_n}$ of vectors. 
Using an approach similar to that in \cref{sec-reduc}, we precisely define our solution multisets $\mathcal{X}$ and $\mathcal{Y}$ by considering the following framework.

\smallskip
$\bullet$ {\bf Solution multisets association.} Let $\kappa\in\set{d +1 ,\ldots, n}$, and let $x^*$, $x$ be two vectors of $\Sigma_q^n$ that differ at exactly $\kappa$ coordinates. We denote the set of indices of these coordinates by $\mathcal{J}(x^*, x) =\set{j_1 ,\ldots, j_\kappa}$. We associate with $(x^*, x)$ and each $u\in\set{0, 1}^\kappa$ the solution $y(x^*, x, u)$ defined by: 
\begin{align}\label{eq-vois--y^u}
	y(x^*, x, u)_j	&=\left\{\begin{array}{rl}
						x^*_j	&\text{if $j =j_c\in\mathcal{J}(x^*, x)$ and $u_c =1$}\\
						x_j		&\text{otherwise}
					\end{array}\right.
\end{align}

Thus, $y(x^*, x, u)$ is obtained from $x$ by changing the coordinate $x_{j_c}$ to $x^*_{j_c}$ for each $c\in[\kappa]$ such that $u_c =1$. Therefore, $y(x^*, x, \mathbf{1}) =x^*$, while the number of non-zero coordinates of $u$ determines the Hamming distance of $y(x^*, x, u)$ from $x$. In particular, $y(x^*, x, \mathbf{0}) =x$.

We then consider two binary arrays $\Psi$ and $\Phi$ with $\kappa$ columns and the same number of rows, denoted $R$. With such a pair $(\Psi, \Phi)$, we associate the solution multisets:
\begin{align}\nonumber 
\mathcal{X}(I, x^*, x)	=\left(y(x^*, x, \Psi_r)\,|\,r\in[R]\right),
	\ \mathcal{Y}(I, x^*, x) =\left(y(x^*, x, \Phi_r)\,|\,r\in[R]\right)
\end{align}

$\bullet$ {\bf Conditions.} To model solutions of $B^d(x)$, the number of 1's in each row of $\Psi$ must be at most $d$. 
To model the solution $x^*$, $\Phi$ must contain at least one row of 1's.
Moreover, since our goal is to relate solution values over $B^d(x)$ to $\mathrm{opt}(I)$ under the assumption that $x^*$ is optimal, we require that $\Psi$ and $\Phi$ satisfy:
\begin{align}\label{eq-vois--P=Q}
	\textstyle\sum_{r =1}^R v(I, y(x^*, x, \Psi_r))	&\textstyle=\sum_{r =1}^R v(I, y(x^*, x, \Phi_r))
\end{align}

Without loss of generality, assume that the goal on $I$ is to maximize. 
We denote by $R^*$ the number of rows of 1's in $\Phi$. Then, under these conditions, solutions that perform the best objective value over $B^d(x)$ satisfy:
\begin{align}
	\textstyle\max_{y\in B^d(x)} v(I, y)
		&\textstyle\geq \sum_{r =1}^R v(I, y(x^*, x, \Psi_r))/R	
										\qquad\text{as $y(x^*, x, \Psi_r) \in B^d(x)$, $r\in [R]$}	\nonumber\\
		&\textstyle= 	\sum_{r =1}^R v(I, y(x^*, x, \Phi_r))/R	
										\qquad\text{by \cref{eq-vois--P=Q}}							\nonumber\\
		&\geq 			R^* \times v(I, x^*)/R +(R -R^*)\mathrm{wor}(I)/R							\label{eq-vois-S-apx}
\end{align}

Such solutions are therefore $R^*/R$-differential approximate if $x^*$ is optimal.

\subsection{Connection to combinatorial designs}\label{sec-vois-CD}

Similar to the previous section, one way to ensure that the arrays $\Psi$ and $\Phi$ satisfy \cref{eq-vois--P=Q} is to require that $\mu^\Psi -\mu^\Phi$ is balanced $k$-wise independent. 
Note that equality \cref{eq-vois--P=Q} follows from the equalities:
\begin{align}\label{eq-vois--P=Q-i}
	\textstyle\sum_{r =1}^R P_i(y(x^*, x, \Psi_r)_{J_i})	&\textstyle=\sum_{r =1}^R P_i(y(x^*, x, \Phi_r)_{J_i}), &i\in[m]
\end{align}

Consider then a constraint $C_i =P_i(x_{J_i})$ of $I$. A sufficient condition for $(\Psi, \Phi)$ to satisfy \cref{eq-vois--P=Q-i} at $i$ is that $C_i$ is evaluated on the same entries over the solution multisets $\left(y(x^*, x, \Psi_r)\,|\,r\in[R]\right)$ and $\left(y(x^*, x, \Phi_r)\,|\,r\in[R]\right)$. 
Let $H$ be the subsequence of indices of $J_i$ belonging to $\mathcal{J}(x^*, x)$, and $t$ the length of $H$ (note that $t$ can be zero). 
Recall that solutions $y(x^*, x, u)$, where $u\in\set{0, 1}^\kappa$, coincide with $x$ for coordinates of index outside $\mathcal{J}(x^*, x)$ and with either $x$ or $x^*$ depending on whether $u_c$ is 0 or 1 for coordinates of index $j_c\in\mathcal{J}(x^*, x)$.
Hence, the two multisubsets $\left(y(x^*, x, \Psi_r)_{J_i}\,|\,r\in[R]\right)$ and $\left(y(x^*, x, \Phi_r)_{J_i}\,|\,r\in[R]\right)$ of $\Sigma_q^{k_i}$ (where $k_i =|J_i|$) coincide {\em if and only if} the two multisubsets $(\Psi^H_r\,|\,r\in[R])$ and $(\Phi^H_r\,|\,r\in[R])$ of $\set{0, 1}^t$ coincide. Considering that $\card{H}\leq k_i\leq k$, this condition is indeed verified if $\mu^\Psi -\mu^\Phi$ is balanced $k$-wise independent. 
 %
We are therefore interested in such pairs $(\Psi, \Phi)$ of arrays that ideally maximize $\mu^\Phi(\mathbf{1})$, which we formalize below.

\begin{definition}\label{def-vois-CD}
Let $k\geq 1$, $d\geq k$, and $n\geq d$ be three integers.
Two arrays $\Psi$ and $\Phi$ with $d$ columns and Boolean coefficients form a {\em $(n, d)$-cover pair of arrays (for short, a $(n, d)$-CPA) of strength $k$} if they satisfy the conditions below:
\begin{enumerate}
	\item \label{it-thm-vois-Q} $\Phi$ contains at least 1 row of the form $(1, 1, \ldots, 1)$;
	\item \label{it-thm-vois-P} the number of 1's in each row of $\Psi$ is at most $d$;
	\item \label{it-thm-vois-P=Q} the function $\mu^\Psi -\mu^\Phi$ is balanced $k$-wise independent.
\end{enumerate} 

For two integers $R^* >0$ and $R\geq R^*$, we denote by $\Delta(R, R^*, n, d, k)$ the (possibly empty) set of $(n, d)$-CPAs of strength $k$ in which the row of 1's has multiplicity $R^*$ and the arrays have $R$ rows each.

Furthermore, we define $\delta(n, d, k)$ as the largest number $\delta$ for which there exist two integers $R^* >0$ and $R\geq R^*$ such that $R^*/R =\delta$ and $\Delta(R, R^*, n, d, k)\neq\emptyset$. 
\end{definition}

The preceding discussion allows us to establish the following connection between CPAs and the approximation guarantees that Hamming balls of radius $k$ might provide for $k$-CSPs:

\begin{theorem}\label{thm-vois-CD}
Let $q\geq 2$, $k\geq 2$, and $d\geq k$ be three integers, $I$ be an instance of $\mathsf{k\,CSP\!-\!q}$, and $x$ be a solution of $I$. Then for $I$, the highest differential ratio reached on $B^d(x)$ and on $\tilde{B}^d(x)$ is at least $\delta(n, d, k)$ and $\delta\left(\lfloor n(q -1)/q\rfloor, d, k\right)$, respectively.
\end{theorem}

\begin{proof}
Let $x^*$ be an optimal solution on $I$, and $\kappa$ be its Hamming distance from $x$. If $x^*\in B^d(x)$, then the largest differential ratio of a solution over $B^d(x)$ is 1, while $1\geq\delta(n, d, k)$. 
So we assume $\kappa >d$. Consider a CPA $(\Psi,\Phi)$ in $\Delta(R, R^*, \kappa, d, k)$, where $R\geq R^* >0$, and assume {\em w.l.o.g.} that the goal on $I$ is to maximize. 
By \cref{it-thm-vois-P=Q} of \cref{def-vois-CD}, $(\Psi,\Phi)$ satisfies \cref{{eq-vois--P=Q-i}} and hence \cref{eq-vois--P=Q}. This, together with \cref{it-thm-vois-Q,it-thm-vois-P} of \cref{def-vois-CD}, implies the inequality \cref{eq-vois-S-apx}. 

Thus, $\delta(\kappa, d, k)$ is a lower bound on the highest differential ratio reached on $B^d(x)$. 
Now, consider $\kappa_a =d_H(x +\mathbf{a}, x^*)$, for $a\in\Sigma_q$. Then for any $a\in\Sigma_q$, $\delta(\kappa_a, d, k)$ is also a lower bound on the highest differential ratio reached on $B^d(x +\mathbf{a})$. Therefore, we have:
$$\begin{array}{rll}
\displaystyle\max_{y\in\tilde{B}^d(x)}\frac{v(I, y) -\mathrm{wor}(I)}{\mathrm{opt}(I) -\mathrm{wor}(I)}
&\displaystyle =\max_{a =0}^{q -1}\left(
	\max_{y\in B^d(x +\mathbf{a})}\frac{v(I, y) -\mathrm{wor}(I)}{\mathrm{opt}(I) -\mathrm{wor}(I)}
\right)		&\geq\max_{a =0}^{q -1}\delta(\kappa_a, d, k)
\end{array}$$

We conclude by noting that the numbers $\delta(n, d, k)$ are naturally non-increasing in $n$. Indeed, extracting the first $\kappa\leq n$ columns of each array of an array pair from $\Delta(R, R^*, n, d, k)$ yields a CPA of $\Delta(R, S^*, \kappa, d, k)$, for some $S^*\geq R^*$. In particular, we have 
$\delta(\kappa, d, k)\geq\delta(n, d, k)$ and
$\max_{a =0}^{q -1}\delta(\kappa_a, d, k)
	=\delta(\min_{a =0}^{q -1}\kappa_a, d, k)
	\geq\delta(\sum_{a =0}^{q -1}\kappa_a/q, d, k)$,
while $\sum_{a =0}^{q -1}\kappa_a =(q -1)n$.
\end{proof}

\subsection{Estimation of numbers $\delta(n, d, k)$ and derived approximation guarantees}

It remains to obtain lower bounds on the numbers $\delta(n, d, k)$. In fact, this can be done by relating CPAs to ARPAs. Namely, an ARPA $(\Psi, \Phi)$ of some family $\Gamma(R, R^*, n, d, k)$ can be interpreted as a CPA of the family $\Delta(R, R^*, n, d, k)$ by interpreting the coefficients $M_r^j$ occurring in a column with index $j\in\Sigma_n$ of $\Psi$ or $\Phi$ as the binary relation $(M_r^j =j)$. 

\begin{proposition}\label{prop-vois-gamma-delta}
For all integers $k\geq 1$, $d\geq k$, and $n\geq d$, we have the inequality $\delta(n, d, k)\geq \gamma(n, d, k)$. 
\end{proposition}%

\begin{proof}[Proof (sketch)]
For each positive integer $n$, we associate a surjective mapping $\sigma_n$ that transforms arrays with $n$ columns and entries from $\Sigma_n$ into arrays with $n$ columns and binary coefficients.  
We index the columns of the arrays from the former set using $\Sigma_n$ and those from the latter set using $[n]$. 
Let $M$ be an $R\times n$ array on $\Sigma_n$, where $R$ is some positive integer. Then $\sigma_n$ maps $M$ to the $R\times n$ array on $\set{0, 1}$ defined by:
\begin{align}\label{eq-delta-gamma-def}
\begin{array}{rll}
	\sigma_n(M)_r^j 	&=\left\{\begin{array}{rl}
							1	&\text{if $M_r^{j -1} =j -1$}\\
							0	&\text{otherwise}\\ 
						\end{array}\right.,		&r\in[R], j\in[n]
\end{array}
\end{align}

For all integers $R^* >0$ and $R\geq R^*$ such that a CPA $(\Psi,\Phi)\in\Gamma(R, R^*, n, d, k)$ exists, we have $(\sigma_n(\Psi), \sigma_n(\Phi))\in\Delta(R, R^*, n, d, k)$  (see \cref{sec-arpa2cpa} for a detailed proof), and hence $\Delta(R, R^*, n, d, k)\neq\emptyset$. This implies $\delta(n, d, k)\geq \gamma(n, d, k)$.
\end{proof}

We notably deduce from \cref{prop-vois-gamma-delta} for the case where $d =k$ that $\delta(n, k, k)$ is at least $\gamma(n, k, k)$, and thus from \cref{thm-gamma_qpk} and the inequality \cref{eq-T-UB} that $\delta(n, k, k)$ is at least $2(k!)/(2n -k)^k$. This bound, together with \cref{thm-vois-CD}, leads to the following structural approximation guarantees for $\mathsf{k\,CSP\!-\!q}$:

\begin{corollary}\label{cor-vois-dapx}
Let $q\geq 2$, $k\geq 2$, and $d\geq k$ be three integers, $I$ be an instance of $\mathsf{k\,CSP\!-\!q}$ with $n\geq k$ variables, and $x$ be a solution of $I$. Then for $I$, the largest differential ratio achieved on $B^d(x)$ and on $\tilde{B}^d(x)$ is at least:
$$\frac{2(k!)}{(2n -k)^k}\text{ and }\frac{2(k!)}{\left(2(q -1)n/q -k\right)^k}$$
\end{corollary}

\subsection{Approximation results for the instance diameter}

We observe that we can even obtain approximation guarantees for the instance diameter, which is a stronger result, provided that the rows of the array $\Phi$, except for its rows of 1's, contain at most $d$ 1's.

\begin{definition}\label{def-vois-CD-diam}
Let $k\geq 1$, $d\geq k$, and $n\geq d$ be three integers.
For two integers $R^* >0$ and $R\geq R^*$, we define $\bar\Delta(R, R^*, n, d, k)$ as the (possibly empty) restriction of $\Delta(R, R^*, n, d, k)$ to pairs of arrays whose rows each contain either $n$ or at most $d$ 1's. 

Furthermore, the largest number $\delta$ for which there exist two integers $R^* >0$ and $R\geq R^*$ such that $R^*/R =\delta$ and $\bar\Delta(R, R^*, n, d, k)\neq\emptyset$ is denoted by $\bar\delta(n, d, k)$.
\end{definition}

\begin{theorem}\label{thm-vois-CD-diam}
Let $q\geq 2$, $k\geq 2$, and $d\geq k$ be three integers, $I$ be an instance of $\mathsf{k\,CSP\!-\!q}$, and $x$ be a solution of $I$. Then for $I$, the ratio of the maximum difference between two solution values over $B^d(x)$ and over $\tilde{B}^d(x)$ to $\abs{\mathrm{opt}(I) -\mathrm{wor}(I)}$ is at least:
$$	\frac{1}{2/\bar\delta(n, d, k) -1}\text{ and }
	\frac{1}{2/\bar\delta\left(\lfloor n(q -1)/q\rfloor, d, k\right) -1}
$$ 
\end{theorem}

\begin{proof}
Let $B =B^d(x)$. We denote by $x^*$ and $x_*$ a best and a worst solution on $I$, and by $\kappa^*$ and $\kappa_*$ their Hamming distance to $x$. We assume {\em w.l.o.g.} that the goal on $I$ is to maximize. 
If $x_*\in B$, then the maximum difference between two solution values over $B$ is the expression $\max_{y\in B} v(I, y) -\mathrm{wor}(I)$, which, according to \cref{thm-vois-CD}, is at least a fraction $\delta(n, d, k)$ of the diameter of $I$. 
Now, by definition, the quantities $\delta(n, d, k)$ and $\bar\delta(n, d, k)$ satisfy $1\geq\delta(n, d, k)\geq\bar\delta(n, d, k)$, and hence $\delta(n, d, k)\geq 1/(2/\bar\delta(n, d, k) -1)$. 
Symmetrically, if $x^*\in B$, then $\max_{y\in B} v(I, y) -\min_{y\in B} v(I, y) =\mathrm{opt}(I) -\min_{y\in B} v(I, y)\geq\delta(n, d, k)\times(\mathrm{opt}(I) -\mathrm{wor}(I))$. 

Therefore, for the rest of the proof, we assume that $\kappa_* >d$ and $\kappa^* >d$, which means that neither $x^*$ nor $x_*$ belongs to $B$. 
 %
Consider then two CPAs $(\Psi^*,\Phi^*)\in\bar\Delta(R, R^*, \kappa^*, d, k)$ and $(\Psi_*,\Phi_*)\in\bar\Delta(S, S_*, \kappa_*, d, k)$, where $(R, S)\geq (R^*, S_*) >(0, 0)$. Since the rows of $\Phi^*$ contain either exactly $\kappa^*$ or at most $d$ 1's, solutions $y(x^*, x, \Phi^*_r)$ that do not coincide with $x^*$ all belong to $B$. 
Accordingly, we can replace $\mathrm{wor}(I)$ in inequality \cref{eq-vois-S-apx} with the minimum solution value over $B$, giving the strengthened inequality:
\begin{align}\label{eq-vois-S-diam+}
\textstyle\max_{y\in B}v(I, y)	
	&\textstyle\geq R^*/R\times v(I, x^*) +(1 -R^*/R)\times\min_{y\in B}v(I, y)
\end{align}

By a symmetric argument, we can derive from $(\Psi_*,\Phi_*)$ the inequality:
\begin{align}\label{eq-vois-S-diam-}
\textstyle\min_{y\in B}v(I, y)	
	&\textstyle\leq S_*/S\times v(I, x_*) +(1 -S_*/S)\times\max_{y\in B}v(I, y)
\end{align}

By subtracting $R^*/R\times\cref{eq-vois-S-diam-}$ from $S_*/S\times\cref{eq-vois-S-diam+}$, we obtain the following lower bound on the ratio of the maximum difference between two solution values over $B$ to the instance diameter: 
\begin{align}\label{eq-diam}
\frac{\max_{y\in B}v(I, y) -\min_{y\in B}v(I, y)}{v(I, x^*) -v(I, x_*)}	&\geq\frac{1}{R/R^* +S/S_* -1}
\end{align}

To establish the approximation guarantee over $B$, we observe that \cref{eq-diam} holds in particular with $R/R^* +S/S_* =2R/R^*$ when $(\Psi^*,\Phi^*)$ and $(\Psi_*,\Phi_*)$ are both obtained by extracting the first columns of the arrays of a CPA realizing $\bar\delta(n, d, k)$.
For the approximation guarantee over $\tilde{B}^d(x)$, we consider a maximizer $x^+$ and a minimizer $x^-$ of $v(I, .)$ on $\tilde{B}^d(x)$. If either $x^- =x_*$ or $x^+=x^*$, then by \cref{thm-vois-CD}, the lower bound announced for the ratio of $v(I, x^+) -v(I, x^-)$ to $\abs{\mathrm{opt}(I) -\mathrm{wor}(I)}$ is valid. 

So we assume that $x^*, x_*\notin\tilde{B}^d(x)$. Let $b, c\in\Sigma_q$, $\kappa^*_b =d_H(x^*, x +\mathbf{b})$, and $\kappa_{c*} =d_H(x^*, x +\mathbf{c})$. We now suppose that $(\Psi^*,\Phi^*)\in\bar\Delta(R, R^*, \kappa^*_b, d, k)$ and $(\Psi_*,\Phi_*)\in\bar\Delta(S, S_*, \kappa_{c*}, d, k)$, in which case relations \cref{eq-vois-S-diam+} and \cref{eq-vois-S-diam-} hold for $B =B^d(x +\mathbf{b})$ and $B =B^d(x +\mathbf{c})$, respectively. 
Considering the inequalities: 
$$\textstyle v(I, x^+)	\geq\max_{y\in B^d(x +\mathbf{b})\cup B^d(x +\mathbf{c})}v(I, y)	
			\geq\min_{y\in B^d(x +\mathbf{b})\cup B^d(x +\mathbf{c})}v(I, y)	\geq v(I, x^-)$$
we deduce from these relations that solutions $x^+$ and $x^-$ satisfy:
\begin{align}
v(I, x^+)	&\geq R^*/R\times\mathrm{opt}(I) +(1 -R^*/R)\times v(I, x^-)	\label{eq-vois-S-diam+2}\\
v(I, x^-)	&\leq S_*/S\times\mathrm{wor}(I) +(1 -S_*/S)\times v(I, x^+)	\label{eq-vois-S-diam-2}
\end{align}

As we previously derived \cref{eq-diam} from $\cref{eq-vois-S-diam+}$ and $\cref{eq-vois-S-diam-}$, we can combine $\cref{eq-vois-S-diam+2}$ and $\cref{eq-vois-S-diam-2}$ so as to obtain the inequality: 
\begin{align}\label{eq-diam2}
\left(v(I, x^+) -v(I, x^-)\right)/\left(\mathrm{opt}(I) -\mathrm{wor}(I)\right)	&\geq 1/(R/R^* +S/S_* -1)
\end{align}

To conclude, we first note that we can choose $b$ and $c$ in such a way that both $\kappa^*_b$ and $\kappa_{c*}$ are at most $\lfloor n(q -1)/q\rfloor$. Then, we observe that the inequality \cref{eq-diam2} holds in particular if $(\Psi^*,\Phi^*)$ and $(\Psi_*,\Phi_*)$ are both extracted from a CPA realizing $\bar\delta(\lfloor n(q -1)/q\rfloor, d, k)$.
\end{proof}

To derive concrete approximation results from \cref{thm-vois-CD-diam}, we need an estimate for the numbers $\bar\delta(n, d, k)$. Again, the ARPAs constructed in the previous section provide such an estimate. We focus on the case where $d =k$ and $n >k$.
We can use \cref{algo-Gamma} to generate a $(n, k)$-ARPA of strength $k$. 
We denote by $R_n$ the number of rows in the arrays of this ARPA, which we know to be $(T(n, k) +1)/2$. We then consider the CPA $(\Psi,\Phi)$ obtained by applying the transformation $\sigma_n$ of \cref{prop-vois-gamma-delta} to the arrays of this pair. \Cref{tab-Delta-rec} shows this CPA when $(k, n)\in\set{(2, 6), (3, 5)}$. According to the proof of \cref{prop-vois-gamma-delta}, $\Psi$ and $\Phi$ form a $(n, k)$-CPA of strength $k$, with $R_n$ rows each, including (in $\Phi$) a single row of 1's.
Upon closer examination, $(\Psi,\Phi)$ is actually an element of $\bar\Delta(R_n, n, 1, k, k)$. In fact, a more detailed analysis (see \cref{sec-vois-id}) reveals that $(\Psi,\Phi)$ can be described as follows:	
\begin{itemize}
	\item the word of 1's occurs exactly once as a row in $\Phi$;
	\item for all integers $a\in\set{0 ,\ldots, k}$ such that $a\equiv k\bmod{2}$, every word $u\in\set{0, 1}^n$ containing $a$ 1's occurs exactly $\binom{n -1 -a}{k -a}$ times as a row in $\Psi$;
	\item for all integers $a\in\set{0 ,\ldots, k}$ such that $a\not\equiv k\bmod{2}$, every word $u\in\set{0, 1}^n$ containing $a$ 1's occurs exactly $\binom{n -1 -a}{k -a}$ times as a row in $\Phi$;
	\item no other word of $\set{0, 1}^n$ occurs in either $\Psi$ or $\Phi$.
\end{itemize}

This CPA provides for $\bar\delta(n, k, k)$ the lower bound of $1/R_n$. Considering the inequalities $1/(2/\bar\delta(n, k, k) -1) >\bar\delta(n, k, k)/2$ and $1/R_n\leq 2(k!)/(2n -k)^k$ (by \cref{eq-T-UB}), it thus implies together with \cref{thm-vois-CD-diam} the following approximation guarantees for the diameter of instances of $k$-CSPs:

\begin{table}\begin{center}{\footnotesize
\setlength\arraycolsep{2pt}
$\begin{array}{cp{4pt}c}
2/(T(6, 2) +1) =1/25		&&2/(T(5, 3) +1) =1/25\\[6pt]
\begin{array}{cp{2pt}c}
\begin{array}{cc|c|c|c|c|}
	\Psi^1	&\Psi^2					&\Psi^3					&\Psi^4					&\Psi^5					&\Psi^6\\
\hline
	1		&1						&0						&0						&0						&0\\
\cline{1-2}
	1		&\multicolumn{1}{c}{0}	&1						&0						&0						&0\\
	0		&\multicolumn{1}{c}{1}	&1						&0						&0						&0\\
	0		&\multicolumn{1}{c}{0}	&0						&0						&0						&0\\
\cline{1-3}
	1		&\multicolumn{1}{c}{0}	&\multicolumn{1}{c}{0}	&1						&0						&0\\
	0		&\multicolumn{1}{c}{1}	&\multicolumn{1}{c}{0}	&1						&0						&0\\
	0		&\multicolumn{1}{c}{0}	&\multicolumn{1}{c}{1}	&1						&0						&0\\
	0		&\multicolumn{1}{c}{0}	&\multicolumn{1}{c}{0}	&0						&0						&0\\
	0		&\multicolumn{1}{c}{0}	&\multicolumn{1}{c}{0}	&0						&0						&0\\
\cline{1-4}
	1		&\multicolumn{1}{c}{0}	&\multicolumn{1}{c}{0}	&\multicolumn{1}{c}{0}	&1						&0\\
	0		&\multicolumn{1}{c}{1}	&\multicolumn{1}{c}{0}	&\multicolumn{1}{c}{0}	&1						&0\\
	0		&\multicolumn{1}{c}{0}	&\multicolumn{1}{c}{1}	&\multicolumn{1}{c}{0}	&1						&0\\
	0		&\multicolumn{1}{c}{0}	&\multicolumn{1}{c}{0}	&\multicolumn{1}{c}{1}	&1						&0\\
	0		&\multicolumn{1}{c}{0}	&\multicolumn{1}{c}{0}	&\multicolumn{1}{c}{0}	&0						&0\\
	0		&\multicolumn{1}{c}{0}	&\multicolumn{1}{c}{0}	&\multicolumn{1}{c}{0}	&0						&0\\
	0		&\multicolumn{1}{c}{0}	&\multicolumn{1}{c}{0}	&\multicolumn{1}{c}{0}	&0						&0\\
\cline{1-5}
	1		&\multicolumn{1}{c}{0}	&\multicolumn{1}{c}{0}	&\multicolumn{1}{c}{0}	&\multicolumn{1}{c}{0}	&1\\
	0		&\multicolumn{1}{c}{1}	&\multicolumn{1}{c}{0}	&\multicolumn{1}{c}{0}	&\multicolumn{1}{c}{0}	&1\\
	0		&\multicolumn{1}{c}{0}	&\multicolumn{1}{c}{1}	&\multicolumn{1}{c}{0}	&\multicolumn{1}{c}{0}	&1\\
	0		&\multicolumn{1}{c}{0}	&\multicolumn{1}{c}{0}	&\multicolumn{1}{c}{1}	&\multicolumn{1}{c}{0}	&1\\
	0		&\multicolumn{1}{c}{0}	&\multicolumn{1}{c}{0}	&\multicolumn{1}{c}{0}	&\multicolumn{1}{c}{1}	&1\\
	0		&\multicolumn{1}{c}{0}	&\multicolumn{1}{c}{0}	&\multicolumn{1}{c}{0}	&\multicolumn{1}{c}{0}	&0\\
	0		&\multicolumn{1}{c}{0}	&\multicolumn{1}{c}{0}	&\multicolumn{1}{c}{0}	&\multicolumn{1}{c}{0}	&0\\
	0		&\multicolumn{1}{c}{0}	&\multicolumn{1}{c}{0}	&\multicolumn{1}{c}{0}	&\multicolumn{1}{c}{0}	&0\\
	0		&\multicolumn{1}{c}{0}	&\multicolumn{1}{c}{0}	&\multicolumn{1}{c}{0}	&\multicolumn{1}{c}{0}	&0\\
\hline
\end{array}&&\begin{array}{cc|c|c|c|c|}
	\Phi^1	&\Phi^2					&\Phi^3					&\Phi^4					&\Phi^5					&\Phi^6\\
\hline
	1		&1						&1						&1						&1						&1\\
\cline{1-2}
	1		&\multicolumn{1}{c}{0}	&0						&0						&0						&0\\
	0		&\multicolumn{1}{c}{1}	&0						&0						&0						&0\\
	0		&\multicolumn{1}{c}{0}	&1						&0						&0						&0\\
\cline{1-3}
	1		&\multicolumn{1}{c}{0}	&\multicolumn{1}{c}{0}	&0						&0						&0\\
	0		&\multicolumn{1}{c}{1}	&\multicolumn{1}{c}{0}	&0						&0						&0\\
	0		&\multicolumn{1}{c}{0}	&\multicolumn{1}{c}{1}	&0						&0						&0\\
	0		&\multicolumn{1}{c}{0}	&\multicolumn{1}{c}{0}	&1						&0						&0\\
	0		&\multicolumn{1}{c}{0}	&\multicolumn{1}{c}{0}	&1						&0						&0\\
\cline{1-4}
	1		&\multicolumn{1}{c}{0}	&\multicolumn{1}{c}{0}	&\multicolumn{1}{c}{0}	&0						&0\\
	0		&\multicolumn{1}{c}{1}	&\multicolumn{1}{c}{0}	&\multicolumn{1}{c}{0}	&0						&0\\
	0		&\multicolumn{1}{c}{0}	&\multicolumn{1}{c}{1}	&\multicolumn{1}{c}{0}	&0						&0\\
	0		&\multicolumn{1}{c}{0}	&\multicolumn{1}{c}{0}	&\multicolumn{1}{c}{1}	&0						&0\\
	0		&\multicolumn{1}{c}{0}	&\multicolumn{1}{c}{0}	&\multicolumn{1}{c}{0}	&1						&0\\
	0		&\multicolumn{1}{c}{0}	&\multicolumn{1}{c}{0}	&\multicolumn{1}{c}{0}	&1						&0\\
	0		&\multicolumn{1}{c}{0}	&\multicolumn{1}{c}{0}	&\multicolumn{1}{c}{0}	&1						&0\\
\cline{1-5}
	1		&\multicolumn{1}{c}{0}	&\multicolumn{1}{c}{0}	&\multicolumn{1}{c}{0}	&\multicolumn{1}{c}{0}	&0\\
	0		&\multicolumn{1}{c}{1}	&\multicolumn{1}{c}{0}	&\multicolumn{1}{c}{0}	&\multicolumn{1}{c}{0}	&0\\
	0		&\multicolumn{1}{c}{0}	&\multicolumn{1}{c}{1}	&\multicolumn{1}{c}{0}	&\multicolumn{1}{c}{0}	&0\\
	0		&\multicolumn{1}{c}{0}	&\multicolumn{1}{c}{0}	&\multicolumn{1}{c}{1}	&\multicolumn{1}{c}{0}	&0\\
	0		&\multicolumn{1}{c}{0}	&\multicolumn{1}{c}{0}	&\multicolumn{1}{c}{0}	&\multicolumn{1}{c}{1}	&0\\
	0		&\multicolumn{1}{c}{0}	&\multicolumn{1}{c}{0}	&\multicolumn{1}{c}{0}	&\multicolumn{1}{c}{0}	&1\\
	0		&\multicolumn{1}{c}{0}	&\multicolumn{1}{c}{0}	&\multicolumn{1}{c}{0}	&\multicolumn{1}{c}{0}	&1\\
	0		&\multicolumn{1}{c}{0}	&\multicolumn{1}{c}{0}	&\multicolumn{1}{c}{0}	&\multicolumn{1}{c}{0}	&1\\
	0		&\multicolumn{1}{c}{0}	&\multicolumn{1}{c}{0}	&\multicolumn{1}{c}{0}	&\multicolumn{1}{c}{0}	&1\\
\hline
\end{array}
\end{array}&&\begin{array}{cp{2pt}c}
\begin{array}{ccc|c|c|}
	\Psi^1	&\Psi^2	&\Psi^3						&\Psi^4						&\Psi^5\\
\hline
	1		&1		&1							&0							&0\\
\cline{1-3}
	1		&1		&\multicolumn{1}{c}{0}		&1							&0\\
	1		&0		&\multicolumn{1}{c}{1}		&1							&0\\
	0		&1		&\multicolumn{1}{c}{1}		&1							&0\\

	1		&0		&\multicolumn{1}{c}{0}		&0							&0\\
	0		&1		&\multicolumn{1}{c}{0}		&0							&0\\
	0		&0		&\multicolumn{1}{c}{1}		&0							&0\\

	0		&0		&\multicolumn{1}{c}{0}		&1							&0\\
\cline{1-4}
	1		&1		&\multicolumn{1}{c}{0}		&\multicolumn{1}{c}{0}		&1\\
	1		&0		&\multicolumn{1}{c}{1}		&\multicolumn{1}{c}{0}		&1\\
	1		&0		&\multicolumn{1}{c}{0}		&\multicolumn{1}{c}{1}		&1\\
	0		&1		&\multicolumn{1}{c}{1}		&\multicolumn{1}{c}{0}		&1\\
	0		&1		&\multicolumn{1}{c}{0}		&\multicolumn{1}{c}{1}		&1\\
	0		&0		&\multicolumn{1}{c}{1}		&\multicolumn{1}{c}{1}		&1\\

	1		&0		&\multicolumn{1}{c}{0}		&\multicolumn{1}{c}{0}		&0\\
	1		&0		&\multicolumn{1}{c}{0}		&\multicolumn{1}{c}{0}		&0\\
	0		&1		&\multicolumn{1}{c}{0}		&\multicolumn{1}{c}{0}		&0\\
	0		&1		&\multicolumn{1}{c}{0}		&\multicolumn{1}{c}{0}		&0\\
	0		&0		&\multicolumn{1}{c}{1}		&\multicolumn{1}{c}{0}		&0\\
	0		&0		&\multicolumn{1}{c}{1}		&\multicolumn{1}{c}{0}		&0\\
	0		&0		&\multicolumn{1}{c}{0}		&\multicolumn{1}{c}{1}		&0\\
	0		&0		&\multicolumn{1}{c}{0}		&\multicolumn{1}{c}{1}		&0\\

	0		&0		&\multicolumn{1}{c}{0}		&\multicolumn{1}{c}{0}		&1\\
	0		&0		&\multicolumn{1}{c}{0}		&\multicolumn{1}{c}{0}		&1\\
	0		&0		&\multicolumn{1}{c}{0}		&\multicolumn{1}{c}{0}		&1\\
\hline
\end{array}&&\begin{array}{ccc|c|c|}
	\Phi^1	&\Phi^2	&\Phi^3						&\Phi^4						&\Phi^5\\
\hline
	1		&1		&1							&1							&1\\
\cline{1-3}
	1		&1		&\multicolumn{1}{c}{0}		&0							&0\\
	1		&0		&\multicolumn{1}{c}{1}		&0							&0\\
	0		&1		&\multicolumn{1}{c}{1}		&0							&0\\

	1		&0		&\multicolumn{1}{c}{0}		&1							&0\\
	0		&1		&\multicolumn{1}{c}{0}		&1							&0\\
	0		&0		&\multicolumn{1}{c}{1}		&1							&0\\

	0		&0		&\multicolumn{1}{c}{0}		&0							&0\\
\cline{1-4}
	1		&1		&\multicolumn{1}{c}{0}		&\multicolumn{1}{c}{0}		&0\\
	1		&0		&\multicolumn{1}{c}{1}		&\multicolumn{1}{c}{0}		&0\\
	1		&0		&\multicolumn{1}{c}{0}		&\multicolumn{1}{c}{1}		&0\\
	0		&1		&\multicolumn{1}{c}{1}		&\multicolumn{1}{c}{0}		&0\\
	0		&1		&\multicolumn{1}{c}{0}		&\multicolumn{1}{c}{1}		&0\\
	0		&0		&\multicolumn{1}{c}{1}		&\multicolumn{1}{c}{1}		&0\\

	1		&0		&\multicolumn{1}{c}{0}		&\multicolumn{1}{c}{0}		&1\\
	1		&0		&\multicolumn{1}{c}{0}		&\multicolumn{1}{c}{0}		&1\\
	0		&1		&\multicolumn{1}{c}{0}		&\multicolumn{1}{c}{0}		&1\\
	0		&1		&\multicolumn{1}{c}{0}		&\multicolumn{1}{c}{0}		&1\\
	0		&0		&\multicolumn{1}{c}{1}		&\multicolumn{1}{c}{0}		&1\\
	0		&0		&\multicolumn{1}{c}{1}		&\multicolumn{1}{c}{0}		&1\\
	0		&0		&\multicolumn{1}{c}{0}		&\multicolumn{1}{c}{1}		&1\\
	0		&0		&\multicolumn{1}{c}{0}		&\multicolumn{1}{c}{1}		&1\\

	0		&0		&\multicolumn{1}{c}{0}		&\multicolumn{1}{c}{0}		&0\\
	0		&0		&\multicolumn{1}{c}{0}		&\multicolumn{1}{c}{0}		&0\\
	0		&0		&\multicolumn{1}{c}{0}		&\multicolumn{1}{c}{0}		&0\\
\hline
\end{array}
\end{array}
\end{array}$
\caption{\label{tab-Delta-rec}Construction for $\bar\Delta\left((T(n, k) +1)/2, 1, n, k, k\right)$: illustration when $(k, n)\in\{(2, 6), (3, 5)\}$. These CPAs are obtained by applying the transformation $\sigma_n$ of \cref{prop-vois-gamma-delta} to the arrays of the ARPAs generated by \cref{algo-Gamma} on parameters (2, 6) and (3, 5).}
}\end{center}\end{table}

\begin{corollary}\label{cor-vois-diam}
Let $q\geq 2$, $k\geq 2$, and $d\geq k$ be three integers, $I$ be an instance of $\mathsf{k\,CSP\!-\!q}$ with $n\geq k$ variables, and $x$ be a solution of $I$. Then the ratio of the maximum difference between two solution values over $B^d(x)$ and over $\tilde{B}^d(x)$ to the diameter of $I$ is at least:
$$\frac{k!}{(2n -k)^k}\text{ and }\frac{k!}{\left(2(q -1)n/q -k\right)^k}$$
\end{corollary}

We emphasize that equality \cref{eq-vois--P=Q}, if $(\Psi,\Phi)$ is a CPA from some family $\bar\Delta(R, R^*, \kappa, d, k)$, allows to express $v(I, x^*)$ as a linear combination of solution values over $B^k(x)$. More generally, given an instance of a $k$-CSP and a solution $x$ of that instance, we can derive from $(\Psi,\Phi)$ an expression of the value of any solution $x^*$ at Hamming distance $\kappa$ from $x$ as a linear combination of solution values over $B^k(x)$. In particular, the specific array pairs used to establish \cref{cor-vois-diam} also yield the following identity:

\begin{theorem}\label{thm-vois-id}
Let $q\geq 2$ and $k\geq 2$ be two integers, $I$ be an instance of $\mathsf{k\,CSP\!-\!q}$, and $x$ and $x^*$ be two solutions of $I$ that are at Hamming distance $\kappa >k$ from each other.

For each $h\in\set{0 ,\ldots, \kappa}$, we denote by $N^h(x^*, x)$ the set of vectors of $\Sigma_q^n$ that coincide with either $x^*$ or $x$ at each of their coordinates, and whose Hamming distance from $x$ is $h$. Formally:
\begin{align}\nonumber
N^h(x^*, x)	&\textstyle:=
	\left\{y\in\set{x^*_1, x_1}\times\ldots\times\set{x^*_n, x_n}\,:\,d_H(x, y) =h\right\}
\end{align}
(In particular, $N^0(x^*, x) =\set{x}$ and $N^\kappa(x^*, x) =\set{x^*}$.) 
Then $v(I, x^*)$ can be expressed as a linear combination of solution values over $N^0(x^*, x)\cup\ldots\cup N^k(x^*, x)$:
\begin{align}\label{eq-vois-id}
	v(I, x^*)	&=\textstyle\sum_{h =0}^k (-1)^{k -h} \binom{\kappa -1 -h}{k -h} \sum_{y \in N^h(x^*, x)} v(I, y)
\end{align}%
\end{theorem}%

\subsection{Concluding remarks}

Hamming balls of radius $k$ therefore approximate the instance diameter of $\mathsf{k\,CSP\!-\!q}$ by a factor $\Omega(1/n^k)$. For $\mathsf{2\,CSP\!-\!2}$, we observe that this bound, $\Omega(1/n^2)$, is relatively weak compared to the (constant) approximation guarantee of $4/\pi -1 >0.273$ for the instance diameter \cite{N98}. 
Nevertheless, when $q$ is a power of 2, the guarantee obtained is, in dense instances, comparable to the gain approximability bound of $\Omega(1/m)$ implied by \cite{HV04}. 
	
The remainder of this section is devoted to assessing the tightness of the analyses performed.

\smallskip
{\bf{\em Hamming balls with radius $k$.}}
 %
Regarding the approximation guarantees offered by Hamming balls of fixed radius $d\geq k$, we cannot expect ratios better than $\Omega(1/n^k)$. We denote by $J^{q, k}_n$ the instance of $\mathsf{CSP(\set{AllZero^{k, q}})}$ that considers all the $k$-ary constraints that can be formulated on a set of $n$ variables, given three positive integers $q$, $k$, and $n \geq k$.
For any $d\in\set{0 ,\ldots, n}$, every vector with exactly $d$ zero coordinates satisfies $\binom{d}{k}$ of the constraints. In particular, we have $\mathrm{opt}(J^{q, k}_n) =\binom{n}{k}$ (the vector of all zeros satisfies all constraints) and $\mathrm{wor}(J^{q, k}_n) =0$ ({\em e.g.}, the vector of all ones satisfies no constraint). Furthermore, for $d\in\set{k ,\ldots, n}$, the maximum solution value over $B^d(\mathbf{1})$ is equal to $\binom{d}{k}$. Hence, the highest differential ratio achieved on  $B^d(\mathbf{1})$ is: 
$$\begin{array}{rll}
	\frac{\binom{d}{k} -0}{\binom{n}{k} -0}
	&=\displaystyle		\frac{d(d -1) \ldots (d -k +1)}{n(n -1) \ldots (n -k +1)}
	&\sim\displaystyle	\frac{k!\binom{d}{k}}{n^k}
\end{array}$$ 

Since $B^d(\mathbf{1})$ contains $\mathbf{1}$, which is a worst solution, this ratio coincides with the ratio of $\max_{y\in B^d(\mathbf{1})} v(J^{q, k}_n, y) -\min_{y\in B^d(\mathbf{1})} v(J^{q, k}_n, y)$ to the diameter of $J^{q, k}_n$. 
When $d =k$, this ratio is asymptotically a factor $2^{k -1}$ and $2^k$ of the lower bounds given by \cref{cor-vois-dapx,cor-vois-diam}, respectively, for the differential approximation of $\mathrm{opt}(I)$ and the approximation of the instance diameter.

\label{sec-vois-tight}	
We now examine the approximation guarantees obtained for the union of the shifts by $\mathbf{a}$, $a\in\Sigma_q$ of Hamming balls of radius $k$. We analyze the ratios reached on the instances $I^{q, k}_n$ (where $n\in\mathbb{N}\backslash\set{0}$) of $\mathsf{Max\,CSP(\set{AllEqual^{k, q}})}$ introduced in \cref{sec-E-tight}. Since the function $AllEqual^{k, q}$ remains invariant under uniform shifts of its inputs, then for any instance $I^{q, k}_n$, a solution $y$ is optimal over some Hamming ball $B^k(x)$ {\em if and only if} it is optimal over $\tilde{B}^k(x)$. 
Consider three positive integers $n$, $q$, and $k$. Given a partition of $[qn]$ into $q$ natural numbers $n_0, \ldots, n_{q -1}$, any vector of $\Sigma_q^{qn}$ with $n_a$ coordinates equal to $a$, $a\in\Sigma_q$ satisfies $\textstyle\sum_{a =0}^{q -1}\binom{n_a}{k}$ constraints. Recall that $\mathrm{opt}(I^{q, k}_n) =\binom{qn}{k}$ and $\mathrm{wor}(I^{q, k}_n) =q\times\binom{n}{k}$.
Moreover, let $x_*$ be a vector with $n$ coordinates equal to $a$ for each $a\in\Sigma_q$. Then, one can easily check (see \cref{sec-I^qk_n} for a complete argument) that for $d\in\set{0 ,\ldots, n}$, the maximum solution value over $\tilde{B}^d(x_*)$ is equal to $\binom{n +d}{k} +\binom{n -d}{k} +(q -2)\binom{n}{k}$. 
The highest differential ratio achieved on $\tilde{B}^d(x_*)$ is therefore:
\begin{align}\label{eq-tildeB^d}
\frac{\binom{n +d}{k} +\binom{n -d}{k} -2\binom{n}{k}}{\binom{qn}{k} -q\times\binom{n}{k}}
	&\sim\frac{2d^2 \binom{k}{2}q}{q^{k -1} -1}\times\frac{1}{(qn)^2}
\end{align}

Now, according to \cref{cor-vois-dapx}, on $I^{q, k}_n$, for any $x\in\Sigma_q^{qn}$, the largest differential ratio of a solution over $\tilde{B}^k(x)$ is at least: 
$$\begin{array}{rll}
	\displaystyle		\frac{2\times k!}{(2(q -1)(qn)/q -k)^k}
	&=\displaystyle		\frac{k!}{2^{k -1}((q -1)n -k/2)^k}
	&\sim\displaystyle	\frac{q^k\times k!}{2^{k -1}(q -1)^k}\times\frac{1}{(qn)^k}
\end{array}$$ 

Thus, when $k =2$, the maximum differential ratio over $\tilde{B}^2(x_*)$ is asymptotically a factor of $8(q -1)/q$ of the lower bound on this ratio given by \cref{cor-vois-dapx}. 
If $k$ and $n$ are constant integers while $q$ can be arbitrarily large, then both this ratio and the lower bound given by \cref{cor-vois-dapx} are in $\Theta(1/q^k)$. 

\smallskip 
{\bf{\em Hamming balls with radius 1.}}
\label{sec-vois-B1_tight}	
The differential approximation guarantees obtained for $\mathsf{E2\,CSP(\mathcal{I}_2^1)}$ --- or, equivalently, $\mathsf{E2\,Lin\!-\!2}$ (see \ref{sec-ap-bib}) --- in \cref{sec-B1} are either tight or asymptotically tight. For this problem, we can derive from \cref{pty-I_q^(k-1)-B^1,thm-E-OA,thm-E-2CSPs} the lower bounds of $1/(2\lceil\nu/2\rceil)$ for the differential ratio of local optima with respect to $B^1$, and of $2/(\lceil\nu/2\rceil\times n)$ for the best differential ratio on an arbitrary Hamming ball of radius 1, where we recall that $\nu$ and $n$ denote the strong coloring number and the number of variables of the instance under consideration.

When $q =k =2$, $I^{2, 2}_n$ is an instance of $\mathsf{E2\,Lin\!-\!2}$, and the left expression in \cref{eq-tildeB^d} evaluates to $d^2/n^2$, which is equal to $1/n^2$ when $d =1$. Considering that in $I^{2, 2}_n$, the strong coloring number and the number of variables both are $2n$, we conclude that the lower bound given in \cref{sec-vois-B1_tight} for the maximum differential ratio over Hamming balls of radius 1 is tight. 

On $I^{2, 2}_n$, local search with respect to $B^1$ will return the optimal solution either $\mathbf{0}$ or $\mathbf{1}$. Therefore, instead of $I^{2, 2}_n$, we consider the instance $\tilde{I}_n$ of $\mathsf{E2\,Lin\!-\!2}$ obtained from $I^{2, 2}_n$ by removing for each $j\in[n]$ the constraint $(x_{2j -1} =x_{2j})$. This instance remains trivially satisfiable ({\em e.g.}, by the zero vector). 
Furthermore, the solution whose non-zero coordinates are the odd-index coordinates is a local optimum with respect to $\tilde{B}^1$. It is easy to verify that the differential ratio of this solution is $1/(2\lceil n/2\rceil)$ (see \cref{sec-tilde{I}_n} for a complete argument). Remembering that $\tilde{I}_n$ is $2n$-partite, this ratio is asymptotically a factor of $2$ of the lower bound our analysis provides for it.

\section{Conclusions and future research directions}\label{sec-conc}
{\bf{\em Combinatorial designs and CSPs.}}
We used different combinatorial structures to analyze CSPs, which allowed us to obtain original approximability bounds for $k$-CSPs. While our investigations span various contexts, the underlying approach remains the same across these contexts. First, we identify partitions $\mathcal{V} =(V_1 ,\ldots, V_\nu)$ of $[n]$ and multisubsets $M$ of $\Sigma_q^\nu$ that have the right properties for the purpose at hand. Then, we consider solutions $y(\mathcal{V}, x, M_r)$ of the form $(x_{V_1} +\mathbf{M_r^1} ,\ldots, x_{V_\nu} +\mathbf{M_r^\nu})$ defined relative to a given solution $x$. 
If no further restrictions are imposed on the CSPs beyond the arity of their constraints, balanced $k$-wise independent functions become the natural tool for their analysis.
We emphasize that in all cases considered, it is sufficient to establish the existence of the pairs $(\mathcal{V}, M)$ that define the multisets of solutions that support the argument, without necessarily making them explicit.

\smallskip
In our opinion, this work highlights the power of combinatorial structures for studying the differential approximability of $k$-CSPs.
Nevertheless, further investigation is needed to assess the tightness of the ratios obtained.
One question concerns the reduction from $\mathsf{k\,CSP\!-\!q}$ to $\mathsf{k\,CSP\!-\!p}$:
on instances of $\mathsf{k\,CSP\!-\!q}$, how close is the differential ratio of the best solutions among those whose coordinates take at most $p$ different values to $\gamma(q, p, k)$ (see relation \cref{eq-reduc-opt})?
Another question concerns the maximum differential ratio over the union of shifts by $(a ,\ldots, a)$, $a\in\{0 ,\ldots, q -1\}$ of Hamming balls with radius $k$, when $k\geq 3$. Namely, we obtain the lower bound of $\Omega(1/n^k)$ for this ratio, and ask whether this result can be improved or not. 
The same question arises when $k\geq 4$ or $q\geq k =3$ for the average differential ratio, for which we have established the lower bounds of $\Omega(1/n^b)$ with $b =\lfloor k/2\rfloor$ for $q =2$ and $b =k -\lceil\log_{\Theta(q)}k\rceil$ for larger integers $q$.

We also emphasize that the existence of a reduction from $\mathsf{k\,CSP\!-\!q}$ to $\mathsf{k\,CSP\!-\!p}$ (where $q >p$) that preserves the differential approximation ratio within some constant multiplicative factor remains completely unresolved when $k >p$. Furthermore, for an integer $k\geq 3$, the question of whether $\mathsf{k\,CSP\!-\!q}$ is approximable within some constant differential approximation ratio remains open for all $q\geq 3$ if $k =3$, otherwise for all $q\geq 2$. 
Nevertheless, we reduce this question to the case where $q\leq k$.

\medskip
{\bf{\em The solution landscape of CSPs.}}
Neighborhood analyses inherently establish that $\rho$-approximate solutions can be found throughout the solution set for a given $\rho$. 
In particular, our analysis for $\mathsf{k\,CSP\!-\!q}$ suggests that $1/n^k$ of the solutions provide a differential approximation guarantee of $\Omega(1/n^k)$. 
For $\mathsf{Ek\,CSP(\mathcal{I}_q^{k -1})}$, $1/(k q)$ of the solutions achieve a differential ratio that is a factor $\Omega(1/n)$ of the average differential ratio. 
In the very special cases of $\mathsf{E(2k +1)Lin\!-\!2}$ and $\mathsf{CSP(\mathcal{O}_q)}$, one half and $1/q$ of the solutions are $1/2$ and $1/q$-differential approximate, respectively.

Thanks to the methodology used, the results obtained for the average differential ratio also provide insights into the distribution of solution values. Namely, if $\mathcal{V}$ and $M$ denote the partition of $[n]$ and the array used to establish the approximation guarantee, then for any $x$, among the solutions $y(\mathcal{V}, x, M_r)$, $r\in[R]$, at least one achieves at least the average differential ratio. 
Therefore, the derived lower bounds for the average differential ratio not only indicate that $\mathbb{E}_X[v(I, X)]$ realizes a certain differential ratio $\rho$, but also suggest that $\rho$-differential approximate solutions are distributed throughout the solution set. 
Thus, our analysis of the average differential ratio teaches us that picking a solution uniformly at random yields a solution with an expected differential ratio of $\Omega(1/\nu^b)$, and that a solution that achieves such a ratio can be found in {\em some} $O(1/\nu^b)$-cardinality neighborhood of every solution.

\medskip
{\bf{\em The average differential ratio.}}
We believe that the average differential ratio has the potential to provide new insights into CSPs. 
In our analysis, we only considered a few characteristics of the input instance, namely: 
	the strong chromatic number of its primary hypergraph, 
	the possible restriction of the functions expressing its constraints to the --- general enough --- function families $\mathcal{E}_q$ and $\mathcal{I}_q^t$, or to the --- rather restrictive --- function family $\mathcal{O}_q$, 
	and the maximum constraint arity. 
Therefore, a next step would be to identify hypergraph structures and function properties that allow the construction of partition-based solution multisets of low cardinality that satisfy \cref{eq-E-M=E}. 
More generally, we think of the characterization of function families $\mathcal{F}$ as the set of submodular functions for which $\mathsf{Max\,CSP(\mathcal{F})}$ or $\mathsf{Min\,CSP(\mathcal{F})}$ admits a constant lower bound on the average differential ratio.

Beyond these aspects, the average differential ratio should be considered as an indicator of the computational complexity of combinatorial optimization problems. In particular, its relationship with approximation measures deserves to be studied. For example, unless $\mathbf{P}\neq\mathbf{NP}$, the diameter of $\mathsf{3\,Sat}$ instances is inapproximable within any constant ratio, \cite{EP05}. The authors of \cite{EP05} were able to derive this result from the hardness result of \cite{H97} for $\mathsf{E3\,Lin\!-\!2}$, based precisely on the fact that, for this particular CSP, the average differential ratio is in $O(1)$.

\medskip
{\bf{\em Combinatorial designs.}}
This work raises new questions about orthogonal arrays and difference schemes while introducing new families of combinatorial structures. 
 %
First, the proposed estimation of the average differential ratio involves an unusual criterion for the construction of orthogonal arrays and difference schemes of given strength and number of factors. Specifically, the analysis in \cref{sec-average} suggests searching for such arrays that maximize their highest frequency (rather than minimizing the number of their rows). 
Recent work has explored this direction for orthogonal arrays of strength 2 \cite{CSV19}.

Second, the reduction of $q$-ary CSPs to $p$-ary CSPs, together with the neighborhood analysis, suggests further investigation of the {\em alphabet reduction} and the {\em cover} pairs of arrays introduced in \cref{sec-reduc,sec-vois}. 
We recently proved that the values of $\gamma(q, p, k)$, $\delta(q, p, k)$, and $\bar\delta(q, p, k)$ all coincide for all triples $(q, p, k)$, and determined their exact value in the case where $k\in\set{p, 2, 1}$\footnote{This result is available in the arXiv preprint {\em J.-F. Culus, S. Toulouse, Optimizing alphabet reduction pairs of arrays, arXiv:2406.10930 (2024).}\label{CT24-reg}}. 
In contrast, when $q >p >k >2$, the naive lower bound of $\delta(q -p +k, k, k)$ --- which we know to be strictly less than $\delta(q, p, k)$\textsuperscript{\ref{CT24-reg}} --- is the only estimate of $\delta(q, p, k)$ that we know. In particular, a thorough study of $\delta(q, p, k)$ for such sets $(q, p, k)$ of parameters would improve our analysis of how well solutions of $\mathsf{k\,CSP\!-\!p}$ instances allow to approximate $\mathsf{k\,CSP\!-\!q}$. 

In addition, we are particularly interested in {\em relaxed} ARPAs and their corresponding numbers $\gamma_E(q, p, k)$, especially in the case when $p =k$. In particular, relaxed $(q, 2)$-ARPAs of strength 2 have the potential to increase the differential approximability of $\mathsf{2\,CSP(\mathcal{E}_q)}$ by a factor of $\Theta(q)$. Although these array pairs are only a slight relaxation of ARPAs, they are unlikely to admit such a similar inductive construction, which makes designing them all the more challenging.

\bibliographystyle{plain}
{\small\bibliography{sat-cd-bib}}

\appendix
In the following, for a positive integer $j$, $e^j$ denotes the $j$th canonical vector (whose dimension depends on the context).

\section{Supplementary material}\label{app-supp}

\subsection{$\mathcal{E}_q$ and $\mathcal{O}_q$ functions families}
\label{sec-func}
\label{sec-CSP_Eq+Oq}

{\bf{\em Function decomposition.}}
Let $q$ and $k$ be two positive integers. Analogous to the concept of {\em even} and {\em odd} functions, any function $P:(\Sigma_q, +)^k \rightarrow \mathbb{R}$ can be decomposed into the sum of a function of $\mathcal{E}_q$ and a function of $\mathcal{O}_q$. To do this, we consider the function $P_E$ defined on $\Sigma_q^k$ by:
\begin{align}\label{eq-P_E}
\begin{array}{rll}
	P_E(y) 	&:= 1/q\times\sum_{a =0}^{q -1} P_\mathbf{a}(y)\\
			& = \sum_{a =0}^{q -1} P(y_1 +a ,\ldots, y_k +a)/q,		&y_1 ,\ldots, y_k\in\Sigma_q
\end{array}
\end{align}

For example, the expression $\sum_{a =0}^{q -1} AllZeros^{k, q}(y_1 +a ,\ldots, y_k +a)$ evaluates to 1 {\em if and only if} there exists in $\set{0 ,\ldots, q -1}$ an integer $a$ for which 
$y_1 +a, y_2 +a ,\ldots, y_k +a$ are all zero modulo $q$, which happens {\em if and only if} $y_1 ,\ldots, y_k$ are all equal. Otherwise, it evaluates to 0. So for $P =AllZeros^{k, q}$, $P_E$ is $1/q\times AllEqual^{k, q}$. 

By construction, $P_E$ is stable under shifting all its inputs by an equal amount $a\in\Sigma_q$, while $P -P_E$ satisfies that $\sum_{a =0}^{q -1}(P -P_E)_{\mathbf{a}}$ is the constant function zero. Note that the function $P -P_E$ can actually be decomposed into the sum of the $q -1$ functions $(P -P_\mathbf{a})/q$, $a\in\Sigma_q$, all of which belong to $\mathcal{O}_q$ and have a mean value of zero. 
Definitions \cref{eq-E_q} of $\mathcal{E}_q$ and \cref{eq-O_q} of $\mathcal{O}_q$ state precisely that $P\in\mathcal{E}_q$ {\em if and only if} $P_E =P$, and $P\in\mathcal{O}_q$ {\em if and only if} $P_E$ is constant (in which case $P_E$ is exactly the constant function $r_P$). 

\bigskip
{\bf{\em Restrictions $\mathsf{CSP(\mathcal{O}_q)}$ and $\mathsf{CSP(\mathcal{E}_q)}$ of $\mathsf{CSP\!-\!q}$.}}
 %
$\mathsf{CSP(\mathcal{O}_q)}$ is remarkable in that it is trivially approximable within a differential factor of $1/q$ (see \cref{sec-E-O_q} and the introduction of \cref{sec-vois-prev}), but $\mathbf{NP-hard}$ to approximate within any constant factor greater than $1/q$, and this even for $\mathsf{E3\,CSP(\mathcal{O}_q)}$ \cite{H97}. 
 %
Regarding $\mathsf{CSP(\mathcal{E}_q)}$, we observe that a function $P$ on $\Sigma_q^k$, where $k$ is a positive integer, can be interpreted as a $(k +1)$-ary function of $\mathcal{E}_q$. Precisely, we associate with $P$ the function $P^E$ defined on $\Sigma_q^{k +1}$ by: 
\begin{align}\label{eq-P^E}
\begin{array}{rll}
	P^E(y_0, y_1 ,\ldots, y_k)
		&:=	P_{\mathbf{-y_0}}(y_1 ,\ldots, y_k)\\
		&=			P(y_1 -y_0 ,\ldots, y_k -y_0),			&y_0, y_1 ,\ldots, y_k\in\Sigma_q
\end{array}
\end{align}

For example, consider the function $AllZeros^{k, q}$. For $y_0, y_1 ,\ldots, y_k\in\Sigma_q$, 
$y_1 -y_0, y_2 -y_0 ,\ldots, y_k -y_0$ are all equal to zero modulo $q$ {\em if and only if}  
$y_1, y_2 ,\ldots, y_k$ are all equal to $y_0$. 
So for $P =AllZeros^{k, q}$, $P^E$ is $AllEqual^{k +1, q}$.  

From the transformation \cref{eq-P^E}, we derive a differential approximation-preserving reduction (see \cref{sec-def-reduc}) $(f, g)$ from $\mathsf{k\,CSP\!-\!q}$ to $\mathsf{(k +1)CSP(\mathcal{E}_q)}$, which induces no loss on the approximation guarantee. 
Given an instance $I$ of $\mathsf{k\,CSP\!-\!q}$, the algorithm $f$ introduces an auxiliary variable $z_0$ and replaces each constraint $P_i(x_{i_1} ,\ldots, x_{i_{k_i}})$ of the input instance with the new constraint $P_i(x_{i_1} -z_0 ,\ldots, x_{i_{k_i}} -z_0)$.  
The function $g(I, .)$ then maps each solution  
	$(x, z_0) =(x_1 ,\ldots, x_n, z_0)$ of $f(I)$ to the solution 
	$x -\mathbf{z_0} =(x_1 -z_0 ,\ldots, x_n -z_0)$ of $I$. 
By definition of $f(I)$, this solution realizes the same objective value on $I$ as $(x, z_0)$ on $f(I)$. 

Considering that $\mathsf{(k +1)CSP(\mathcal{E}_q)}$ is a special case of $\mathsf{(k +1)CSP\!-\!q}$, this reduction somehow indicates that $\mathsf{(k +1)CSP(\mathcal{E}_q)}$ can be seen as an intermediate problem between $\mathsf{k\,CSP\!-\!q}$ and $\mathsf{(k +1)CSP\!-\!q}$. 
\Cref{pty-2CSP-B^1} of \cref{sec-B1} analyzes this reduction, when applied to $\mathsf{2\,CSP\!-\!2}$, in terms of differential approximation guarantees related to the $\tilde{B}^1$ neighborhood function.
 %
Examples of this reduction can be found in \cite{MPT03A,EG04}, where the constraints are either  disjunctions on Boolean variables, or their generalization to $q$-ary alphabets. More precisely, $\mathsf{NAE\,Sat\!-\!q}$ and $\mathsf{Sat\!-\!q}$ are the $q$-ary CSPs in which a constraint requires that a set of literals are not all equal for the former problem, are not all zero for the latter problem, where a literal $\ell_j$ is either the variable $x_j$ or its shift $x_j +a$ by some  $a\in[q -1]$ (see, {\em e.g.}, \cite{AM08}). Then for all positive integers $k$, $\mathsf{k\,Sat\!-\!q}$ $D$-reduces to $\mathsf{(k +1)NAE\,Sat\!-\!q}$ with no loss on the approximation guarantee \cite{MPT03A,EG04}. 
We denote symmetrically by $\mathsf{AE\,Sat\!-\!q}$ and $\mathsf{CCSP\!-\!q}$ the $q$-ary CSPs in which a constraint requires a set of literals to be all equal or all zero, respectively. Because of the connection between the two functions $AllZeros^{h, q}$ and $AllEqual^{h +1, q}$ for each $h\in[k]$, $\mathsf{k\,CCSP\!-\!q}$ also $D$-reduces to $\mathsf{(k +1)AE\,Sat\!-\!q}$ with no loss on the approximation guarantee.

\subsection{Approximability bounds of the literature} 
\label{sec-ap-bib}

In \cref{sec-dapx}, we claim that the 6-gadget of \cite{H97}, which reduces  $\mathsf{E3\,Lin\!-\!2}$ to $\mathsf{E2\,Lin\!-\!2}$, implies a differential approximability upper bound of 7/8 for the restriction of $\mathsf{E2\,Lin\!-\!2}$ to bipartite instances, denoted by $\mathsf{Bipartite\,E2\,Lin\!-\!2}$. In addition, \cref{tab-E_cmp} reports approximability bounds for the restriction of $\mathsf{E2\,CSP(\mathcal{I}_2^1)}$ to bipartite instances, as well as for $\mathsf{E3\,CSP(\mathcal{I}_2^2)}$. We show that these bounds are correct.

\medskip
First, for $\mathsf{Bipartite\,E2\,Lin\!-\!2}$, approximating the optimal gain over a random assignment or approximating the optimal gain over a worst solution somehow reduces to the same thing:
\begin{property}\label{pty-BipE2Lin2}
A solution of a bipartite instance of $\mathsf{E2\,Lin\!-\!2}$ is $g$-gain approximate {\em if and only if} it is $(1/2 +g/2)$-differential approximate.
\end{property}%

\begin{proof}
Let $I$ be an instance of $\mathsf{Bipartite\,E2\,Lin\!-\!2}$, and $(L, R)$ be a 2-coloring of $I$. 
For any two solutions $x$ and $y$ such that $y_L =x_L$ and $y_R =\bar x_R$, we have $v(I, x) +v(I, y) =\sum_{i =1}^m w_i$. This implies that 
	$\mathrm{opt}(I) +\mathrm{wor}(I) =\sum_{i =1}^m w_i =2\times\mathbb{E}_X[v(I, X)]$. Equivalently:
\begin{align}\label{eq-BipE2Lin2}
\begin{array}{rll}
	\mathbb{E}_X[v(I, X)] -\mathrm{wor}(I)
	&=\mathrm{opt}(I) -\mathbb{E}_X[v(I, X)]
	&=\displaystyle\left(\mathrm{opt}(I) -\mathrm{wor}(I)\right)/2
\end{array}
\end{align}

We deduce:
$$\begin{array}{rl}
\displaystyle
	\frac{v(I, x) -\mathrm{wor}(I)}{\mathrm{opt}(I) -\mathrm{wor}(I)}
&=\displaystyle
	\frac{v(I, x) -\mathbb{E}_X[v(I, X)]}{2\left(\mathrm{opt}(I) -\mathbb{E}_X[v(I, X)]\right)}
	+\frac{\mathbb{E}_X[v(I, X)] -\mathrm{wor}(I)}{2\left(\mathbb{E}_X[v(I, X)] -\mathrm{wor}(I)\right)}
\end{array}$$

The result is straightforward.
\end{proof}

In \cite{AN04}, Alon and Naor show that $\mathsf{Bipartite\,E2\,Lin\!-\!2}$ is approximable within gain approximation ratio $2\ln(1 +\sqrt{2})/\pi$. According to \cref{pty-BipE2Lin2}, equivalently, they show that $\mathsf{Bipartite\,E2\,Lin\!-\!2}$ is approximable within differential approximation ratio $1/2 +\ln(1 +\sqrt{2})/\pi$. 

\medskip
Second, on the Boolean alphabet, $\mathsf{Ek\,CSP(\mathcal{I}_2^{k -1})}$ is nothing more than $\mathsf{Ek\,Lin\!-\!2}$. This is because, for any positive integer $k$, the $k$-ary balanced $(k -1)$-wise independent Boolean functions are functions of the form $a XNOR^k +b$ for two constants $a$ and $b$:

\begin{property}\label{pty-I_2^(k -1)-XOR}
Let $k$ be a positive integer. Then a pseudo-Boolean function $P:\set{0, 1}^k\rightarrow\mathbb{R}$ is balanced $(k -1)$-wise independent {\em if and only if} it coincides with $XNOR^k$ up to an affine transformation. 
\end{property}

\begin{proof}
Consider two Boolean vectors $u, v\in\set{0, 1}^k$ with the same number of non-zero coordinates. We denote by $J =\set{j_1 ,\ldots, j_\kappa}$ the set of coordinate indices where $u$ and $v$ differ. Thus, $v$ can be described as the vector $u +\sum_{r =1}^\kappa e^{j_r}$ where, by the assumption $XNOR^k(u) =XNOR^k(v)$, $\kappa$ is even. Therefore, we can write $P(v) -P(u)$ as: 
$$\begin{array}{rll}
P(v) -P(u)	&=P\left(u +\sum_{r =1}^\kappa e^{j_r}\right)	-P(u)\\
			&=\sum_{s =1}^{\kappa/2} \left(
				P(u +\sum_{r =1}^{2s} e^{j_r}) -P(u +\sum_{r =1}^{2s -2} e^{j_r}) 
			\right)
\end{array}$$

Let $s\in[\kappa/2]$. If $P\in\mathcal{I}_k^{k -1}$, then we can deduce from \cref{eq-I_q^t} that we have:
$$\begin{array}{rll}
P(u +\sum_{r =1}^{2s} e^{j_r})
	&=P(u +\sum_{r =1}^{2s -2} e^{j_r})
	&=2r_P -P(u +\sum_{r =1}^{2s -1} e^{j_r})	
\end{array}$$

We conclude that $P$ evaluates to the same value on all vectors with an even number of non-zero coordinates on the one hand, and on all vectors with an odd number of non-zero coordinates on the other hand. In other words, there exist two real numbers $a$ and $b$ such that $P$ is the function $a \times XNOR^k +b \times XOR^k$ or, equivalently, $(a -b)XNOR^k +b$. 
\end{proof}

According to \cref{pty-I_2^(k -1)-XOR}, the approximability bound established in \cite{AN04} for bipartite instances of $\mathsf{E2\,Lin\!-\!2}$ actually holds for bipartite instances of $\mathsf{E2\,CSP(\mathcal{I}_2^1)}$. 
Moreover, Khot and Naor in \cite{KN08} show that on instances of $\mathsf{E3\,Lin\!-\!2}$, the optimal gain over a random assignment is approximable within an expected factor of $\Omega(\sqrt{\ln n/n})$. Similarly, this fact actually holds for $\mathsf{E3\,CSP(\mathcal{I}_2^2)}$.

\begin{figure}[t]{\footnotesize\begin{center} 
\includegraphics[trim=0pt 0pt 0pt 0pt]{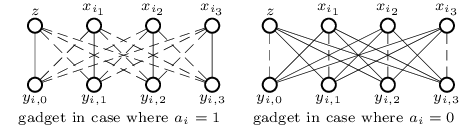}
\caption{\label{fig-gadget}The $6$-gadget from \cite{H97} that transforms each constraint 
$(x_{i_1} +x_{i_2} +x_{i_3}\equiv a_i\bmod{2})$ of an instance $I$ of $\mathsf{Max\,E3\,Lin\!-\!2}$ into a set of $XNOR^2$ (shown as solid lines) and $XOR^2$ (shown as dashed lines) constraints.}
\end{center}}\end{figure} 

\medskip
Finally, the 6-gadget of \cite{H97}, which  reduces $\mathsf{E3\,Lin\!-\!2}$ to $\mathsf{E2\,Lin\!-\!2}$, implies a gain approximability upper bound of 3/4 for bipartite instances of $\mathsf{E2\,Lin\!-\!2}$:

\begin{proposition}\label{prop-BipE2Lin2}	
$\mathsf{Bipartite\,E2\,Lin\!-\!2}$ is inapproximable within any constant gain approximation ratio better than $3/4$, unless $\mathbf{P} =\mathbf{NP}$.
\end{proposition}%

\begin{proof}
Consider an instance $I$ of $\mathsf{Max\,E3\,Lin\!-\!2}$.
The reduction of \cite{H97} first introduces $4m +1$ binary auxiliary variables $y_{i,0}, y_{i,1}, y_{i,2}, y_{i,3}, i\in[m]$, and $z$. Then, for each constraint 
	$(x_{i_1} +x_{i_2} +x_{i_3}\equiv a_i\bmod{2})$ 
of $I$, it generates sixteen $XOR^2$ or $XNOR^2$ constraints, all of weight $w_i/2$. \Cref{fig-gadget} illustrates these constraints. 

To a solution $(x, y, z)$ of the resulting instance of $\mathsf{Max\,E2\,Lin\!-\!2}$, the reduction associates the solution $x$ if $z =0$, and $\bar x$ otherwise, of $I$.  
We denote by $I'$ the instance produced, and by $w(I)$ and $w(I')$ the sum of the constraint weights on $I$ and $I'$, respectively. $I'$ is obviously bipartite. It also satisfies (see \cite{H97}):
\begin{align}
\label{eq-gadget-w}
	w(I')				&=8\times w(I)												\\[-2pt]
\label{eq-gadget-apx0}
	v(I, x)				&\geq v(I', (x, y, 0)) -5w(I),	&(x, y)\in\set{0, 1}^{n+4m}	\\[-2pt]
\label{eq-gadget-apx1}
	v(I, \bar x)		&\geq v(I', (x, y, 1)) -5w(I),	&(x, y)\in\set{0, 1}^{n+4m}	\\[-2pt]
\label{eq-gadget-opt}
	\mathrm{opt}(I')	&=\mathrm{opt}(I) +5 w(I)
\end{align}

Suppose we can compute a solution $(x, y, z)$ on $I'$ that is $\varepsilon$-gain approximate, for some $\varepsilon >0$. Since the two solutions $(x, y, z)$ and $(\bar x, \bar y, \bar z)$ both perform the same objective value on $I'$, we can assume without loss of generality that $z =0$. So consider the solution $x$ of $I$. We observe successively:
\begin{align}
v(I, x)	&\geq 	v(I', (x, y, z)) -5w(I)	\qquad\hfill\text{by \cref{eq-gadget-apx0}}	\nonumber\\
		&\geq 	\varepsilon\,\mathrm{opt}(I') +(1 -\varepsilon)\times w(I')/2	-5w(I)	
							\qquad\hfill\text{by assumption on $(x, y, z)$}			\nonumber\\
		&=		\varepsilon\left(\mathrm{opt}(I) +5 w(I)\right) +(1 -\varepsilon) 4w(I) -5 w(I)	
							\qquad\hfill\text{by \cref{eq-gadget-opt,eq-gadget-w}}	\nonumber\\
		&=		\varepsilon\,\mathrm{opt}(I) -(1 -\varepsilon) w(I)			\label{eq-gadget}
\end{align}

Now, for all constants $\delta >0$, $\mathsf{Gap_{(1 -\delta, 1/2 +\delta)}E3\,Lin\!-\!2}$ is $\mathbf{NP}$-hard \cite{H97}. This means that, given an instance $I$ of $\mathsf{Max\,E3\,Lin\!-\!2}$ that verifies either $\mathrm{opt}(I)\geq (1 -\delta) w(I)$ or $\mathrm{opt}(I)\leq(1/2 +\delta) w(I)$, it is $\mathbf{NP}$-hard to decide which of these two cases occurs.
Let $\delta >0$, and consider such an instance $I$. 
If $\mathrm{opt}(I)\geq (1 -\delta) w(I)$, then by \cref{eq-gadget}, we have:
$$\begin{array}{rll}
	v(I, x)	&\geq 	\varepsilon\times(1 -\delta) w(I) -(1 -\varepsilon) w(I)
			&=		w(I)\times\left( (2 -\delta)\varepsilon -1 \right)
\end{array}$$

Otherwise $\mathrm{opt}(I)\leq(1/2 +\delta) w(I)$, and hence $v(I, x)\leq(1/2 +\delta) w(I)$. 
Note that $(2 -\delta)\varepsilon -1> 1/2 +\delta$ {\em iff} $\delta <(2\varepsilon -3/2)/(1 +\varepsilon)$, while $(2\varepsilon -3/2)/(1 +\varepsilon) >0$ {\em iff} $\varepsilon >3/4$. 
Thus, if $\varepsilon >3/4$, then for sufficiently small $\delta$, we can decide whether $\mathrm{opt}(I)\geq (1 -\delta) w(I)$ or $\mathrm{opt}(I)\leq(1/2 +\delta) w(I)$ by comparing $v(I, x)$ with $(1/2 +\delta) w(I)$, implying $\mathbf{P} =\mathbf{NP}$.
\end{proof}

According to \cref{pty-BipE2Lin2}, \cref{prop-BipE2Lin2} equivalently states that $\mathsf{Bipartite\,E2\,Lin\!-\!2}$ is inapproximable within any constant differential factor greater than 7/8, unless $\mathbf{P} =\mathbf{NP}$.

\subsection{Computation of optimal designs}
\label{sec-cd}

We explain how we calculated the arrays and the values shown in 
\cref{tab-OA-DS-ex-q=3,tab-rho,tab-Gamma-ex,tab-Gamma_E-ex,tab-gamma,tab-gamma_E}. 
\begin{table}\footnotesize{\begin{center}
$\begin{array}{rl|rl}
\rho(\nu, q, t)
	&=\left\{\begin{array}{rrl}
		\multicolumn{3}{l}{\max_{P:\Sigma_q^\nu\rightarrow[0, 1], R} P(\mathbf{0})}		\\
		s.t.	&\multicolumn{2}{c}{\cref{eq-PL-defR,eq-PL-rho-0max,eq-PL-rho-Iqt}}		\\
				&R				&=		1												\\
	\end{array}\right.
&F(\nu, q, t)
	&=\left\{\begin{array}{rrl}
		\multicolumn{3}{l}{\min_{P:\Sigma_q^\nu\rightarrow\mathbb{N}, R} 	R}			\\
		s.t.	&\multicolumn{2}{c}{\cref{eq-PL-defR,eq-PL-rho-0max,eq-PL-rho-Iqt}}		\\
				&R 				&\geq	1												\\
	\end{array}\right.\\
R(\nu, q, t)
	&=\left\{\begin{array}{rrl}
		\multicolumn{3}{l}{\min_{P:\Sigma_q^\nu\rightarrow\mathbb{N}, R} 	R}			\\
		s.t.	&\multicolumn{2}{c}{\cref{eq-PL-defR,eq-PL-rho-0max,eq-PL-rho-Iqt}}		\\
				&P(\mathbf{0})	&\geq	\rho(\nu, q, t)\times R							\\
				&R 				&\geq	1												\\
\end{array}\right.
&R^*(\nu, q, t)
	&=\left\{\begin{array}{rrl}
		\multicolumn{3}{l}{\max_{P:\Sigma_q^\nu\rightarrow\mathbb{N}, R} P(\mathbf{0})}	\\
		s.t.	&\multicolumn{2}{c}{\cref{eq-PL-defR,eq-PL-rho-0max,eq-PL-rho-Iqt}}		\\
				&R 				&=		F(\nu, q, t)									\\ 
	\end{array}\right.
\end{array}$
\caption{\label{tab-PL-rho}Linear programs for orthogonal arrays.}
\end{center}}\end{table}

\medskip
{\bf{\em Orthogonal arrays and difference schemes of \cref{sec-average}.}}
Let $q\geq 1$, $t\geq 1$, and $\nu\geq t$ be three integers. To model orthogonal arrays of strength $t$ with $\nu$ columns and entries from the symbol set $\Sigma_q$, we associate with each $u\in\Sigma_q^\nu$ a variable $P(u)$ that represents either the number of occurrences or the frequency of $u$ in the array, depending on whether we are modeling the array itself or the measure it induces on $\Sigma_q^\nu$. Thus, these variables have domain $\mathbb{N}$ or $[0, 1]$, depending on the context. We use an additional (continuous) variable $R$ to represent either the number of rows in the array (in which case $R$ must be greater than or equal to 1), or the total frequency of the words of $\Sigma_q^\nu$ in the array (in which case $R$ must be 1). To avoid symmetries, we consider only arrays in which the row of 0's has the maximum frequency. 

The variables $P(u)$, $u\in\Sigma_q^\nu$ and $R$ are always subject to the following constraints:
\begin{align}
\textstyle\label{eq-PL-defR}
	\sum\limits_{u\in\Sigma_q^\nu} P(u)						&=		R		&					\\
\textstyle\label{eq-PL-rho-0max}
	P(\mathbf{0})											&\geq	P(u),	&u\in\Sigma_q^\nu	\\		\textstyle\label{eq-PL-rho-Iqt}
	\sum\limits_{u =(u_1 ,\ldots, u_\nu)\in\Sigma_q^\nu:u_J =v} P(u) 	&= 	\frac{1}{q^t}R,	
		&J =(j_1 ,\ldots, j_t)\in[\nu]^t,\,j_1 <\ldots<j_t,\ v\in\Sigma_q^t		
\end{align}

Constraint \cref{eq-PL-defR} defines $R$ as the total frequency or number of occurrences of the words from $\Sigma_q^\nu$ in the array. Inequalities \cref{eq-PL-rho-0max} force the word of all zeros to have the highest frequency. Finally, relations \cref{eq-PL-rho-Iqt} ensure that the corresponding array induces a balanced $t$-wise independent distribution over $\Sigma_q^\nu$. 

When the variables are integer, depending on the optimization objective, we consider the additional constraint $R \geq 1$ to prohibit the trivial solution $R =0 =P(u), u\in\Sigma_q^\nu$.
We consider two optimization criteria: the number of rows (which we want to minimize), and the maximum frequency of a word (which we want to maximize). Specifically, we are interested in:
\begin{enumerate}
	\item computing $\rho(\nu, q, t)$, which can be done by solving the upper-left linear program in continuous variables of \cref{tab-PL-rho};
	\item minimizing the number of rows in an OA that realizes $\rho(\nu, q, t)$, which can be done by solving the lower-left program of \cref{tab-PL-rho};
	\item computing $F(\nu, q, t)$, which can be done by solving the upper-right program of \cref{tab-PL-rho};
	\item maximizing the maximum frequency of a word in an $OA(F(\nu, q, t), \nu, q, t)$, which can be done by solving the lower-right program of \cref{tab-PL-rho}.
\end{enumerate}

\smallskip
The case of difference schemes is rather similar.
First, in order to avoid symmetries, we associate a variable $P(u)$ only with the words $u\in\Sigma_q^\nu$ with a zero first coordinate. 
Second, instead of relations \cref{eq-PL-rho-Iqt}, for each 
$J =(j_1 ,\ldots, j_t)\in[\nu]^t$ with $j_1 <\ldots<j_t$ 
and each $v\in\set{0}\times\Sigma_q^{t -1}$, we consider the constraint: 
\begin{align}\textstyle\label{eq-PL-rho_E-Iqt}
\sum_{a =0}^{q -1}\sum_{u\in\set{0}\times\Sigma_q^{\nu -1}:u_J =v +\mathbf{a}} P(u)	&= R/q^{t -1}
\end{align}

\medskip
{\bf{\em Array pairs of \cref{sec-reduc}.}}
Let $k \geq 2$, $p \geq k$, and $q >p$ be three integers. We denote by $\mathcal{U}$ the set of words $u\in\Sigma_q^q$ with at most $p$ pairwise distinct coordinates. To compute the number $\gamma(q, p, k)$ or an ARPA that realizes it, we consider the variables $P(u), u\in\mathcal{U}$, $Q(u), u\in\Sigma_q^q$ and $R$, so as to model the array $\Psi$ (or the frequencies in $\Psi$), the array $\Phi$ (or the frequencies in $\Phi$), and the number of rows in these arrays (or the total frequency of words of $\Sigma_q^q$ in these arrays, in which case $R$ is equal to 1).
These variables are subject to the constraints that $R$ must coincide with $\sum_{u\in\mathcal{U}} P(u)$, and the constraints
\begin{align}\label{eq-PL-gamma-Iqt}
\begin{array}{rl}
\sum_{u =(u_0 ,\ldots, u_{q -1})\in\mathcal{U}:u_J =v} P(u) &=\textstyle
\sum_{u =(u_0 ,\ldots, u_{q -1})\in\Sigma_q^q:u_J =v} Q(u),\qquad\qquad\\
\multicolumn{2}{r}{J =(j_1 ,\ldots, j_k)\in\Sigma_q^k,\,j_1 <\ldots<j_k,\ v\in\Sigma_q^k}
\end{array}
\end{align}
which ensure that the corresponding pair $(\Psi, \Phi)$ of arrays satisfies $\mu_\Psi -\mu_\Phi\in\mathcal{I}_q^k$.

\smallskip
The case of relaxed ARPAs is quite similar. First, we eliminate symmetric solutions by restricting the variables $P(u)$ and $Q(u)$ to words $u$ of $\Sigma_q^q$ such that $u_0 =0$. Second, instead of the constraints \cref{eq-PL-gamma-Iqt}, for each $J =(j_1 ,\ldots, j_k)\in\Sigma_q^k$ with $j_1 <\ldots<j_k$ and each $v\in\set{0}\times\Sigma_q^{k -1}$, we consider the constraint: 
\begin{align}\label{eq-PL-gamma_E-Iqt}
\textstyle	
	\sum_{a =0}^{q -1}\sum_{u\in\mathcal{U}: u_0 =0 \wedge  u_J =v +\mathbf{a}} P(u) 
&=\textstyle	
	\sum_{a =0}^{q -1}\sum_{u\in\Sigma_q^q: u_0 =0 \wedge u_J =v +\mathbf{a}} Q(u)
\end{align}

For both problems, the goal is to maximize the ratio $Q(0, 1 ,\ldots, q -1)/R$. 

\section{Extended proofs}\label{app-pv}

\subsection{Lower bounds given in \cref{sec-E-dapx} for the average differential ratio}
\label{sec-ap-avd}

In the following, we denote by $\mathrm{avd}(I)$ the average differential ratio on the considered $k$-CSP instance $I$.

\medskip
{\bf{\em $\bullet$ \Cref{cor-E-k_partite}.}}
For \cref{it-B-nu}: $\mathrm{avd}(I)\geq\rho(\nu, q, k)$ by \cref{thm-E-OA}, 
	while $\rho(\nu, q, k)\geq 1/F(\nu, q, k)$, which is $1/q^k$ by \cref{thm-Bush}. 

For \cref{it-B-qpp}: $\mathrm{avd}(I)\geq\rho(\nu, p^\kappa, k)$ by \cref{thm-E-OA,thm-reduc} and thus $\mathrm{avd}(I)\geq 1/F(\nu, p^\kappa, k)$, while $F(\nu, p^\kappa, k) =(p^\kappa)^k$ by \cref{thm-Bush}. 
Furthermore, we have the inequalities $p^\kappa\leq 2^{\lceil\log_2 q\rceil}$ (by definition of $p^\kappa$) and $2^{\lceil\log_2 q\rceil}\leq 2(q -1)$ (by definition of $\lceil\log_2 q\rceil$). So the inequality $p^\kappa\leq 2(q -1)$ and thus $1/p^{\kappa k}\geq 1/(2q -2)^k$ holds.

The argument for \cref{it-B-qp2} is similar.

\medskip
{\bf{\em $\bullet$ \Cref{cor-E-2CSPs,cor-E-3CSPs}.}}
By \cref{thm-E-OA}, we have $\mathrm{avd}(I)\geq\rho(\nu, q, 2)$ and $\mathrm{avd}(I)\geq\rho(\nu, 2, 3)$ for $\mathsf{2\,CSP\!-\!q}$ and $\mathsf{3\,CSP\!-\!2}$, respectively. 

For \cref{cor-E-2CSPs}: let $s =\lceil(\nu -1)/q\rceil q +1$. Note that $s$ verifies $s\in\set{\nu ,\ldots, \nu +q -1}$ and, if $\nu\geq 2$, $s >q$. Since $s\geq\nu$, we have $\rho(\nu, q, 2)\geq\rho(s, q, 2)$. Since $\nu\geq 2$, then $s\geq q$, and thus we know from \cref{thm-E-2CSPs} that the value of $\rho(s, q, 2)$ is $1/(s(q -1) +1)$. 

For \cref{cor-E-3CSPs}: we know from \cref{cor-E-3CSP2} that for $\nu\geq3$, $1/\rho(\nu, 2, 3)$ is the quantity $4\lceil\nu/2\rceil$.

\medskip
{\bf{\em $\bullet$ \Cref{cor-E-kCSP-2}.}}
First, assume that $k$ is even, and let $s =2^{\lceil\log_2(\nu +1)\rceil} -1$. 
Since $s\geq\nu$, we deduce from \cref{thm-E-OA} that $\mathrm{avd}(I)$ is at least $\rho(s, 2, k)$ and hence at least $1/F(s, 2, k)$. 
Furthermore, by the assumptions $\nu\geq k$ and $k\geq 4$ we have $\lceil\log_2(\nu +1)\rceil\geq\lceil\log_2 5\rceil =3$, and hence $s\geq 7$. 
Moreover, since $s\geq\nu$ and $\nu\geq k$ while $s$ is odd and $k$ is even, then $s$ is greater that $k$. The integers $k/2$ and $s$ therefore satisfy the conditions of \cref{thm-BCH}, implying that $F(s, 2, k)$ is at most $(s +1)^{k/2}$. Finally, $s$ by construction satisfies $s +1\leq 2\nu$.

For odd values of $k$, we deduce from \cref{thm-E-OA} and \cref{prop-OA-E-F} that $\mathrm{avd}(I)$ is at least $1/2\times\rho(\nu -1, 2, k -1)$. Since $\nu -1\geq k -1\geq 4$, following the same argument as for the even case, we get the inequalities: 
$$\begin{array}{rll}
	1/\rho(\nu -1, 2, k -1)	&\leq F(2^{\lceil\log_2 \nu\rceil}-1, 2, k -1)
							&\leq 2^{\lceil\log_2\nu\rceil\times(k -1)/2}
\end{array}$$

The conclusion follows.

\medskip
{\bf{\em $\bullet$ \Cref{cor-E-kCSP-q}.}}
Given that $\nu\geq k$, we know from \cref{thm-E-OA,thm-reduc} that $\rho(\nu, p^\kappa, k)$, and hence $1/F(\nu, p^\kappa, k)$, is a valid lower bound for $\mathrm{avd}(I)$. 
Now for $s =(p^\kappa)^{\lceil\log_{p^\kappa}\nu\rceil}$ we have $F(\nu, p^\kappa, k)\leq F(s, p^\kappa, k)$, while by \cref{thm-Bierb} $F(s, p^\kappa, k)$ is at most $p^\kappa\times s^{k -\lceil\log_{p^\kappa}(k)\rceil}$.

\medskip
{\bf{\em $\bullet$ \Cref{cor-E-kCSP(I_q^t)-Bush}.}}
The argument is basically the same as for the general case, except that we consider the lower bound of $\rho(\nu -t, q, \min\set{\nu -t, k})$ --- instead of $\rho(\nu, q, \min\set{\nu, k})$ --- for $\mathrm{avd}(I)$.

\medskip
{\bf{\em $\bullet$ \Cref{cor-E-kCSP(E_q)}.}}
According to \cref{thm-E-OA}, we can consider the lower bound $\rho_E(\nu, q, k)$ and hence $1/E(\nu, q, k)$ for $\mathrm{avd}(I)$.

First assume that $q$ is 2 and $k$ is an odd integer at least 3.
If $k =3$, then by \cref{cor-E-3CSP2} we have $\rho_E(\nu, 2, 3) =1/(2\lceil\nu/2\rceil)$. 
If $\nu\in\{k, k +1\}$, we successively deduce from \cref{prop-OA-E-F} and \cref{thm-Bush} that $\rho_E(\nu, 2, k)$ is equal to $2\times\rho(\nu, 2, k)$ and thus to $2\times1/2^k$.
If $k\geq 5$, we know from \cref{prop-OA-E-F} that $E(\nu, 2, k)$ coincides with $F(\nu -1, 2, k -1)$, where $\nu -1\geq k -1\geq 4$. So the argument for this case is essentially the same as for $\mathsf{k\,CSP\!-\!2}$. 

Now assume $q\geq 3$ and let $s =q^{\lceil\log_q\nu\rceil}$. For this case, we consider the inequality $E(\nu, q, k)\leq E(s, q, k)$ and the upper bound \cref{thm-Bierb} provides for $E(s, q, k)$.

\subsection{Array pairs of \cref{sec-reduc}}
\label{sec-qkk}
\label{sec-T-UB}
\label{sec-gEq22}

\subsubsection{\Cref{lem-gamma-p=k}}

Consider \cref{algo-Gamma_rec}. Our goal is to prove that, at the end of the algorithm, the difference $\mu_\Psi -\mu_\Phi$ of the frequencies of rows occurring in $\Psi$ and $\Phi$ is balanced $k$-wise independent. To this end, we first establish a technical lemma. 

\begin{lemma}
For three natural numbers $a$, $b$, and $c\leq b$, we define:
\begin{align}\label{eq-id-combi-S-def}
S(a, b, c)	&\textstyle:=	
				\sum_{r\geq 0}	(-1)^r \binom{a}{r} \binom{b -r}{c -r}
\end{align}

These numbers satisfy:
\begin{align}\label{eq-S}
S(a, b, c)	&=\textstyle\binom{b -a}{c},	&a, b, c\in\mathbb{N},\ b\geq a,\,c\leq b
\end{align}
\end{lemma}

\begin{proof}
Let $a\in\mathbb{N}$. We show by induction on $b$ that \cref{eq-S} is satisfied at $(a, b, c)$ for all pairs $(b, c)$ of natural numbers such that $b\geq a$ and $c\leq b$.
 %
If $b =a$, then for $c\in\set{0 ,\ldots, a}$, given that $\binom{a}{r}\binom{a -r}{c -r} =\binom{a}{c}\binom{c}{r}$, $r\in\set{0 ,\ldots, c}$, we have:
$$\begin{array}{rl}
S(a, a, c)	&=\binom{a}{c}\times\sum_{r =0}^c(-1)^r\binom{c}{r}	
\end{array}$$

So $S(a, a, c)$ is 1 if $c =0$ and 0 otherwise, as well as $\binom{0}{c}$. 
Identity \cref{eq-S} is therefore satisfied at $(a, a, c)$ for all $c\in\set{0 ,\ldots, a}$. 
Now assume that it is satisfied at $(a, b -1, c)$ for every $c\in\set{0 ,\ldots, b -1}$, where $b$ is some integer strictly greater than $a$, and consider an integer $c\in\set{0 ,\ldots, b}$. We want to show that \cref{eq-S} is satisfied at $(a, b, c)$. 
If $c =0$, then $S(a, b, 0) =(-1)^0\binom{a}{0}\binom{b}{0} =1$.
If $c =b$, then $S(a, b, b) =\sum_{r =0}^a(-1)^r\binom{a}{r}$.
In both cases, $S(a, b, c)$ actually coincides with $\binom{b -a}{c}$. 
Now suppose $c\in[b -1]$. In this case, we deduce successively:
$$\begin{array}{rlll}
S(a, b, c)	
	&=\sum_{r =0}^{\min\set{a, c}}(-1)^r\binom{a}{r}
			\left(\binom{b -1 -r}{c -r} +\binom{b -1 -r}{c -1 -r}\right)
													&\text{by Pascal's rule}			\\
	&=S(a, b -1, c) +S(a, b -1, c -1)	&\text{according to \cref{eq-id-combi-S-def}}	\\
	&=\binom{b -a -1}{c} +\binom{b -a -1}{c -1}		&\text{by induction hypothesis}		\\
\end{array}$$

So $S(a, b, c) =\binom{b -a}{c}$, which completes the argument.
\end{proof}

We now prove that, at the end of \cref{algo-Gamma_rec}, the arrays $\Psi$ and $\Phi$ verify:
$$\begin{array}{rl}
\begin{array}{rl}
\card{\set{r\in[R']\,|\,\Psi_r^J =v}} 
	-\card{\set{r\in[R']\,|\,\Phi_r^J =v}}	&=0,
\ \ \ \ \\\multicolumn{2}{r}{
	J =(j_1 ,\ldots, j_k)\in\Sigma_q^k,\ j_1 <\ldots< j_k,\ v\in\Sigma_q^k}
\end{array}	&\cref{eq-qkk}
\end{array}$$
where we recall that $R'$ refers to the final number of rows in $\Psi$ and $\Phi$.

\begin{proof}
Consider a subsequence $J$ of length $k$ of $(0 ,\ldots, q -2)$. Since $(\Psi,\Phi)$ initially is an ARPA of strength $k$, $(\Psi^J_r\,|\,r \in[R])$ and $(\Phi^J_r\,|\,r \in[R])$ define the same multiset of rows. The same holds for the two subarrays $(\Psi^J_r\,|\,r\in\set{R +1 ,\ldots, R'})$ and $(\Phi^J_r\,|\,r\in\set{R +1 ,\ldots, R'})$, due to the shape of the rows inserted by the algorithm. Thus, it remains to show that for any subsequence $J =(j_1 ,\ldots, j_{k -1})$ of length $k -1$ of  $(0 ,\ldots, q -2)$ and any $v\in\Sigma_q^k$, the two subarrays $(\Psi^J,\Psi^{q -1})$ and $(\Phi^J,\Phi^{q -1})$ coincide with $v$ on the same number of rows. 
We consider three cases:

$\bullet$ $v\notin\set{j_1, q -1}\times\ldots\times\set{j_{k -1}, q -1}\times\set{0, q -1}$: 
by construction, given $M\in\set{\Psi, \Phi}$, $(M_r^J, M_r^{q -1}) =v$ cannot occur unless $r\leq R$ and $(M_r^J, M_r^{q -1}) =(M_r^J, M_r^0)$. Subarrays $(\Psi^J, \Psi^{q -1})$ and $(\Phi^J, \Phi^{q -1})$ therefore coincide with $v$ on the same number of rows, due to the initial assumption on $(\Psi, \Phi)$.

$\bullet$ $(v_1 ,\ldots, v_{k -1}) =J$ and $v_k\in\set{0, q -1}$. 
If $v_k =q -1$, then the $R^*$ rows of the form $(0, 1 ,\ldots, q -1)$ in $\Phi$ and the $R^*$ rows of the form $(\alpha(J), q -1)$ in $\Psi$ are the only rows of $\Psi$ and $\Phi$ that coincide with $v$ on their coordinates of index $(J, q -1)$. 

Now suppose that $v_k =0$. 
In $\Psi$, the rows $\Psi_r$ satisfying $(\Psi^J_r, \Psi^{q -1}_r) =v$ are exactly those rows initially satisfying $(\Psi^J_r, \Psi^0_r) =v$. 
Since $\mu_\Psi -\mu_\Phi$ is initially balanced $k$-wise independent, the number of such rows in $\Phi$ is the same as in $\Psi$. 
Now, the rows $\Phi_r$ satisfying $(\Phi^J_r, \Phi^{q -1}_r) =v$ are all but $R^*$ of the rows initially satisfying $(\Phi_r^J, \Phi_r^0) =v$, and the $R^*$ rows of the form $(\alpha(J), 0)$. 

$\bullet$ $v\in\set{j_1, q -1}\times\ldots\times\set{j_{k -1}, q -1}\times\set{0, q -1}$ and $(v_1 ,\ldots, v_{k -1}) \neq J$. Since $(v_1 ,\ldots, v_{k -1})$ has at least one coordinate equal to $q -1$, given $M\in\set{\Psi, \Phi}$, $(M_r^J, M_r^{q -1}) =v$ cannot occur unless $r >R$. 
So we count in both arrays the number of rows of the form $(\alpha(H), v_k)$ satisfying $\alpha(H)_J =(v_1 ,\ldots, v_{k -1})$, where $H$ is a subset of size at most $k -1$ of $\Sigma_{q -1}$. 
Let $L$ denote the set of indices $j_s\in J$ for which $v_s =j_s$. For a subset $H$ of $\Sigma_{q -1}$, we have $\alpha(H)_J =(v_1 ,\ldots, v_{k -1})$ {\em if and only if} $L\subseteq H$ and $H\cap(J\backslash L) =\emptyset$. If $|L| =\ell$, then the number of such subsets $H$ of a given size $h <k$ is equal to $\binom{q -k}{h -\ell}$. 

We deduce that for each $h\in\set{0 ,\ldots, k -1}$, the construction generates 
$$\textstyle R^*\times\binom{q -h -2}{k -h -1}\times\binom{q -k}{h -\ell}$$
rows of the form $(\alpha(H), v_k)$ where $|H| =h$ and $\alpha(H)_J =(v_1 ,\ldots, v_{k -1})$. These rows are inserted into $\Psi$ if either $h$ has the same parity as $k -1$ and $v_k =q -1$, or $h$ does not have the same parity as $k -1$ and $v_k =0$; otherwise, they are inserted into $\Phi$. Hence, we have:
\begin{align}\label{eq-Nqkk-fin}
\begin{array}{l}
	\card{\set{r \in [R']\,|\,\Psi_r^J =v}} 
	-\card{\set{r \in [R']\,|\,\Phi_r^J =v}}\\
	\qquad=	R^*\times\left\{\begin{array}{rl}
		\sum_{h =\ell}^{k -1} (-1)^{k -1 -h} \binom{q -k}{h -\ell} \binom{q -2 -h}{k -1 -h} &\qquad\text{if $v_k =q -1$}\\
		-\sum_{h =\ell}^{k -1} (-1)^{k -1 -h} \binom{q -k}{h -\ell} \binom{q -2 -h}{k -1 -h}&\qquad\text{if $v_k =0$}
	\end{array}\right.
\end{array}
\end{align}

By the assumption that $(v_1 ,\ldots, v_{k -1})\neq J$, $\ell$ is some integer in $\set{0 ,\ldots, k -2}$. Given any such integer $\ell$, we observe:
$$\begin{array}{rl}
	\sum_{h =\ell}^{k -1} (-1)^{k -1 -h} \binom{q -k}{h -\ell} \binom{q -2 -h}{k -1 -h}
		&=\sum_{j =0}^{k -1 -\ell} (-1)^{k -1 -\ell -j} \binom{q -k}{j} \binom{q -2 -\ell -j}{k -1 -\ell -j}\\
		&=(-1)^{k -1 -\ell} \times S(q -k, q -2 -\ell, k -1 -\ell)
\end{array}$$

For $\ell\in\set{0 ,\ldots, k -2}$, we have $q -2 -\ell\geq q -k$, and $q -2 -\ell\geq k -1 -\ell$ (since $q >k$). We thus know from \cref{eq-S} that $S(q -k, q -2 -\ell, k -1 -\ell)$ sums to the binomial coefficient $\binom{k -2 -\ell}{k -1 -\ell}$, which is 0. We conclude that $\Psi$ and $\Phi$ do satisfy \cref{eq-qkk}, which means that $\mu^\Psi -\mu^\Phi$ is balanced $k$-wise independent: the proof is  complete. 
\end{proof}

\subsubsection{\Cref{cor-reduc-dapx}}

Our goal is to establish the inequality:
$$\begin{array}{rlll}
(T(a, b) +1)/2	&\leq (2a -b)^b/(2\times b!),	&a, b\in\mathbb{N},\ a >b	&\cref{eq-T-UB}
\end{array}$$
where given any two integers $b\geq 0$ and $a >b$, $T(a, b)$ is defined by:
$$\begin{array}{rlc}
	T(a, b) 		&=\sum_{r =0}^b \binom{a}{r} \binom{a -1 -r}{b -r}	&\cref{eq-T-def}
\end{array}$$

\begin{proof}
Applying first recursion \cref{eq-T_b+} to $T(a, b)$, and then recursion \cref{eq-T_a+b+} to $T(a, b -1)$, we get:
$$\begin{array}{rll}
	T(a, b) +1	
		&=2^b \binom{a}{b} -T(a, b -1)	+1									&\text{by \cref{eq-T_b+}}\\
		&=2^b \binom{a}{b} -2^{b -1} \binom{a -1}{b -1} -T(a -1, b -2) +1	&\text{by \cref{eq-T_a+b+}}
\end{array}$$

First, we deduce from recursion \cref{eq-T_a+b+} and the definition \cref{eq-T-def} of $T(a, b)$ that we have:
$$\begin{array}{rllll}
	T(a -1, b -2) -1 &\geq T(a -b +1, 0) -1 &=\binom{a -b +1}{0}\binom{a -b -0}{0} -1 &=0
\end{array}$$

So $T(a, b) +1 \leq 2^b\binom{a}{b} -2^{b -1}\binom{a -1}{b -1}$.
Now we rewrite the expression $2^b\binom{a}{b} -2^{b -1}\binom{a -1}{b -1}$ as:
$$\begin{array}{rl}
	2^b \binom{a}{b} -2^{b -1} \binom{a -1}{b -1}
		&= 		2^b\binom{a -1}{b -1}\times(a/b -1/2)\\
		&= 		(2^b/b!) \times (a -1)\times\ldots\times(a -b +1) \times (a -b/2)
\end{array}$$
Due to the inequality of arithmetic and geometric means, the following relation holds:
$$\begin{array}{rlll}
\prod_{i =1}^{b -1} (a -i)	&\leq \left(\sum_{i =1}^{b -1} (a -i)/(b -1)\right)^{b -1}	&=(a -b/2)^{b -1}
\end{array}$$

Putting it all together, we finally get:
$$\begin{array}{rll}
(T(a, b) +1)/2		&\leq 1/2\times	2^b/b! \times (a -b/2)^{b -1} \times (a -b/2)
					&=(2a -b)^b/(2\times b!)
\end{array}$$

Inequality \cref{eq-T-UB} thus holds.
\end{proof}

\subsubsection{\Cref{cor-reduc-2E_q}}
\begin{table}\footnotesize{\begin{center}
$\begin{array}{c|c}
\begin{array}{c} 
	\begin{array}{cc}
		\multicolumn{2}{c}{(\Psi,\Phi)\in\Gamma_E(3, 1, 3, 2, 2):}\\[3pt]
		\setlength{\arraycolsep}{2pt}\begin{array}{ccc}	 
			\Psi^1	&\Psi^2	&\Psi^3\\\hline
			0 		&1		&0\\
			0 		&0		&1\\
			0 		&2		&2\\
		\end{array}&\setlength{\arraycolsep}{2pt}\begin{array}{ccc}
			\Phi^1	&\Phi^2	&\Phi^3\\\hline
			0 		&1		&2\\
			0 		&2		&1\\
			0 		&0		&0\\
		\end{array}
	\end{array}\\\\[-2pt]\begin{array}{cc}
		\multicolumn{2}{c}{(\Psi,\Phi)\in\Gamma_E(4, 1, 4, 2, 2):}\\[3pt]
		\setlength{\arraycolsep}{2pt}\begin{array}{cccc}
			\Psi^1	&\Psi^2	&\Psi^3	&\Psi^4\\\hline
			0 		&1		&0		&1\\
			0 		&3		&0		&3\\
			0 		&0		&2		&2\\
			0 		&2		&2		&0\\
		\end{array}&\setlength{\arraycolsep}{2pt}\begin{array}{cccc}
			\Phi^1	&\Phi^2	&\Phi^3	&\Phi^4\\\hline
			0 		&1		&2		&3\\
			0 		&3		&2		&1\\
			0 		&2		&0		&2\\
			0 		&0		&0		&0\\
		\end{array}
	\end{array}\\\\[-2pt]\begin{array}{cc}
		\multicolumn{2}{c}{(\Psi,\Phi)\in\Gamma_E(10, 2, 5, 2, 2):}\\[3pt]
		\setlength{\arraycolsep}{2pt}\begin{array}{ccccc}
			\Psi^1	&\Psi^2	&\Psi^3	&\Psi^4	&\Psi^5\\\hline
			0		&1		&0		&1		&0\\
			0		&0		&1		&0		&1\\
			0		&4		&4		&0		&4\\
			0		&1		&0		&0		&1\\
			0		&4		&0		&4		&4\\
			0		&0		&2		&2		&0\\
			0		&0		&0		&2		&2\\
			0		&3		&3		&3		&0\\
			0		&0		&3		&3		&3\\
			0		&2		&2		&0		&0\\
		\end{array}&\setlength{\arraycolsep}{2pt}\begin{array}{ccccc}
			\Phi^1	&\Phi^2	&\Phi^3	&\Phi^4	&\Phi^5\\\hline
			0		&1		&2		&3		&4\\
			0		&1		&2		&3		&4\\
			0		&4		&3		&2		&1\\
			0		&4		&3		&2		&1\\
			0		&2		&4		&1		&3\\
			0		&3		&1		&4		&2\\
			0		&0		&0		&0		&0\\
			0		&0		&0		&0		&0\\
			0		&0		&0		&0		&0\\
			0		&0		&0		&0		&0\\
		\end{array}
	\end{array}
\end{array}& 
\begin{array}{cc}
	\multicolumn{2}{c}{(\Psi,\Phi)\in\Gamma_E(21, 3, 7, 2, 2):}\\[3pt]
	\setlength{\arraycolsep}{2pt}\begin{array}{ccccccc}
		\Psi^1	&\Psi^2	&\Psi^3	&\Psi^4	&\Psi^5	&\Phi^6	&\Phi^7\\\hline
		0		&1		&0		&1		&0		&1		&0	\\
		0		&0		&1		&0		&1		&0		&1	\\
		0		&6		&6		&0		&6		&0		&6	\\
		0		&1		&0		&0		&1		&0		&1	\\
		0		&6		&0		&6		&6		&0		&6	\\
		0		&1		&0		&1		&0		&0		&1	\\
		0		&6		&0		&6		&0		&6		&6	\\
		0		&0		&2		&2		&0		&0		&2	\\
		0		&5		&5		&0		&0		&5		&5	\\
		0		&2		&0		&0		&2		&2		&0	\\
		0		&0		&2		&0		&0		&2		&2	\\
		0		&5		&5		&0		&5		&5		&0	\\
		0		&0		&5		&5		&0		&5		&5	\\
		0		&2		&2		&0		&0		&2		&0	\\
		0		&0		&0		&3		&3		&3		&0	\\
		0		&0		&0		&0		&3		&3		&3	\\
		0		&4		&4		&4		&4		&0		&0	\\
		0		&0		&4		&4		&4		&4		&0	\\
		0		&0		&0		&4		&4		&4		&4	\\
		0		&3		&3		&3		&0		&0		&0	\\
		0		&0		&3		&3		&3		&0		&0	\\
	\end{array}&\setlength{\arraycolsep}{2pt}\begin{array}{ccccccc}
		\Phi^1	&\Phi^2	&\Phi^3	&\Phi^4	&\Phi^5	&\Phi^6	&\Phi^7\\\hline
		0		&1		&2		&3		&4		&5		&6	\\
		0		&1		&2		&3		&4		&5		&6	\\
		0		&1		&2		&3		&4		&5		&6	\\
		0		&6		&5		&4		&3		&2		&1	\\
		0		&6		&5		&4		&3		&2		&1	\\
		0		&6		&5		&4		&3		&2		&1	\\
		0		&2		&4		&6		&1		&3		&5	\\
		0		&2		&4		&6		&1		&3		&5	\\
		0		&5		&3		&1		&6		&4		&2	\\
		0		&5		&3		&1		&6		&4		&2	\\
		0		&3		&6		&2		&5		&1		&4	\\
		0		&4		&1		&5		&2		&6		&3	\\
		0		&0		&0		&0		&0		&0		&0	\\
		0		&0		&0		&0		&0		&0		&0	\\
		0		&0		&0		&0		&0		&0		&0	\\
		0		&0		&0		&0		&0		&0		&0	\\
		0		&0		&0		&0		&0		&0		&0	\\
		0		&0		&0		&0		&0		&0		&0	\\
		0		&0		&0		&0		&0		&0		&0	\\
		0		&0		&0		&0		&0		&0		&0	\\
		0		&0		&0		&0		&0		&0		&0	\\
	\end{array}
\end{array}
\end{array}$
\caption{\label{tab-gamma_E-q22}
Relaxed $(q, 2)$-ARPAs of strength 2 with a ratio of $R^*$ to $R$ equal to $1/q$, for $q\in\set{3, 4, 5, 7}$.}
\end{center}}\end{table}

\Cref{tab-gamma_E-q22} shows the relaxed ARPAs of $\Gamma_E(R^*,R,2,2)$ with $R^*/R =1/q$ for $q\in\set{3, 4, 5, 7)}$. For $q =8$, the linear program (in lp format) in continuous variables modeling $\gamma_E(8, 2, 2)$ (see \cref{sec-cd}) and the optimal solution given by CPLEX for this problem are available at the URL: 
\href{https://lipn.univ-paris13.fr/~toulouse/cd/RARPA822/}{https://lipn.univ-paris13.fr/~toulouse/cd/RARPA822/}.

\subsection{Array pairs of \cref{sec-vois}}
\label{sec-arpa2cpa}
\label{sec-vois-id}

\subsubsection{\Cref{prop-vois-gamma-delta}}

Consider a CPA $(\Psi,\Phi)\in\Gamma(R, R^*, n, d, k)$, for five positive integers $R$, $R^*$, $n$, $d$, and $k$ such that $n\geq d\geq k$ and $R \geq R^*$. Our goal is to show that $(\sigma_n(\Psi), \sigma_n(\Phi))\in\Delta(R, R^*, n, d, k)$, where $\sigma_n$ is the mapping defined by \cref{eq-delta-gamma-def}.

\begin{proof}
First, by definition of $\sigma_n$, a row $\sigma_n(\Phi)_r$ (where $r\in[R]$) consists entirely of 1's {\em if and only if} $\Phi_r$ is of the form $(0, 1,\ldots, n -1)$. 
 %
Second, the number of 1's in $\sigma_n(\Psi)_r$ corresponds to the number of indices $j\in\Sigma_n$ for which $\Psi_r^j =j$. Since each row $\Psi_r$ of $\Psi$ has at most $d$ different values, the equality  $\Psi_r^j =j$ cannot occur for more than $d$ values of $j\in\Sigma_n$. Thus, each row of $\sigma_n(\Psi)$ contains at most $d$ 1's, ensuring the second condition of
\cref{def-vois-CD}.
 %
Finally, consider a $k$-length subsequence $J =(j_1 ,\ldots, j_k)$ of $(1 ,\ldots, n)$, along with a vector $u\in\set{0, 1}^k$. We denote by $\mathcal{V}(u)$ the set of vectors $v\in\Sigma_n^k$ such that:
$$\begin{array}{ll}
	\left\{\begin{array}{rll}
		v_s	&=		j_s -1	&\text{if $u_s =1$}\\ 
		v_s	&\neq 	j_s -1	&\text{if $u_s =0$}\\ 
	\end{array}\right.,		&s\in[k]
\end{array}$$

Let $H =(j_1 -1 ,\ldots, j_k -1)$ and $M\in\set{\Psi, \Phi}$.
By definition of $\sigma_n$, the frequency of $u$ in the subarray $\sigma_n(M)^J$ corresponds to the total frequency of the vectors from $\mathcal{V}(u)$ in the subarray $M^H$. Since $\mu^\Psi -\mu^\Phi$ is balanced $k$-wise independent, each $v\in\mathcal{V}(u)$ occurs equally often as a row in both subarrays $\Psi^H$ and $\Phi^H$. The total frequency of the vectors from $\mathcal{V}(u)$ is therefore the same in both subarrays $\Psi^H$ and $\Phi^H$. Equivalently, $u$ occurs the same number of times as a row in both subarrays $\sigma_n(\Psi)^J$ and $\sigma_n(\Phi)^J$. 

Thus, we have shown that $(\sigma_n(\Psi), \sigma_n(\Phi))$ does indeed satisfy all three conditions of \cref{def-vois-CD}, with the multiplicity $R^*$ for the row of 1's. 
\end{proof}

\subsubsection{\Cref{thm-vois-id}}

Let $k\geq 1$ and $n >k$ be two integers, and $(\Psi,\Phi)$ be the pair of Boolean arrays obtained by applying the map $\sigma_n$ of \cref{prop-vois-gamma-delta} to the arrays of the ARPA returned by \cref{algo-Gamma} on input $(k, n)$. We argue that $(\Psi,\Phi)$ can be described as follows:
\begin{itemize}
	\item the word of 1's occurs exactly once as a row in $\Phi$;
	\item for all $a\in\set{0 ,\ldots, k}$ with $a\equiv k\bmod{2}$, every word $u\in\set{0, 1}^n$ containing $a$ 1's occurs exactly $\binom{n -1 -a}{k -a}$ times as a row in $\Psi$;
	\item for all $a\in\set{0 ,\ldots, k}$ with $a\not\equiv k\bmod{2}$, every word $u\in\set{0, 1}^n$ containing $a$ 1's occurs exactly $\binom{n -1 -a}{k -a}$ times as a row in $\Phi$;
	\item no other word of $\set{0, 1}^n$ occurs in either $\Psi$ or $\Phi$.
\end{itemize}

\begin{algorithm}
\begin{algorithmic}[1]
\State\label{Delta-init}$\Psi, \Phi\gets\set{\beta(k, [k])}$
\For{$i =k +1$ {\bf to} $n$}
	\State\label{Delta-it_init}Insert an $i$th column of zeros in $\Psi$ and $\Phi$
	\State\label{Delta-it_init_Phi}Set the $i$th coefficient of the first row of $\Phi$ to 1
	\ForAll{$J\subseteq[i -1]$ such that $|J| <k$}\label{Delta-it_J}
		 \State\label{Delta_it_J1}Insert $\binom{i -2 -|J|}{k -1 -|J|}$ rows of the form 
		 	$\beta(i, J\cup\set{i})$ in $\Psi$ if $|J|\not\equiv k\bmod{2}$, in $\Phi$ otherwise
		 \State\label{Delta_it_J0}insert $\binom{i -2 -|J|}{k -1 -|J|}$ rows of the form
		 	$\beta(i, J)$ in $\Psi$ if $|J|\equiv k\bmod{2}$, in $\Phi$ otherwise
	\EndFor
\EndFor
\State\Return $(\Psi, \Phi)$ 
\end{algorithmic}
\caption{Constructing a CPA of $\Delta\left((T(n, k) +1)/2, n, k, k \right)$ 
	given two positive integers $k$ and $n >k$. 
$\beta(k, [k])$, $\beta(i, J\cup\set{i})$ and $\beta(i, J)$ are the Boolean words defined by \cref{eq-def-beta}.}
\label{algo-Delta}
\end{algorithm}

\begin{proof}
For an integer $i\in\set{k ,\ldots, n -1}$ and a subset $J$ of $[i]$, we denote by $\beta(i, J)$ the incidence vector of $J$ seen as a subset of $[i]$. Namely, $\beta(i, J)$ is the word of $\set{0, 1}^i$ defined by:
\begin{align}\label{eq-def-beta}
	\beta(i, J)_j	&=\left\{\begin{array}{rl}
						1	&\text{if $j \in J$}\\
						0	&\text{otherwise}
					\end{array}\right.,		&j\in[i]
\end{align}

In particular, $\beta(k, [k])$ and $\beta(n, [n])$ are the words of $k$ and $n$ 1's, respectively.
Applying the transformation $\sigma_n$ to the arrays returned by \cref{algo-Gamma} on input $(k, n)$ reduces to running \cref{algo-Delta} on input $(k, n)$. \Cref{tab-Delta-rec} shows the construction when $(k, n)\in\set{(2, 6), (3, 5)}$. 

To establish identity \cref{eq-vois-id}, we count the number of occurrences of each word of $\set{0, 1}^n$ in the resulting arrays $\Psi$ and $\Phi$. 
In \cref{algo-Delta}, \cref{Delta-init} first inserts a single row of the form $\beta(k, [k])$ into each of the two arrays. \Cref{Delta-it_init,Delta-it_init_Phi} then expand these rows into the rows $\beta(n, [k])$ and $\beta(n, [n])$ in $\Psi$ and $\Phi$, respectively. 
During a given iteration $i\in\set{k +1 ,\ldots, n}$, \cref{Delta_it_J0,Delta_it_J1} generate rows of the form $\beta(i, J)$, where $J$ is a subset of cardinality at most $k$ of $[i]$ such that $|J| <k$ or $i\in J$; \cref{Delta-it_init} then extends each of these rows into the row $\beta(n, J)$.

So consider a subset $J$ of $[n]$ and the associated word $\beta(n, J)$. According to the previous observations, if $|J| =n$, then the word $\beta(n, J)$ occurs once, in $\Phi$. If $J =[k]$, then $\beta(n, J)$ occurs once, in $\Psi$. If $|J|\in\set{k +1 ,\ldots, n -1}$, then $\beta(n, J)$ occurs neither in $\Psi$ nor in $\Phi$.
So we assume that $|J|\leq k$ and $J\neq[k]$. We denote by $i^*$ the value 0 if $J =\emptyset$, the largest integer in $J$ otherwise. 
Rows of the form $\beta(n, J)$ in the final arrays can result either from inserting rows of the form $\beta(i^*, J)$ by \cref{Delta_it_J1} during iteration $i^*$ (if $i^* >k$), or from inserting rows of the form $\beta(i, J)$ by \cref{Delta_it_J0} during an iteration $i\in\set{\max\set{i^*, k} +1,\ldots, n}$ (if $i^* <n$). These rows all occur in $\Psi$ if $|J|\equiv k\bmod{2}$; otherwise they all occur in $\Phi$. 

Thus, on the one hand, rows of the form $\beta(n, J)$ all occur in the same array. 
On the other hand, the exact number of these rows in $\Psi$ or $\Phi$ is equal to:
$$\left\{\begin{array}{rl} 
\binom{n -|J| -1}{k -|J|}									&\text{if $i^* =n$}\\
							 \sum_{i =k +1}^n \binom{i -|J| -2}{k -|J| -1}
															&\text{if $i^*\in\set{0 ,\ldots, k}$}\\
\binom{i^* -|J| -1}{k -|J|}	+\sum_{i =i^* +1}^n \binom{i -|J| -2}{k -|J| -1}
															&\text{if $i^*\in\set{k +1 ,\ldots, n -1}$}
\end{array}\right.$$

Now we have trivially for each $t\in\set{k +1 ,\ldots, n -1}$:
$$\begin{array}{rlll} 
	 	\sum_{i =t}^n \binom{i -|J| -2}{k -|J| -1}
	&=	\sum_{i =t}^n \left(\binom{i -|J| -1}{k -|J|} -\binom{i -|J| -2}{k -|J|}\right)	
	&=	\binom{n -|J| -1}{k -|J|} -\binom{t -|J| -2}{k -|J|}
\end{array}$$

Thus, for $a\in\set{0 ,\ldots, k}$, each Boolean word of length $n$ with $a$ 1's is generated $\binom{n -a -1}{k -a}$ times, and occurs in the same array as the word of $n$ 1's {\em if and only if} $a\not\equiv k\bmod{2}$. This concludes the proof.
\end{proof}

\subsection{Instance families of \cref{sec-E-tight,sec-vois-tight}}
\label{sec-I^qk_n}
\label{sec-tilde{I}_n}

\subsubsection{Instances $I^{q, k}_n$ discussed in \cref{sec-E-tight,sec-vois-tight}}


First, we establish the following inequality, which will be used in further arguments:
\begin{align}\label{eq-binom}
\textstyle\binom{x +z}{k} +\binom{y}{k}	&\textstyle>\binom{y +z}{k} +\binom{x}{k},
	&k, x, y, z\in\mathbb{N},\ x >y,\ z >0,\ x +z\geq k
\end{align}
\begin{proof}
Equivalent to \cref{eq-binom}, $k$, $x$, $y$ and $z$ should satisfy 
$\binom{x +z}{k} -\binom{x}{k}	>\binom{y +z}{k} -\binom{y}{k}$. We write:
$$\begin{array}{rll}
	\binom{x +z}{k} -\binom{x}{k}	&=\sum_{i =0}^{z -1} \binom{x +i}{k -1}\\
	\binom{y +z}{k} -\binom{y}{k}	&=\sum_{i =0}^{z -1} \binom{y +i}{k -1}\\
	\binom{x +i}{k -1}				&\geq\binom{y +i}{k -1},	&i\in\set{0 ,\ldots, z -1}\\
\end{array}$$

By the assumption $x +z \geq k$, the above inequality is strict (at least) at rank $i =z -1$, implying \cref{eq-binom}. 
\end{proof}

Let $q\geq 2$, $k \geq 2$, and $n\geq k$ be three integers. 
Our goal is to establish that, on $I^{q, k}_n$, the worst solution value satisfies:
\begin{align}\label{eq-Iqnk-wor}
	\mathrm{wor}(I^{q, k}_n) &\textstyle =q\times\binom{n}{k}
\end{align}

\begin{proof}
Let $x_*$ denote a worst solution on $I^{q, k}_n$ and $n_0 ,\ldots, n_{q -1}$ denote the number of its coordinates equal to $0 ,\ldots, q -1$, respectively. 
If $n_0 =\ldots= n_{q -1} =n$, then $\sum_{a =0}^{q -1} \binom{n_a}{k} =q \times \binom{n}{k}$. 
Otherwise, there exist two indices $a, b\in\set{0 ,\ldots, q -1}$ such that $n_a >n$ and $n_b <n$. Taking \cref{eq-binom} at $(k, x, y, z) =(k, n, n_b, n_a -n)$, we get the inequality:
$$\begin{array}{rl}
\binom{n_a}{k} +\binom{n_b}{k}	&>\binom{n_a +n_b -n}{k} +\binom{n}{k}
\end{array}$$

We deduce that assigning the value $b$ to $n_a -n$ of the coordinates of $x_*$ equal to $a$ results in an objective value strictly lower than $v(I^{q, k}_n, x_*)$, which is a contradiction.
\end{proof}

Let $q \geq 2$, $k \geq 2$, $n \geq k$, and $d \leq n$ be four natural numbers. We denote by $x$ a best solution on $I^{q, k}_n$, among those that are at Hamming distance at most $d$ from a worst solution. Our goal is to establish that the objective value reached at $x$ has the following expression:
$$\begin{array}{rl}
	v(I^{q, k}_n, x) &=\binom{n +d}{k} +\binom{n -d}{k} +(q -2)\binom{n}{k}
\end{array}$$

\begin{proof}
We denote by $n_0 ,\ldots, n_{q -1}$ the number of coordinates of $x$ equal to $0 ,\ldots, q -1$. In follows from \cref{eq-Iqnk-wor} that the Hamming distance of $x$ to a worst solution is the quantity:
\begin{align}\label{eq-Bd_*}
\begin{array}{rll}
\sum_{a\in\Sigma_q: n_a >n} (n_a -n)
	&=\sum_{a\in\Sigma_q: n_a <n} (n -n_a)
	&=\sum_{a =0}^{q -1}\abs{n_a -n}/2
\end{array}
\end{align}	

Suppose there are two different integers $a, b$ such that $n_a, n_b <n$. Since the Hamming distance of $x$ from a worst solution is at most $d$, while $d\leq n$, we have $2n -n_a -n_b\leq d\leq n$ and thus, $n_b \geq n -n_a$. We can therefore consider the solution obtained from $x$ by assigning the value $a$ to $n -n_a$ of its coordinates equal to $b$. Taking \cref{eq-binom} at $(k, x, y, z) =(k, n_a, n_a +n_b -n, n -n_a)$, we obtain the inequality:
$$\begin{array}{rl}
\binom{n}{k} +\binom{n_a +n_b -n}{k}	&>\binom{n_b}{k} +\binom{n_a}{k}
\end{array}$$

So this solution, which is at the same distance from a worst solution as $x$, performs an objective value strictly better than $x$: contradiction. We deduce that $n_a <n$ cannot occur for more than one integer $a\in\Sigma_q$. By a similar argument, $n_b >n$ cannot occur for more than one integer $b\in\Sigma_q$. We deduce that $v(I, x)$ is $\binom{n +i}{k} +\binom{n -i}{k} +(q -2)\binom{n}{k}$ for some $i\in\set{0 ,\ldots, d}$, while this expression is maximized at $i =d$.
\end{proof}

\subsubsection{Instances $\tilde{I}_n$ discussed in \cref{sec-vois-tight}}

$\tilde{I}_n$ is the instance of $\mathsf{E2\,Lin\!-\!2}$ that considers $2n$ Boolean variables $x_1 ,\ldots, x_{2n}$ and the $m =\binom{2n}{2} -n$ constraints $(x_j =x_h)$ where $j, h\in[2n]$, $j <h$, and $(j, h)\notin\set{(2\ell -1, 2\ell)\,|\,\ell\in[n]}$. 
This instance is trivially satisfiable ({\em e.g.}, by the zero vector), thus $\mathrm{opt}(\tilde{I}_n) =m$. 
We now prove that the worst solution value on $\tilde{I}_n$ has the following expression:
\begin{align}\label{eq-Jn-wor}
\mathrm{wor}(\tilde{I}_n) 	&\textstyle=m -n^2 +(n\bmod{2}), &n\in\mathbb{N}\backslash\set{0}
\end{align}

\begin{proof}
Consider a solution $x$ of $\tilde{I}_n$ with $a\in\set{0 ,\ldots, 2n}$ non-zero coordinates. 
Without loss of generality, we assume that $a\leq n$ (otherwise, we consider $\bar x$ instead of $x$). The objective value taken at such a vector $x$ can be expressed as:
\begin{align}\nonumber
v(\tilde{I}_n, x)	
	&\textstyle=m -a(2n -a) +\card{\set{j\in[n]\,|\,x_{2j -1}\neq x_{2j}}}
\end{align}

For a fixed $a\in\set{0 ,\ldots, n}$, $\card{\set{j\in[n]\,|\,x_{2j -1}\neq x_{2j}}}$ is minimized, {\em e.g.} by setting: 
$$\begin{array}{rl}
	x_j	&= 	\left\{\begin{array}{rll}
				1		&\text{if $j\leq a$}\\
				0		&\text{otherwise}
			\end{array}\right.
\end{array}$$

For such a vector $x$, $x_{2j -1}\neq x_{2j}$ folds for the single index $j =(a -1)/2$ if $a$ is odd, otherwise for no index, and hence:
$$\begin{array}{rl}
	v(\tilde{I}_n, x)	&=m -a(2n -a) +(a \bmod{2})
\end{array}$$

Finally, we observe that $n$ is a maximizer of $a(2n -a) -(a\bmod{2})$ over $\set{0 ,\ldots, n}$.
\end{proof}

Clearly, on $\tilde{I}_n$, the solution whose non-zero coordinates are the odd-index coordinates is a local optimum with respect to $\tilde{B}^1$, with value $m -n^2 +n$. This solution realizes a differential ratio of:
$$\begin{array}{rll}
	\displaystyle\frac{m -n^2 +n -\left(m -n^2 +(n\bmod{2})\right)}
							 {m -\left(m -n^2 +(n\bmod{2})\right)}
		&\displaystyle=	\frac{n -(n\bmod 2)}{n^2 -(n\bmod 2)}
		&\displaystyle=	\frac{1}{n +(n\bmod 2)}	
\end{array}$$

\end{document}